\input ./style/arxiv-general.cfg
\documentclass[aop,MSNbibl,dvips]{arximspdf}
\makeatletter
   \@ifpackageloaded{graphicx}{}{\usepackage{graphicx}}
\makeatother

\doi{10.1214/15-AOP1005}
\volume{44}
\issue{2}
\pubyear{2016}
\firstpage{1488}
\lastpage{1534}
\docsubty{FLA}

\makeatletter
\newcommand{\eqref}[1]{(\ref{#1})}
\newtheorem{theorem}{Theorem}[section]
\newtheorem{corollary}[theorem]{Corollary}
\newproclaim{example}[theorem]{Example}
\newtheorem{proposition}[theorem]{Proposition}
\newtheorem{lemma}[theorem]{Lemma}
\newproclaim{definition}[theorem]{Definition}
\newproclaim{remark}[theorem]{Remark}

\def\cB{\mathcal{B}}
\def\cG{\mathcal{G}}
\def\cH{\mathcal{H}}
\def\cF{\mathcal{F}}
\def\cS{\mathcal{S}}
\def\cX{\mathcal{X}}
\def\cY{\mathcal{Y}}

\def\bD{\mathbb{D}}
\def\bN{\mathbb{N}}
\def\bR{\mathbb{R}}
\def\bZ{\mathbb{Z}}

\makeatother

\begin{document}
\begin{frontmatter}

\title{Intermittency for the wave and heat equations with fractional
noise in time}
\runtitle{Intermittency for wave and heat equations}

\begin{aug}
\author[A]{\fnms{Raluca M.}~\snm{Balan}\corref{}\ead[label=e1]{rbalan@uottawa.ca}\ead[label=u1,url]{http://aix1.uottawa.ca/\textasciitilde rbalan/}\thanksref{T1}}
\and
\author[B]{\fnms{Daniel}~\snm{Conus}\ead[label=e2]{daniel.conus@lehigh.edu}\ead[label=u2,url]{http://www.lehigh.edu/\textasciitilde dac311/}}
\runauthor{R.~M. Balan and D. Conus}
\thankstext{T1}{Supported by a grant from the Natural Sciences and
Engineering Research Council of Canada.}
\address[A]{Department of Mathematics and Statistics\\
University of Ottawa\\
585 King Edward Avenue\\
Ottawa, Ontario K1N 6N5\\
Canada\\
\printead{e1}\\
\printead{u1}}
\address[B]{Department of Mathematics\\
Lehigh University\\
14 East Packer Avenue\\
Bethlehem, Pennsylvania 18109\\
USA\\
\printead{e2}\\
\printead{u2}}
\affiliation{University of Ottawa and Lehigh University}
\end{aug}

%
\received{\smonth{10} \syear{2013}}
%
\revised{\smonth{6} \syear{2014}}

%
\begin{abstract}
In this article, we consider the stochastic wave and heat equations
driven by a Gaussian noise which is spatially homogeneous and behaves
in time like a fractional Brownian motion with Hurst index $H>1/2$. The
solutions of these equations are interpreted in the Skorohod sense.
Using Malliavin calculus techniques, we obtain an upper bound for the
moments of order $p \geq2$ of the solution. In the case of the wave
equation, we derive a Feynman--Kac-type formula for the second moment
of the solution, based on the points of a planar Poisson process. This
is an extension of the formula given by Dalang, Mueller and Tribe
[\emph{Trans. Amer. Math. Soc.} \textbf{360} (2008) 4681--4703], in the
case $H=1/2$, and allows us to obtain a lower bound for the second
moment of the solution.
These results suggest that the moments of the solution grow much faster
in the case of the fractional noise in time than in the case of the
white noise in time.
\end{abstract}

%
%
\end{frontmatter}
%
\section{Introduction}
\label{sec:Intro}

In this article, we consider the stochastic wave equation
%
\renewcommand{\theequation}{SWE}
\begin{equation}\label{wave}
\cases{ %
\displaystyle\frac{\partial^2 u}{\partial t^2}(t,x) =  \Delta u(t,x) +
u(t,x)\dot{W}(t,x),&\quad $\bigl(t>0,x \in\bR^d \bigr),$
\vspace*{2pt}\cr
u(0,x) =  u_0,
\vspace*{2pt}\cr
\displaystyle\frac{\partial u}{\partial t}(0,x)  =  v_0, }
\end{equation}
and the stochastic heat equation
%
\renewcommand{\theequation}{SHE}
\begin{equation}\label{heat}
\cases{ %
\displaystyle\frac{\partial u}{\partial t}(t,x)  =  \frac
{1}{2}
\Delta u(t,x) + u(t,x)\dot{W}(t,x),&\quad $\bigl(t>0,x \in\bR^d \bigr),$
\vspace*{2pt}\cr
u(0,x)  =  u_0,}
\end{equation}
where $\Delta$ stands for the Laplacian operator on $\bR^d$,
and $\dot{W}$ denotes the formal derivative of a Gaussian noise $W$
(whose rigorous definition is given below). The definition of the
solution to equations (\ref{wave}) and (\ref{heat}) is given in Section~\ref{sec:framework} below, using the Skorohod integral with respect to $W$.
Intuitively, the noise $\dot{W}$ is homogeneous in space (with spatial
covariance kernel $f$) and behaves in time like a fractional Brownian
motion (fBm) with Hurst index $H>1/2$. The initial conditions $u_0$ and
$v_0$ are nonnegative constants. In the case of the wave equation
\eqref
{wave}, we assume that $d \leq3$, while for the heat equation \eqref
{heat}, $d \geq1$ can be arbitrary.

There is a large amount of literature dedicated to the case $H=1/2$,
when the noise behaves in time like the Brownian motion. In this case,
we say that the noise is white in time.
We refer the reader to
the lecture notes \cite{spde_book} for an introduction to the subject,
as well as \cite{dalang-frangos98,dalang,millet-sanzsole99,sanzsole-sarra02,QS-SanzSole04a,Nualart-QS07,CD08,dalang-sanzsole08,FK13} for a sample
of relevant references. The case $H\neq1/2$ has to be treated by
different methods, since the noise is not a semi-martingale in time.
In recent years, there has been a growing interest in studying
equations with general Gaussian noise, and in particular equations
driven by a noise which behaves in time like a fBm with Hurst parameter
$H\neq1/2$; see
\cite
{hu01,hu_nualart,BT10,BT10-SPA,B12-Fourier,HLN12,HNS11,HHNT14,CHS14,chen14}.



In the present article, the noise is introduced by a zero-mean Gaussian process
$W=\{W(\varphi);\varphi\in\cH\}$ with covariance
%
\renewcommand{\theequation}{\arabic{equation}}
\setcounter{equation}{0}
\begin{equation}
\label{cov-W} E \bigl(W(\varphi)W(\psi) \bigr)=\langle\varphi,\psi
\rangle_{\cH}.
\end{equation}
Here $\cH$ is a Hilbert space defined as the completion of the space
$C_0^{\infty}(\bR_{+} \times\bR^d)$ of infinitely differentiable
functions with compact support on $\bR_{+} \times\bR^d$, with respect
to the inner product $\langle\cdot, \cdot\rangle_{\cH}$ defined by
%
\begin{equation}
\label{def-cov} \langle\varphi,\psi\rangle_{\cH}= \alpha_H
\int_{(\bR_{+}
\times\bR
^d)^2} \varphi(t,x)\psi(s,y)|t-s|^{2H-2}f(x-y)\,dt
\,dx \,ds \,dy,
\end{equation}
where $\alpha_H=H(2H-1)$.
We assume that $H \in(\frac{1}{2},1)$, and $f$ is the Fourier
transform in $\cS'(\bR^d)$ of a tempered measure $\mu$ on $\bR^d$,
where $\cS'(\bR^d)$ is the dual of the space $\cS(\bR^d)$ of rapidly
decreasing infinitely differentiable functions on $\bR^d$.

Using the fact that
\[
\int_{\bR^d} \int_{\bR^d}\varphi(x)
\psi(y)f(x-y)\,dx\,dy=\int_{\bR
^d}\cF \varphi(\xi)\overline{\cF\psi(
\xi)}\mu(d\xi) \qquad\forall\varphi ,\psi\in\cS \bigl(\bR^d \bigr),
\]
we arrive at the following alternative expression for the inner product
$\langle\cdot, \cdot\rangle_{\cH}$:
%
\begin{eqnarray}
\label{def-cov2} &&\langle\varphi, \psi\rangle_{\cH}
\nonumber
\\[-8pt]
\\[-8pt]
\nonumber
&&\qquad=\alpha_H
\int_{0}^{\infty} \int_{0}^{\infty}
\int_{\bR^d} |t-s|^{2H-2} \cF\varphi(t,\cdot) (\xi)
\overline{\cF\psi(s,\cdot) (\xi)}\mu(d\xi)\,dt\,ds,
\end{eqnarray}
where $\cF$ denotes the Fourier transform in the $x$-variable.

In the present article, we consider the following four cases:
\begin{longlist}[(iii)]
\item[(i)] $f(0)<\infty$ (i.e., $\mu$ is a finite measure);
\item[(ii)] $f(x)=|x|^{-\alpha}$ for some $0<\alpha<d$ [i.e., $\mu
(d\xi
)=c_{\alpha,d}|\xi|^{-(d-\alpha)}\,d\xi$];
\item[(iii)] $ f(x)=\prod_{j=1}^{d}|x_j|^{-\alpha_j}$ for
some $\alpha_j \in(0,1)$ [i.e., $\mu(d\xi)=c_{(\alpha_j)_j}\times\break
\prod_{j=1}^{d}|\xi_j|^{\alpha_j-1}\,d\xi$];
\item[(iv)] $d=1$ and $f=\delta_0$ (i.e., $\mu$ is the Lebesgue measure).
\end{longlist}
Here we denote by $|x|$ the Euclidean norm of $x \in\bR^d$.

Case (i) corresponds to a spatially smooth noise $\dot{W}$. In case
(ii), $f$ is called the \emph{Riesz kernel} with exponent $\alpha$. Case
(iii) with the parametrization $\alpha_j=2-2H_j$ for some $H_j \in
(\frac{1}{2},1)$ leads to a noise $\dot{W}$ which is called a \emph{fractional Brownian sheet} with indices $(H,H_1, \ldots,H_d)$. Finally,
case (iv) corresponds to a (rougher) noise $\dot{W}$ which is ``white
in space.'' This describes the spatial behavior of the noise in the
four cases. On the other hand, in time, the noise is smoother than the
white noise (the Brownian motion), since $H > 1/2$.
We note in passing that the results of the present article can be
extended to $H=1/2$, recovering results which are already known for
equations \eqref{heat} and \eqref{wave} with white noise in time. To
ease the exposition, we discuss only the case $H>1/2$.

The stochastic heat equation \eqref{heat} driven by space--time white
noise $\dot{W}$ arises in different contexts and has been studied by
many authors. This equation is the continuous form of the parabolic
Anderson model studied by Carmona and Molchanov in \cite
{carmona_molchanov}, and plays a major role in the study of the
KPZ equation in physics; see \cite{KPZ}. The connection between the
stochastic heat equation and the KPZ equation (via the Hopf--Cole
transformation) was known informally by physicists for quite some time;
see, for example, \cite{bertini_cancrini}. Recently, this connection
has been made rigorous by Hairer in \cite{hairer}, using the theory of
rough paths; see also \cite{BQS}. Equation \eqref{heat} with fractional
noise in time has been studied in \cite{hu01,hu_nualart,BT10}.
References \cite{caithamer05,QT07,B12} are dedicated to the wave
equation with fractional noise.

In this article, we consider the Malliavin calculus approach for
defining a solution to equations \eqref{wave} and \eqref{heat}, as in
\cite{B12}, respectively \cite{BT10}. In particular, we introduce the
following assumption, known as \emph{Dalang's condition}: 
\renewcommand{\theequation}{DC}
\begin{equation}\label{cond-mu}
\int_{\bR^d}\frac{\mu(d\xi)}{1+|\xi|^2}<\infty.  
\end{equation}
This condition is necessary and sufficient for the existence of the
solution to equations \eqref{heat} and \eqref{wave}, when the noise is
white in time; see \cite{dalang}. It is also sufficient for the
existence of the solution to these equations, when the noise is
fractional in time and has spatial covariance given by the Riesz
kernel; see \cite{BT10,B12}.
The necessity of (\ref{cond-mu}) in the case of \eqref{heat} has been proved in
\cite{BC13}.

Note that (\ref{cond-mu}) is satisfied in cases (i) and (iv). In cases
(ii) and (iii), it holds if and only if $a<2$, where $a$ is defined by
\eqref{def-a} below.

The purpose of this paper is to study intermittency properties for the
solutions to equations \eqref{heat} and \eqref{wave}. Intuitively,
a space--time random field is called \emph{physically intermittent} if it
develops very high peaks concentrated on small spatial islands, as time
becomes large. To give a formal mathematical definition of
intermittency for a random-field $u=\{u(t,x); t \geq0, x \in\bR^d\}$,
we consider the \emph{upper Lyapunov exponent}
%
\renewcommand{\theequation}{\arabic{equation}}
\setcounter{equation}{3}
\begin{equation}
\gamma(p) := \limsup_{t \rightarrow\infty} \frac{1}{t} \log
E\bigl|u(t,x)\bigr|^p
\end{equation}
for any $p \geq1$ [assuming that $\gamma(p)$ does not depend on $x$].
Traditionally, in the literature, the random-field $u$ is called \emph
{weakly intermittent} if
%
\begin{equation}
\label{def-intermit} \gamma(2) > 0\quad \mbox{and}\quad \gamma(p) < \infty\qquad\mbox{for all } p
\geq2.
\end{equation}
If $\gamma(1) = 0$, and $u(t,x) \geq0$, then weak intermittency
implies full intermittency. Recall that a random field $u$ is \emph
{fully intermittent} if $p \mapsto\gamma(p)/p$ is strictly increasing;
see \cite{carmona_molchanov}. Intuitively, full intermittency shows
that for $p > q$,
%
\begin{equation}
\label{mom-asymp} \limsup_{t \rightarrow\infty} \frac{\|u(t,x)\|_p}{\|u(t,x)\|_q} = \infty,
\end{equation}
where $\|\cdot\|_p$ denotes the norm in $L^p(\Omega)$. In other words,
asymptotically, the $p$th moment of $u(t,x)$ is significantly larger
than its $q$th moment. This suggests that the random variable $u(t,x)$
may take very large values with small 
(but significant) probabilities,
and therefore it develops high peaks, when $t$ is large. We refer to
\cite{bertini_cancrini}, Section~2.4, for a detailed explanation of this
phenomenon.

Intermittency for the spatially-discrete heat equation was studied in
\cite{carmona_molchanov}. In~\cite{FK09}, Foondun and Khoshnevisan
proved weak intermittency for the solution to equation \eqref{heat}
driven by space--time white noise, assuming that the initial condition
$u_0$ is bounded away from 0. Similar investigations have been carried
out in \cite{bertini_cancrini,CK12,Chen_Dalang}.
In the recent article \cite{CHSX}, Chen, Hu, Song and Xing have given
the exact asymptotics for the moments of the solution to equation
\eqref{heat} driven by a fractional noise in time, with spatial
covariance kernel given by cases (ii)--(iv) above.
Intermittency for the solution of the stochastic wave equation driven
by a Gaussian noise which is white in time was studied in \cite
{dalang-mueller09,CJKS}.

The fractional aspect of the noise in time leads to a different notion
of weak intermittency, which is obtained by a slight modification of
the Lyapunov exponent. More precisely, for $\rho>0$ and $p \geq1$, we
define the \emph{modified upper Lyapunov exponent} (of index $\rho$) by
%
\begin{equation}
\label{upper-Lyapunov} \gamma_{\rho}(p) := \limsup_{t \rightarrow\infty}
\frac
{1}{t^{\rho}} \log E\bigl|u(t,x)\bigr|^p.
\end{equation}
By analogy with \eqref{def-intermit}, we say that the random-field $u$
is \emph{weakly $\rho$-intermittent} if
\[
\gamma_{\rho}(2) > 0 \quad\mbox{and}\quad \gamma_{\rho}(p) < \infty \qquad
\mbox{for all } p \geq2.
\]
Also, we say that $u$ is \emph{fully $\rho$-intermittent} if $p \mapsto
\gamma_{\rho}(p)/p$ is strictly increasing. These definitions guarantee
that for a fully $\rho$-intermittent random-field $u$, relation \eqref
{mom-asymp} still holds, and so the intuitive (physical) notion of
intermittency remains valid. A similar argument as the one developed in
\cite{bertini_cancrini} still applies to explain the existence of the
high peaks and the islands. Moreover, it remains true that weak $\rho
$-intermittency of $u$ implies its full $\rho$-intermittency, provided
that $u(t,x) \geq0$ and $\gamma_{\rho}(1) = 0$. (This can be proved by
convexity arguments which do not depend on the exponent of $t$ used in
the definition of $\rho$-intermittency.)


This article is organized as follows. In Section~\ref{sec:main}, we
describe our main results and introduce the exponents $\rho$ for
equations \eqref{wave} and \eqref{heat}. 
Section~\ref{sec:framework} contains a review of some Malliavin
calculus techniques which are needed for the definition of the
solution. 
In Section~\ref{existence-section}, we prove the existence of the
solution to equation \eqref{wave} in \emph{any} spatial dimension $d
\geq1$, and we give an upper bound for its second moment. An upper
bound for its $p$th moment is given in Section~\ref{UB-section}. In
Section~\ref{FK-section}, we obtain a Feynman--Kac-type representation
for the second moment of the solution of \eqref{wave} with $d \leq3$,
based on the points of a planar Poisson process. This result is used in
Section~\ref{LB-section} to obtain a lower bound for the second moment
of the solution to \eqref{wave}. Section~\ref{sec:PAM} is dedicated to
the equation \eqref{heat}.
An elementary estimate is given in Appendix \ref{sec:app}. Appendix \ref{sec:appB} contains the proof of an inequality 
which is used in Section~\ref{existence-section}.

\section{Main results} \label{sec:main}

In this section, we discuss the two main results of this article.

The following exponents are used
for the weak $\rho$-intermittency of the solutions to equations \eqref
{wave}, respectively \eqref{heat}:
%
\begin{equation}
\rho_\mathrm{ w} = \frac{2H+2-a}{3-a}, \qquad\rho_\mathrm{ h} =
\frac
{4H-a}{2-a}, \label{def-rho}
\end{equation}
 where
%
\begin{equation}
\label{def-a} a =\cases{ %
0, &\quad $\mbox{in case (i),}$
\vspace*{2pt}\cr
\alpha,&\quad $\mbox{in case (ii),}$
\vspace*{2pt}\cr
\displaystyle\sum_{j=1}^{d}\alpha_j, &\quad
$\mbox{in case (iii),}$
\vspace*{2pt}\cr
1, & \quad$\mbox{in case (iv).}$ }
\end{equation}

We are now ready to state the first result about equation \eqref{wave}.
We refer to \eqref{def-sol} below for the definition of the solution.

\begin{theorem}
\label{wave-main}
Let $f$ be a kernel of cases \textup{(i)--(iv}). Let $\rho_\mathrm{ w}$ and $a$ be
the constants given by \eqref{def-rho}, respectively \eqref{def-a}.
Assume that condition \eqref{cond-mu} holds.
\begin{longlist}[(a)]
\item[(a)] For any $d \geq1$, equation \eqref{wave} has a solution $\{
u(t,x);t \geq0, x \in\bR^d\}$, given by relation \eqref
{Wiener-chaos-sol} below. If $d \leq2$, the solution is unique.

\item[(b)] For any $d \geq1$, $p \geq2$, $x \in\bR^d$ and for any $t>0$,
\begin{equation}
\label{wave-Lp-bound} E\bigl|u(t,x)\bigr|^p \leq c_1^p(u_0+tv_0)^{p}
\exp \bigl(c_2 p^{(4-a)/(3-a)} t^{\rho_\mathrm{ w}} \bigr),
\end{equation}
where $c_1>0$ is a constant depending on $a$, and $c_2>0$ is a constant
depending on $H$ and $a$.

\item[(c)] Suppose that $d \leq3$. Then for any $x \in\bR^d$ and for any $t>0$,
%
\begin{equation}
E\bigl|u(t,x)\bigr|^2 \geq c_3 u_0^{2}
\exp \bigl(c_4 t^{\rho_\mathrm{ w}} \bigr),
\end{equation}
where $c_3>0$ and $c_4>0$ are constants depending on $H$ and $a$. 
\end{longlist}
\end{theorem}

A similar result holds for the parabolic equation \eqref{heat}.

\begin{theorem}
\label{heat-main}
Let $f$ be a kernel of cases \textup{(i)--(iv)}. Let $\rho_\mathrm{ h}$ and $a$ be
the constants given by \eqref{def-rho}, respectively \eqref{def-a}.
Assume that condition \eqref{cond-mu} holds. Let $d \geq1$ be arbitrary.
\begin{longlist}[(a)]
\item[(a)] Equation \eqref{heat} has a unique solution $\{u(t,x); t \geq0,x
\in\bR^d\}$.

\item[(b)] For any $p \geq2$, for any $x \in\bR^d$ and for any $t>0$, 
\begin{equation}
\label{heat-Lp-bound} E\bigl|u(t,x)\bigr|^p \leq c_1^p
u_0^{p} \exp \bigl(c_2 p^{(4-a)/(2-a)}
t^{\rho
_\mathrm{ h}} \bigr),
\end{equation}
where $c_1>0$ is a constant depending on $a$, and $c_2>0$ is a constant
depending on $H$ and $a$. 

\item[(c)] For any $x \in\bR^d$ and for any $t>0$,
%
\begin{equation}
E\bigl|u(t,x)\bigr|^2 \geq c_3 u_0^{2}
\exp \bigl(c_4 t^{\rho_\mathrm{ h}} \bigr),
\end{equation}
where $c_3>0$ and $c_4>0$ are constants depending on $H$ and $a$.
\end{longlist}
\end{theorem}

Most moment estimates for solutions to s.p.d.e.'s with white noise in
time rely on martingale properties of stochastic integrals. Since the
fBm is not a semi-martingale, different techniques have to be used when
the noise is fractional in time. In the case of equations \eqref{wave}
and \eqref{heat}, one can give explicitly the Wiener chaos
representation of the solution. The upper bounds \eqref{wave-Lp-bound}
and \eqref{heat-Lp-bound} are obtained directly using the equivalence
of $L^2(\Omega)$- and $L^p(\Omega)$-norms on each Wiener chaos.
The lower bounds require more work. For this, we follow the approach of
Dalang and Mueller \cite{dalang-mueller09}, which consists of using a
Feynman--Kac (FK) type representation for the second moment of the
solution, based on a Poisson process. Such a representation was
originally developed in \cite{DMT08} for equations driven by a noise
that is white in time. It was extended to the heat equation driven by
fractional noise in time by the first author of this article in \cite
{B09}. The extension to the wave equation with fractional noise in time
is given in Section~\ref{FK-section} below.

Article \cite{DMT08} contains also a FK representation for the $n$th
moment of the solution of the wave (or heat) equation, for any integer
$n \geq2$ (Theorem~5.1 of \cite{DMT08}). The proof of this result uses
the fact that the stochastic integral with respect to the noise $W$ is
a martingale in time, which allows the authors of \cite{DMT08} to apply
It\^o's formula. In the case of the fractional noise in time, the
stochastic integral is not a semi-martingale. There exists an It\^o's
formula for the Skorohod integral with respect to the classical fBm
(Theorem~8 of \cite{alos-nualart03}), which could probably be
generalized to the case of the noise $W$. However, this formula
contains an extra correction term involving the Malliavin derivative of
the integrand process, which is difficult to handle. For this reason,
we could not apply the method of Dalang, Mueller and Tribe \cite{DMT08}
to obtain an FK representation (and an exponential lower bound) for the
moment of order $n \geq2$ of the solution to either wave of heat
equation. We note that a lower bound for the $n$th moment of the
solution to the heat equation has been recently obtained in \cite
{HHNT14}, using an FK representation for the moments which is specific
to the parabolic case (see Theorem~3.6 of \cite{HHNT14}), and is
different than the one used in the present paper. The lower bound for
the $n$th moment of the solution to the wave equation remains an open problem.

As in \cite{dalang-mueller09}, we focus mainly on the hyperbolic case
(Theorem~\ref{wave-main}). The proof of Theorem~\ref{heat-main} is very
similar, and we only point out the differences in comparison to the
hyperbolic case in Section~\ref{sec:PAM}. We made this choice since the
results for the wave equation are completely new, in particular the
second moment FK type representation.

In \cite{CHSX}, Chen, Hu, Song and Xing obtained stronger results than
our Theorem~\ref{heat-main}, by computing the exact Lyapunov exponent
for the solution of
equation \eqref{heat}, defined as the limit when $t \to\infty$,
instead of the $\limsup$ in \eqref{upper-Lyapunov}; see Theorem~6.1 of
\cite{CHSX}.
In \cite{CHSX}, the solution is defined in the weak sense (i.e., using
multiplication against test functions), and the stochastic integral is
interpreted in the Stratonivich sense, according to Definition~4.2 of
\cite{HNS11}.
However, their method requires the additional assumptions
$a<4H-2$ in cases (ii)--(iii), and $H>3/4$ in case (iv), which are not
needed in the present article. The proofs of \cite{CHSX} rely on a
Feynman--Kac representation for the weak solution and its moments (due
to \cite{HNS11}), which can only be proved under the above-mentioned
additional assumptions.
By Theorem~7.2 of \cite{HNS11}, a similar Feynman--Kac representation
exists for the mild solution (defined using the Skorohod integral, as
in the present work), under the same assumptions mentioned above. Using
this representation and under the same assumptions, it may be possible
to compute the exact Lyapunov exponent for the mild solution, although
this is not proved in \cite{CHSX}. We believe that in the absence of
these assumptions, the methods of \cite{HNS11} and \cite{CHSX} cannot
be applied for equation \eqref{heat}, even when it is interpreted in
the Shorohod sense. These assumptions appear also in the recent
preprint \cite{HHNT14} for the Feynman--Kac representation for the
solution of equation \eqref{heat} interpreted in the Stratonovich sense
(see Hypothesis~4.1 of \cite{HHNT14}), but are not needed for obtaining
exponential upper and lower bounds for the moments of the solution of
\eqref{heat}, interpreted in the Stratonovich or Skorohod sense, as
shown by Theorem~6.4 of \cite{HHNT14}. In the case of the heat equation
with noise as in case (iii) above, some exponential upper and lower
bounds for the first moment of the solution (interpreted in the
Stratonovich sense) have been obtained in \cite{song12}.


The appropriate exponents $\rho$ are different in the hyperbolic and
parabolic cases. Nevertheless, since $H>1/2$, $\rho_\mathrm{ h}>\rho_\mathrm{
w}>1$. Therefore, the lower bounds in Theorems \ref{wave-main} and
\ref
{heat-main} imply that $\gamma(2) = \infty$, which shows that the
solutions to \eqref{wave} and \eqref{heat} are not weakly intermittent
in the classical sense. However, these solutions are weakly $\rho
$-intermittent (in the sense defined in Section~\ref{sec:Intro}) with
$\rho= \rho_\mathrm{ w}$ for the wave equation and $\rho= \rho_\mathrm{ h}$
for the heat equation. The results of Theorems~\ref{wave-main} and
\ref
{heat-main} do not provide full $\rho$-intermittency. When the noise is
white in time, one typically obtains full $\rho$-intermittency by
proving that $\gamma(1) = 0$. According to Song~\cite{song12} (using
the Stratonovitch integral), it appears that this may not be true in
the case of the fractional noise in time. In this case, an alternative
method is to obtain a sharp lower bound on the moments of order $p > 2$
of the solution, as in \cite{dalang-mueller09}. This is subject of
ongoing research.

In the case when $H = 1/2$, $\rho_\mathrm{ w} = \rho_\mathrm{ h} = 1$, and we
recover some of the known results of intermittency for the heat and
wave equations with white noise in time. For instance, intermittency
for the heat equation was studied in \cite{FK09}. For the wave
equation, full intermittency was obtained in \cite{dalang-mueller09}
with the spatial covariance of case~(i). Some upper bounds were
obtained in \cite{CJKS}.

As mentioned before, when $H > 1/2$, $\rho_\mathrm{ h}>\rho_\mathrm{ w}>1$. A
consequence of this is that the moments of the solution at some fixed
time $t$ are typically larger in the case of the fractional noise in
time compared to the white-noise case. This would imply that the size
of the peaks would be larger in the fractional case. Since $H > 1/2$,
the noise is positively correlated in time, which explains why peaks
build up larger values. Indeed, the fractional noise, when large, tends
to remain large for a longer period of time, which then results in a
higher build-up for the random-field $u$.

The upper and lower bounds given by Theorems \ref{wave-main} and \ref
{heat-main} show that the exponents $\rho_\mathrm{ w}$ and $\rho_\mathrm{ h}$
are sharp. A lower bound result on the moments of order $p > 2$ would
be needed in order to get the sharp behavior of the exponent $\gamma
_{\rho}(p)$ as a function of $p$. This remains an open problem in the
case of the wave equation. (In the case of the heat equation, this has
been recently proved in the preprint \cite{HHNT14}.) We note that in
our results, the behavior of $\gamma_{\rho}(p)$ as a function of $p$
does not depend on $H$. For the wave equation, we obtain that $\gamma
_{\rho_\mathrm{ w}}(p) \leq C p^{4/3}$ in case (i) (spatially smooth
noise), and $\gamma_{\rho_\mathrm{ w}}(p) \leq C p^{3/2}$ in case (iv)
(spatial white-noise), where $C>0$ denotes a constant which does not
depend on $p$. These confirm the behavior in the order $p$ obtained in
\cite{dalang-mueller09} for case (i) and in \cite{CJKS} for case (iv).
For the heat equation, we obtain that $\gamma_{\rho_\mathrm{ h}}(p) \leq C
p^{2}$ in case (i) and $\gamma_{\rho_\mathrm{ h}}(p) \leq C p^{3}$ in case
(iv), for a constant $C>0$ which does not depend on $p$, which
correspond to the sharp order for white noise in time.

Finally, we would like to point out that Theorems \ref{wave-main} and
\ref{heat-main} constitute a first step toward a more careful study of
the intermittent behavior of the solution to the stochastic heat and
wave equations. Indeed, following the program developed in \cite{CJK,CJKS2}, sharp Lyapunov exponents for the moments of solutions to SPDEs
(in particular their behavior as a function of $p$) are key ingredients
for obtaining quantitative results regarding some physical properties
of the solution, such as the height of the peaks, the size of the
peak-islands and some space--time scaling results for the behavior of
the peaks. These could lead to a careful understanding of the impact of
the temporal (and spatial) correlation of the noise on the physical
behavior of the solution. In particular, observing how the modified
Lyapunov exponents impact the physical properties would be an important
step in the understanding of mathematical intermittency. In the case of
the white-noise in time, the existence of the sharp Lyapunov exponents
has allowed the authors of \cite{CJK,CJKS2} to obtain KPZ-type scaling
exponents for the solution to the stochastic heat equation.

\section{Framework}
\label{sec:framework}

In this section, we introduce the framework, and we give a brief
summary of the results of Balan \cite{B12}, which are needed in the
present article.

We denote by $G_\mathrm{ w}$, respectively $G_\mathrm{ h}$, the fundamental
solution of the wave equation, respectively the heat equation.
In the case of the wave equation, recall that when $d \leq2$, $G_\mathrm{
w}(t,\cdot)$ is a function given by
%
\begin{eqnarray}
\label{def-Gw} G_\mathrm{ w}(t,x)&=&\frac{1}{2}1_{\{|x| \leq t\}}
\qquad\mbox{if } d=1 \quad\mbox{and}
\nonumber
\\[-8pt]
\\[-8pt]
\nonumber
 G_\mathrm{ w}(t,x)&=&\frac{1}{2\pi}
\frac{1}{\sqrt{t^2-|x|^2}}1_{\{|x| <
t\}
} \qquad\mbox{if } d=2.
\end{eqnarray}
In both cases, $\int_{\bR^d}G_\mathrm{ w}(t,x)\,dx=t$.
When $d=3$, $G_\mathrm{ w}(t,\cdot)$ is a finite measure on $\bR^3$ given
by
\[
G_\mathrm{ w}(t,\cdot)=\frac{1}{4\pi t}\sigma_t,
\]
where $\sigma_t$ is the surface measure on $\partial B(0,t)$, and
$G_\mathrm{ w}(t,\bR^3)=t$. When $d \geq4$, $G_\mathrm{ w}(t,\cdot)$ is a
distribution. For any $d \geq1$, the Fourier transform of $G_\mathrm{
w}(t,\cdot)$ is given by
%
\begin{equation}
\label{Fourier-G} \cF G_\mathrm{ w}(t,\cdot) (\xi)=\frac{\sin(t|\xi|)}{|\xi|},\qquad \xi
\in \bR^d.
\end{equation}

In the case of the heat equation, for any dimension $d \geq1$, $G_\mathrm{
h}(t, \cdot)$ is a function, known as the \emph{heat kernel}. More precisely,
\[
G_\mathrm{ h}(t,x)=\frac{1}{(2\pi t)^{d/2}}\exp \biggl(-\frac
{|x|^2}{2t}
\biggr) \quad\mbox{and}\quad \cF G_\mathrm{ h}(t,\cdot) (\xi)=\exp \biggl(-
\frac{t|\xi|^2}{2} \biggr).
\]
Below, we write $G$ when the results apply for both $G_\mathrm{ w}$ or
$G_\mathrm{ h}$.

We denote by $w_\mathrm{ w}$ (resp., $w_\mathrm{ h}$) the solution of the
homogeneous wave (resp., heat) equation with the same initial condition
as \eqref{wave} [resp., \eqref{heat}], that is,
\begin{equation}
\label{def-w} w_\mathrm{ w}(t,x) = u_0+t
v_0 \quad\mbox{and}\quad w_\mathrm{ h}(t,x) = u_0.
\end{equation}
We write $w$ when the results apply for both $w_\mathrm{ w}$ and $w_\mathrm{
h}$. Note that $w(t,x)$ does not depend on $x$ in either case,
and $w(t,x) \geq u_0 \geq0$ for all $t>0$ and $x \in\bR^d$.

We now discuss the concept of solution. Informally, a (mild) solution
of \eqref{wave} or \eqref{heat} should be a process $\{u(t,x);t \geq
0,x \in\bR^d\}$ which satisfies
%
\begin{equation}
\label{def-sol-1} u(t,x) = w(t,x) + \int_{0}^{t}
\int_{\bR^d}G(t-s,x-y)u(s,y)W(ds,dy),
\end{equation}
provided the stochastic integral on the right-hand side is well defined
(in a certain sense).
Still informally, replacing $u(s,y)$ on the right-hand side of (\ref
{def-sol-1}) by its definition and iterating this procedure, we
conclude that the solution of \eqref{wave} or \eqref{heat} should be
given by the following series of iterated integrals:
%
\begin{eqnarray}
\label{informal-series} u(t,x)&=&w(t,x)+\int_0^t \int
_{\bR^d}G(t-s,x-y)W(ds,dy)
\nonumber
\\
&&{}+
\int_0^t \int_{\bR^d} \int
_0^s \int_{\bR
^d}G(t-s,x-y)\\
&&\hspace*{77pt}{}\times G(s-r,y-z)W(dr,dz)W(ds,dy)+
\cdots.\nonumber
\end{eqnarray}

To give a rigorous meaning to this procedure, we use an approach based
on Malliavin calculus with respect to the isonormal Gaussian process
$W=\{W(\varphi);\varphi\in\cH\}$ with covariance specified by
\eqref
{cov-W}, where $\cH$ is the Hilbert space defined as the completion of
$C_{0}^{\infty}(\bR_{+} \times\bR^d)$ with respect to the inner
product $\langle\cdot, \cdot\rangle_{\cH}$ given by \eqref{def-cov}.

We recall the basic elements of Malliavin calculus; see \cite
{nualart06} for more details. It is known that every square-integrable
random variable $F$ which is measurable with respect to $W$, has the
Wiener chaos expansion
\[
F=E(F)+\sum_{n \geq1}F_n \qquad\mbox{with }
F_n \in\cH_n,
\]
where $\cH_n$ is the $n$th Wiener chaos space associated to $W$.
Moreover, each $F_n$ can be represented as $F_n=I_n(f_n)$ for some $f_n
\in\cH^{\otimes n}$, where $\cH^{\otimes n}$ is the $n$th tensor
product of $\cH$, and
$I_n\dvtx \cH^{\otimes n} \to\cH_n$ is the multiple Wiener integral with
respect to~$W$.
By the orthogonality of the Wiener chaos spaces and an isometry-type
property of $I_n$, we obtain that
\[
E|F|^2=(EF)^2+\sum_{n \geq1}E\bigl|I_n(f_n)\bigr|^2=(EF)^2+
\sum_{n \geq1}n! \| \widetilde{f}_n
\|_{\cH^{\otimes n}}^{2},
\]
where $\widetilde{f}_{n}$ is the symmetrization of $f_n$ in all $n$ variables
\[
\widetilde{f}_n(t_1,x_1,
\ldots,t_n,x_n)=\frac{1}{n!}\sum
_{\rho\in
S_n}f_n(t_{\rho(1)},x_{\rho(1)},
\ldots,t_{\rho(n)},x_{\rho(n)}),
\]
where $S_n$ is the set of all permutations of $\{1, \ldots,n\}$.

Let $\cS$ be the class of
smooth random variables of the form
%
\begin{equation}
\label{form-F}F=f \bigl(W(\varphi_1),\ldots, W(
\varphi_n) \bigr),
\end{equation}
where $f \in C_{b}^{\infty}(\bR^n)$, $\varphi_i \in\cH$, $n \geq
1$ and
$C_b^{\infty}(\bR^n)$ is the class of bounded $C^{\infty}$-functions
on $\bR^n$, whose partial derivatives are bounded. The \textit{Malliavin
derivative} of $F$ of the form (\ref{form-F}) is an $\cH$-valued
random variable given by
\[
DF:=\sum_{i=1}^{n}\frac{\partial f}{\partial x_i}
\bigl(W(\varphi _1),\ldots, W(\varphi_n) \bigr)
\varphi_i.
\]
%

We endow $\cS$ with the norm $\|F\|_{\bD^{1,2}}^{2}:=E|F|^2+E\|D F
\|_{\cH}^{2}$. The operator $D$ can be extended to the space
$\bD^{1,2}$, the completion of $\cS$ with respect to $\|\cdot
\|_{\bD^{1,2}}$.

The \textit{divergence operator} $\delta$ is defined as the adjoint of the
operator $D$. The domain of $\delta$, denoted by $\operatorname{Dom}
\delta$, is the set of $u \in L^2(\Omega;\cH)$ such that
\[
\bigl|E \langle DF,u \rangle_{\cH}\bigr| \leq c \bigl(E|F|^2
\bigr)^{1/2}\qquad \forall F \in\bD^{1,2},
\]
where $c$ is a constant depending on $u$. If $u \in\operatorname{ Dom}
\delta$, then $\delta(u)$ is the element of $L^2(\Omega)$
characterized by the following duality relation:
%
\begin{equation}
\label{duality} E \bigl(F \delta(u) \bigr)=E\langle DF,u \rangle_{\cH}\qquad
\forall F \in\bD^{1,2}.
\end{equation}
In particular, $E[\delta(u)]=0$. If $u \in\operatorname{Dom}  \delta$, we use
the notation
\[
\delta(u)=\int_0^{\infty} \int
_{\bR^d}u(t,x) W(\delta t, \delta x),
\]
and we say that
$\delta(u)$ is the \textit{Skorohod integral} of $u$ with respect to $W$.

We recall the following criterion for Skorohod integrability; see also
Proposition~1.3.7 of \cite{nualart06}. 

\begin{proposition}[(Proposition~2.5 of \cite{B12})]
\label{integr-crit}
Assume that $u \in L^2(\Omega;\cH)$ has the Wiener chaos expansion
%
\begin{equation}
\label{chaos-exp} u(t,x)=\sum_{n \geq0}I_n
\bigl(f_n(\cdot,t,x) \bigr),
\end{equation}
where $f_0(t,x)=E(u(t,x))$, $I_0(x)=x$ and $f_n(\cdot,t,x) \in\cH
^{\otimes n}$ for any $n \geq1$. 
Then $u \in\operatorname{ Dom} \delta$ if and only if the series $\sum_{n
\geq
0}I_{n+1}(f_n)$ converges in $L^2(\Omega)$, and in this case
$\delta(u)=\sum_{n \geq0}I_{n+1}(f_{n})$.
\end{proposition}

We are now ready to give the rigorous definition of the solution to
equations \eqref{heat} and \eqref{wave}. Let $\cF_t$ be the $\sigma
$-field generated by $W(1_{[0,s] \times A})$ for $s \in[0,t],A \in\cB
_b(\bR^d)$, where $\cB_b(\bR^d)$ is the class of all bounded Borel sets
in~$\bR^d$.

\begin{definition}
\label{definition-solution}
An $(\cF_t)_t$-adapted square-integrable process $u=\{u(t,x);\break t
\geq0, x \in\bR^d\}$ is called a \textit{\textup{(}mild\textup{)} solution} of \eqref
{wave} or \eqref{heat} if it satisfies the following integral equation:
%
\begin{equation}
\label{def-sol} u(t,x) = w(t,x) + \int_{0}^{t}
\int_{\bR^d}G(t-s,x-y)u(s,y)W(\delta s,\delta y),
\end{equation}
that is, $v^{(t,x)} \in\operatorname{ Dom}  \delta$ and $u(t,x)=w(t,x)+\delta
(v^{(t,x)})$ for all $(t,x) \in\bR_{+} \times\bR^d$, where
%
\begin{equation}
\label{def-v} v^{(t,x)}(s,\cdot)=1_{[0,t]}(s)G(t-s,x-\cdot)u(s,
\cdot),\qquad s \geq0
\end{equation}
and $\cdot$ denotes the missing $y$-variable.
\end{definition}

In the case of equation \eqref{wave} in dimension $d \leq2$ or
equation \eqref{heat} in any dimension $d$, $G(t,\cdot)$ is a function,
and the existence and uniqueness of the solution can be proved
similar to page 303 of \cite{hu_nualart}. We recall this argument
here. Assume that a solution $u(t,x)$ exists and has the Wiener chaos
expansion \eqref{chaos-exp} for some functions $f_n(\cdot,t,x) \in
\cH
^{\otimes n}$. Since $G$ is a deterministic function, it follows that
the process $v^{(t,x)}$ given by \eqref{def-v} has the Wiener chaos expansion
$v^{(t,x)}(s,y)=\sum_{n \geq0}I_{n}(g_n^{(t,x)}(\cdot,s,y))$, with kernels
%
\begin{equation}
\label{def-gn} g_{n}^{(t,x)}(\cdot,s,y)=1_{[0,t]}(s)G(t-s,x-y)f_n(
\cdot,s,y).
\end{equation}
By Proposition~\ref{integr-crit}, $v^{(t,x)} \in\operatorname{ Dom}  \delta$ if
and only if $\sum_{n \geq0}I_{n+1}(g_n^{(t,x)})$ converges in
$L^2(\Omega)$. In this case,
$\delta(v^{(t,x)})=\sum_{n \geq0}I_{n+1}(g_{n}^{(t,x)})$,
and relation $u(t,x)=w(t,x)+\delta(v^{(t,x)})$ becomes
\[
\sum_{n \geq0}I_n \bigl(f_n(
\cdot,t,x) \bigr)=w(t,x)+\sum_{n \geq
0}I_{n+1}
\bigl(g_n^{(t,x)} \bigr).
\]
By the uniqueness of the Wiener chaos expansion, we infer that
$f_0(t,x)=w(t,x)$ and $f_{n+1}(\cdot,t,x)=g_n^{(t,x)}$ for any $n \geq
0$. This allows us to find $f_n$ recursively:
%
\begin{eqnarray}\label{def-fn}
&&
f_n(t_1,x_1,
\ldots,t_n,x_n,t,x) \nonumber\\
&&\qquad =  G(t-t_n,x-x_n)G(t_n-t_{n-1},x_n-x_{n-1})
\cdots
\\
 & & \qquad\quad{}\times G(t_2-t_1,x_2-x_1)w(t_1,x_1)
1_{\{0<t_1<\cdots<t_n<t\}}.\nonumber
\end{eqnarray}
Therefore, if the series
$\sum_{n \geq0}I_{n+1}(g_n^{(t,x)})=\sum_{n \geq
0}I_{n+1}(f_{n+1}(\cdot,t,x))$ converges in $L^2(\Omega)$,
then the solution $u$ exists and is unique, with the Wiener chaos
expansion~\eqref{chaos-exp} with kernels $f_n(\cdot,t,x)$ given by
\eqref{def-fn}.
This coincides with the informal interpretation (\ref{informal-series}).

In the case of equation \eqref{wave} with $d \geq3$, the procedure for
constructing a solution is more complicated, since $G_\mathrm{ w}(t,\cdot)$
is a distribution in $\bR^d$. We describe below the steps of this
procedure, following \cite{B12}.

\emph{Step} 1. Define the kernel $f_n(\cdot,t,x)$ as a distribution in
$\cS'(\bR^{nd})$, identifying its action on a test function, as in
Section~2.1 of \cite{B12}. By Proposition~2.1 of \cite{B12}, for any
$0<t_1<\cdots<t_n<t$, $f_n(t_1,\cdot, \ldots,t_n,\cdot,t,x)$ is a
distribution in $\bR^{nd}$ whose Fourier transform [in $\cS'(\bR
^{nd})$] is the function
\begin{eqnarray}\label{Fourier-f}
&&\cF f_n(t_1,\cdot,\ldots,t_n,\cdot,t,x) (
\xi_1,\ldots,\xi_n)\nonumber\\
&&\qquad=(u_0+t_1
v_0)e^{-i (\xi_1+\cdots+\xi_n)\cdot x} \overline{\cF G_\mathrm{
w}(t_2-t_1,\cdot) (\xi_1)}
\\
 &&\qquad\quad{}\times \overline{\cF G_\mathrm{ w}(t_3-t_2,
\cdot) (\xi_1+\xi_2)} \cdots \overline {\cF
G_\mathrm{ w}(t-t_n,\cdot) (\xi_1+\cdots+
\xi_n)}.\nonumber
\end{eqnarray}
$f_n(t_1,\cdot,\ldots,t_n,\cdot,t,x)$ is defined to be $0$ for $(t_1,
\ldots,t_n)\in[0,t]^n \setminus T_n(t)$ where $T_n(t)=\{0<t_1<\cdots
<t_n<t\}$. Note that in Proposition~2.1 of \cite{B12}, it is assumed
that $u_0=1$ and $v_0=0$, so that $w_\mathrm{ w}=1$. This result continues
to hold when the function $w_\mathrm{ w}$ is given by (\ref{def-w}), since
$w_\mathrm{ w}$ does not depend on $x$.

\emph{Step} 2. Let $\widetilde{f}_{n}(\cdot,t,x)$ be the symmetrization\vspace*{1pt}
of $f_n(\cdot,t,x)$. By Remark~2.3 of \cite{B12}, if $\|\widetilde
{f}_n(\cdot,t,x)\|_{\cH^{\otimes n}}^2<\infty$, then $\widetilde
{f}_n(\cdot,t,x) \in\cH^{\otimes n}$ and the multiple Wiener integral
$I_n(f_n(\cdot,t,x))=I_n(\widetilde{f}_n(\cdot,t,x))$ is a well-defined
element of $\cH_n$.

\emph{Step} 3. Suppose that the series $\sum_{n \geq1}I_n(f_n(\cdot
,t,x))$ converges in $L^2(\Omega)$, that is,
%
\begin{equation}
\label{series-finite} \sum_{n \geq1}n! \bigl\|\widetilde{f}_n(
\cdot,t,x)\bigr\|_{\cH^{\otimes
n}}^{2}<\infty.
\end{equation}
Let
%
\begin{equation}
\label{Wiener-chaos-sol} u(t,x):=w(t,x)+\sum_{n \geq1}I_n
\bigl(f_n(\cdot,t,x) \bigr).
\end{equation}

\emph{Step} 4. Define
$v^{(t,x)}(s,\cdot)$ by relation \eqref{def-v}. This is a product
between the distribution $G_\mathrm{ w}(t-s,x-\cdot)$ and the function
$u(s,\cdot)$. The process $v^{(t,x)}$ has the Wiener chaos expansion
$v^{(t,x)}(\bullet)=\sum_{n \geq0}I_n(f_{n+1}(\cdot,\bullet,t,x))$,
where $\bullet$ denotes the missing $(s,y)$ variable; see the proof of
Theorem~2.8 in \cite{B12}. By Proposition~\ref{integr-crit}, $v^{(t,x)}
\in\operatorname{ Dom}  \delta$ and $\delta(v^{(t,x)})=\sum_{n \geq
0}I_{n+1}(f_{n+1}(\cdot,t,x))=u(t,x)-w(t,x)$. Hence the process $u=\{
u(t,x); t \geq0,x \in\bR^d\}$ with the Wiener chaos expansion \eqref
{Wiener-chaos-sol} is a solution of (\ref{wave}). Moreover,
%
\begin{equation}
\label{moment-2-sol} E\bigl|u(t,x)\bigr|^2=w(t,x)^2+\sum
_{n \geq1}\frac{1}{n!} \alpha_n(t),
\end{equation}
where $\alpha_n(t)=n!E|I_n(f_n(\cdot,t,x))|^2=(n!)^2 \|\widetilde
{f}_n(\cdot,t,x)\|_{\cH^{\otimes n}}^2$.

\emph{Step} 5. It remains to prove \eqref{series-finite}. When the
spatial covariance function $f$ is given by case (ii) above, this
follows by Proposition~3.4 of \cite{B12}. A similar argument can be
used for cases (i), (iii) and (iv); see Proposition~\ref{exist-sol} below.

Summarizing, to prove that a solution of \eqref{wave} exists in the
case $d \geq3$, we only need to show that the series $\sum_{n \geq
1}I_n(f_n(\cdot,t,x))$ converges in $L^2(\Omega)$; that is, \eqref
{series-finite} holds. In this case, one such solution is given by
\eqref{Wiener-chaos-sol}.

\begin{remark}
The uniqueness of the solution of \eqref{heat} for $d \geq3$ was
not treated in \cite{B12}. It may be possible to show that the solution
is unique in this case too. This would require significant
modifications to the method described above for the case $d \leq2$,
since both terms $G(t-s,x-y)$ and $f_n(\cdot,s,y)$ encountered in
definition \eqref{def-gn} of $g_n^{(t,x)}(\cdot,s,y)$ are distributions
in $y$. We do not investigate this problem here.
We note in passing that the classical method for proving uniqueness
does not seem to work for equations \eqref{wave} or \eqref{heat} when
the solution is interpreted in the sense of Definition~\ref
{definition-solution}. To see this, assume that there are two solutions
$u$ and $v$, and let $d=u-v$. Then
\[
d(t,x)=\int_0^t \int_{\bR^d}G(t-s,x-y)\,d(s,y)W(
\delta s,\delta y).
\]
The $L^{2}(\Omega)$-norm of the Skorohod integral above is a sum of two
terms, the second one involving the Malliavin derivative of $d$; see
\cite{nualart98}, relation (1.11). This second term vanishes when the
noise is white in time, 
but when the noise is fractional in time, it is not clear how to treat
this term.
\end{remark}

\begin{remark}
After examining (\ref{Fourier-f}), we infer that
in the case of equation \eqref{wave} with $d=3$, $f_n(t_1,\cdot
,\ldots
,t_n,\cdot,t,x)$ is a finite measure on $\bR^{3n}$ given by
\begin{eqnarray*}
&&f_n(t_1,\cdot,\ldots,t_n,\cdot,t,x)\nonumber\\
&&\qquad=
G(t-t_n,x-dx_n)G(t_n-t_{n-1},x_n-dx_{n-1})
\cdots
\\
& &\qquad\quad{}\times G(t_2-t_1,x_2-dx_1)w(t_1,x_1)1_{\{0<t_1<\cdots<t_n<t\}},\nonumber
\end{eqnarray*}
where for fixed $a \in\bR^3$, we denote by $G(t,a-\cdot)$ the measure
defined by
$G(t,a-\cdot)(A)=G(t,a-A)$ for all $A \in\cB(\bR^3)$.
\end{remark}

\begin{remark}
Notice that in both the hyperbolic and parabolic cases, the
function (or distribution) $f_n$ is stationary in the sense that, for
all $t_1,\ldots,t_n \in[0,t]$ and for any $x_1,\ldots,x_n,x \in\bR^d$,
\[
f_n(t_1,x_1,\ldots,t_n,x_n,t,x)
= f_n(t_1,x_1-x,\ldots,t_n,
x_n-x, t, 0).
\]
This remains valid for $\widetilde{f}_n$. A direct consequence is that
$\|\widetilde{f}_n(\cdot,t,x)\|_{\cH^{\otimes n}}^2$, and hence
$\alpha
_n(t)$, do not depend on $x$. Since the initial conditions are
constant, $w$ does not depend on $x$ either and the moments of $u$ are
independent of $x$. This justifies the definition of Lyapunov exponent
independent of $x$. Also, notice that it is possible to show that the
law of $u(t,x)$ is independent of $x$; see, for instance, \cite{dalang}
in the white noise case. These remarks are not true if the initial
conditions are not constant.
\end{remark}

We return now to series \eqref{series-finite}, which is also related to
the second moment of the solution $u(t,x)$; see \eqref{moment-2-sol}.
An important role in the present paper is played by the $n$th term of
this series, which depends on $\alpha_n(t)$. First, note that an
expression similar
to (\ref{def-cov2}) exists for the $n$-fold inner product $\langle
\cdot
, \cdot\rangle_{\cH^{\otimes n}}$. Using this expression, we have
%
\begin{equation}
\label{def-alpha} \alpha_n(t)= \alpha_H^n
\int_{[0,t]^{2n}} \prod_{j=1}^{n}|t_j-s_j|^{2H-2}
\psi_{n}(\mathbf{t},\mathbf{s})\,d\mathbf{t}\,d\mathbf{s},
\end{equation}
where we denote $\mathbf{t}=(t_1,\ldots,t_n)$ and $\mathbf{s}=(s_1, \ldots
,s_n)$, and we define
%
\begin{eqnarray}
\label{def-psi}&& \psi_{n}(\mathbf{t},\mathbf{s})=\int
_{\bR^{nd}} \cF g_\mathbf{t}^{(n)}(\cdot ,t,x) (
\xi_1, \ldots,\xi_n)
\nonumber
\\[-8pt]
\\[-8pt]
\nonumber
&&\hspace*{64pt}{}\times
\overline{\cF g_\mathbf{s}^{(n)}(
\cdot ,t,x) (\xi _1, \ldots,\xi_n)}\mu(d
\xi_1)\cdots\mu(d\xi_n)
\end{eqnarray}
with $g_\mathbf{t}^{(n)}(\cdot,t,x)=n! \widetilde{f}_n(t_1,\cdot,\ldots
,t_n,\cdot,t,x)$. Note that $\psi_{n}(\mathbf{t},\mathbf{s})$ depends also on
$t$, so that the correct notation should be $\psi_{n}(\mathbf{t},\mathbf{
s},t)$. To simplify the notation, we omit writing $t$ in $\psi
_{n}(\mathbf{
t},\mathbf{s},t)$.

An alternative calculation of the function $\psi_n(\mathbf{t},\mathbf{s})$ is
needed in Section~\ref{FK-section} below, for equation \eqref{wave}
with $d \leq3$. For this, let $\rho,\sigma\in S_n$ be such that
\[
0<t_{\rho(1)}<\cdots< t_{\rho(n)} <t \quad\mbox{and}\quad 0 <
s_{\sigma(1)}< \cdots< s_{\sigma(n)} < t,
\]
and denote $t_{\rho(n+1)}=s_{\sigma(n+1)}=t$. Then if $d \leq2$, we have
%
\begin{eqnarray}
\label{def-psi1} \psi_{n}(\mathbf{t},\mathbf{s})&=& \int
_{\bR^{2nd}} \prod_{j=1}^{n}
G(t_{\rho(j+1)}-t_{\rho
(j)},x_{\rho
(j+1)}-x_{\rho(j)})
w(t_{\rho(1)},x_{\rho(1)})
\nonumber\\
& &\hspace*{22pt} {}\times\prod_{j=1}^{n}G(s_{\sigma(j+1)}-s_{\sigma(j)},y_{\sigma
(j+1)}-y_{\sigma(j)})
w(s_{\sigma(1)},y_{\sigma(1)})\\
&&\hspace*{22pt}{}\times \prod_{j=1}^{n}f(x_j-y_j)
\,d\mathbf{x} \,d\mathbf{y},\nonumber
\end{eqnarray}
with the notation $\mathbf{x}=(x_1, \ldots,x_n),\mathbf{y}=(y_1, \ldots,y_n)
\in\bR^{nd}$, whereas if $d=3$,
%
\begin{eqnarray}
\label{def-psi2}  \psi_{n}(\mathbf{t},\mathbf{s})&=& \int
_{\bR^{2nd}} \prod_{j=1}^{n}
G(t_{\rho(j+1)}-t_{\rho
(j)},x_{\rho
(j+1)}-dx_{\rho(j)})
w(t_{\rho(1)},x_{\rho(1)})
\nonumber\\
& &\hspace*{22pt}{}\times \prod_{j=1}^{n}G(s_{\sigma(j+1)}-s_{\sigma(j)},y_{\sigma
(j+1)}-dy_{\sigma(j)})
w(s_{\sigma(1)},y_{\sigma(1)})\\
&&\hspace*{22pt}{}\times \prod_{j=1}^{n}f(x_j-y_j).\nonumber
\end{eqnarray}
(In both integrals above, we use the notation $x_{\rho(n+1)}=y_{\sigma
(n+1)}=x$.)

This concludes the summary of the results of \cite{B12} which are
needed here.

\section{Hyperbolic case: Existence of the solution}
\label{existence-section}

In this section, we prove the existence of a solution of equation
\eqref
{wave} [given by \eqref{Wiener-chaos-sol}] in any space dimension $d
\geq1$, when $f$ is a kernel of cases (i)--(iv). This yields the
conclusion of Theorem~\ref{wave-main}(a) and (b) (with $p=2$).

We let $G = G_\mathrm{ w}$ and $w = w_\mathrm{ w}$.
We introduce the following constant:
%
\begin{equation}
\label{def-K-mu} K({\mu}):=\sup_{\eta\in\bR^d} \int_{\bR^d}
\frac{1}{1+|\xi
-\eta|^2} \mu(d\xi).
\end{equation}
Note that $K({\mu})<\infty$ if and only if \eqref{cond-mu} holds; see
the proof of Lemma~8 in \cite{dalang-mueller03}.

Note that (\ref{cond-mu}) is satisfied in cases (i) and (iv). In cases
(ii) and (iii), (\ref{cond-mu}) holds if and only if
%
\begin{equation}
\label{cond-mu-1} a<2.
\end{equation}
We define a constant $K_\mathrm{ w}$ by
%
\begin{equation}
\label{def-K} K_\mathrm{ w} =\cases{ %
\mu
\bigl( \bR^d \bigr), & \quad$\mbox{in case (i),}$
\vspace*{2pt}\cr
4K({\mu}), & \quad $\mbox{in cases (ii) and (iii),}$
\vspace*{2pt}\cr
\pi,&\quad  $\mbox{in case (iv).}$}
\end{equation}

We have the following preliminary result.

\begin{lemma}
\label{estimate-psi}
Let $f$ be a kernels of cases \textup{(i)--(iv)}. Assume that (\ref{cond-mu}) holds.
For any $t>0$ and for any $\mathbf{t}=(t_1,\ldots,t_n)$ in $[0,t]^n$,
\[
\psi_{n}(\mathbf{t},\mathbf{t}) \leq(u_0+t
v_0)^2 K_\mathrm{ w}^n
(u_1,\ldots, u_n)^{2-a},
\]
where $a$ is given by (\ref{def-a}), $u_j=t_{\rho(j+1)}-t_{\rho(j)}$
for $j=1, \ldots,n$, $t_{\rho(1)} < \cdots< t_{\rho(n)}$ for some
$\rho\in S_n$, $t_{\rho(n+1)}=t$, and
$K_\mathrm{ w}$ is the constant defined in \eqref{def-K}.
\end{lemma}

\begin{pf}
As in the proof of Lemma~3.2 of \cite{B12},
by (\ref{def-psi}), (\ref{Fourier-f}) and (\ref{Fourier-G}), we obtain
\begin{eqnarray*}
&&\psi_n(\mathbf{t},\mathbf{t})
\\
&& \qquad =(u_0+t_{\rho(1)} v_0)^2 \int
_{\bR^{nd}} \frac{\sin^2(u_1|\xi
_1|)}{|\xi_1|^2} \cdots\\
&&\hspace*{86pt}\qquad\quad{}\times \frac{\sin^2(u_n|\xi_1+\cdots+\xi
_n|)}{|\xi
_1+\cdots+\xi_n|^2}\mu(d
\xi_1) \cdots\mu(d\xi_n).
\end{eqnarray*}
We consider separately the four cases:
\begin{itemize}
\item\emph{Case} (i). Using the fact that $|x^{-1}\sin x| \leq1$, we have
\[
\psi_n(\mathbf{t},\mathbf{t}) \leq(u_0+t
v_0)^2 \bigl[\mu \bigl(\bR^d \bigr)
\bigr]^n (u_1 ,\ldots, u_n)^2.
\]

\item\emph{Case} (ii). This case was treated in Lemma~3.2 of \cite{B12}.

\item\emph{Case} (iii). Let $c=c_{(\alpha_j)_j}$. Using the change of
variables $\eta_j=\xi_1+\cdots+\xi_j$,
\begin{eqnarray*}
\psi_n(\mathbf{t},\mathbf{t}) & = & c^n
(u_0+t_{\rho(1)} v_0)^2 \int
_{\bR
^{d}} \,d\eta_1 \frac{\sin^2(u_1|\eta_1|)}{|\eta_1|^2} \prod
_{j=1}^{d} |\eta_{1,j}|^{\alpha_j-1}
\\
&& {}\times\int_{\bR^{d}} \,d\eta_2 \frac{\sin^2(u_2|\eta
_2|)}{|\eta_2|^2}
\prod_{j=1}^{d} |\eta_{2,j}-
\eta_{1,j}|^{\alpha_j-1}
\\
&&{} \vdots
\\
&&{} \times\int_{\bR^{d}} \,d\eta_n \frac{\sin^2(u_n|\eta
_n|)}{|\eta_n|^2}
\prod_{j=1}^{d} |\eta_{n,j}-
\eta_{n-1,j}|^{\alpha_j-1},
\end{eqnarray*}
where $\eta_i=(\eta_{i,j})_{j=1, \ldots,d}$ with $\eta_{i,j} \in
\bR$.
Note that for any $t>0$ and $\eta\in\bR^d$,
\begin{eqnarray*}
&&c\int_{\bR^d}\frac{\sin^2(t|\xi|)}{|\xi|^2} \prod
_{j=1}^{d} |\xi _j-\eta
_j|^{\alpha_j-1}\,d\xi\\
&&\qquad= ct^{2-a} \int_{\bR^d}
\frac{\sin^2(|\xi|)}{|\xi|^2} \prod_{j=1}^{d} |\xi
_j-t\eta _j|^{\alpha_j-1}\,d\xi
\\
&&\qquad=t^{2-a} \int_{\bR^d}\frac{\sin^2(|\xi+t\eta|)}{|\xi+t\eta|^2}\mu(d\xi )
\\
&&\qquad\leq4 t^{2-a}\int_{\bR^d} \frac{\mu(d\xi)}{1+|\xi+t\eta|^2} \leq4
t^{2-a}K(\mu),
\end{eqnarray*}
since $(\sin(x)/x)^2 \leq4/(1+x^2)$ for all $x > 0$. Hence
\[
\psi_n(\mathbf{t},\mathbf{t}) \leq(u_0+tv_0)^2
\bigl(4K(\mu) \bigr)^n (u_1 ,\ldots,
u_n)^{2-a}.
\]

\item\emph{Case} (iv). Using the change of variables $\eta_j=\xi
_1+\cdots
+\xi_j$, we have
\[
\psi_n(\mathbf{t},\mathbf{t})=(u_0+t_{\rho(1)}v_0)^2
\int_{\bR^{n}} \frac
{\sin
^2(u_1|\eta_1|)}{|\eta_1|^2}\cdots\frac{\sin^2(u_n|\eta
_n|)}{|\eta
_n|^2} \,d
\eta_1 \cdots\, d\eta_n.
\]
Using \eqref{def-Gw}, \eqref{Fourier-G} and Plancherel's theorem, we
obtain that for any $t>0$,
\[
\int_{\bR}\frac{\sin^2(t|\xi|)}{|\xi|^2}\,d\xi=\pi t.
\]
%
Hence
\[
\psi_n(\mathbf{t},\mathbf{t}) = (u_0+tv_0)^2
\pi^n u_1 ,\ldots, u_n.
\]
\end{itemize}
\upqed\end{pf}

The following result is an extension of Proposition~3.1 of \cite
{dalang-mueller09} to the case of the fractional noise in time.

\begin{proposition}
\label{exist-sol}
Let $f$ be a kernel of cases \textup{(i)--(iv)}, and $\rho_\mathrm{ w},a,K_\mathrm{ w}$
be the constants given by (\ref{def-rho}), (\ref{def-a}), respectively
(\ref{def-K}). Assume that (\ref{cond-mu}) holds. Then:
\begin{longlist}[(a)]
\item[(a)] for any $t>0$ and for any integer $n \geq1$,
%
\begin{equation}
\label{UB-alpha} \alpha_n(t) \leq(u_0+t
v_0)^2 c^n K_\mathrm{
w}^n \frac
{t^{(2H+2-a)n}}{(n!)^{2-a}},
\end{equation}
where $\alpha_n(t)$ is given by \eqref{def-alpha} and $c$ is a constant
depending on $H$ and $a$;

\item[(b)] for any $d \geq1$, equation (\ref{wave}) has a solution
$u(t,x)$ [given by \eqref{Wiener-chaos-sol}] which has the following
property: for any $x \in\bR^d$ and for any $t>0$,
\[
E\bigl|u(t,x)\bigr|^2 \leq c_1(u_0+tv_0)^2
\exp \bigl(c_2 K_\mathrm{ w}^{1/(3-a)}t^{\rho
_\mathrm{ w}}
\bigr),
\]
where $c_1>0$ is a constant depending on $a$, and $c_2>0$ is a constant
depending on $H$ and $a$.
\end{longlist}
\end{proposition}


\begin{pf}
(a) We proceed as in the proof of Proposition~3.3 of \cite{B12}.

For any $\mathbf{t}=(t_1,\ldots,t_n) \in[0,t]^n$, we define $\beta(\mathbf{
t})=\prod_{j=1}^{n}u_j$, where $u_j=t_{\rho(j+1)}-t_{\rho(j)}$, and
$\rho\in S_n$ is chosen such that $t_{\rho(1)}<\cdots<t_{\rho(n)}$,
and $t_{\rho(n+1)}=t$.

By the Cauchy--Schwarz inequality, $\psi_n(\mathbf{t},\mathbf{s}) \leq\psi
_n(\mathbf{t},\mathbf{t})^{1/2}\psi_n(\mathbf{s},\mathbf{s})^{1/2}$. By Lem\-ma~\ref
{estimate-psi}, it follows that
%
\begin{equation}
\label{bound-psi-ts} \psi_n(\mathbf{t},\mathbf{s}) \leq
(u_0+t v_0)^2 K_\mathrm{
w}^n \bigl[\beta(\mathbf{t}) \beta(\mathbf{s}) \bigr]^{(2-a)/2}.
\end{equation}

Using definition (\ref{def-alpha}) of $\alpha_n(t)$ and (\ref
{bound-psi-ts}), we obtain
\[
\alpha_n(t) \leq(u_0+tv_0)^2
K_\mathrm{ w}^n \alpha_H^n \int
_{[0,t]^{2n}} \prod_{j=1}^{n}|t_j-s_j|^{2H-2}
\bigl[\beta(\mathbf{t}) \beta(\mathbf{ s}) \bigr]^{(2-a)/2} \,d\mathbf{t}\,d
\mathbf{s}.
\]

We now use the fact that for any $\varphi\in L^{1/H}(\bR^n)$,
%
\begin{equation}
\label{MMV-ineq} \alpha_H^n \int_{\bR^{2n}}
\prod_{j=1}^{n}|t_j-s_j|^{2H-2}
\bigl|\varphi (\mathbf{ t})\bigr|\bigl|\varphi(\mathbf{s})\bigr| \,d\mathbf{t}\,d\mathbf{s} \leq
b_{H}^n \biggl(\int_{\bR
^n}\bigl|\varphi(
\mathbf{t})\bigr|^{1/H}\,d\mathbf{t} \biggr)^{2H}\hspace*{-10pt}
\end{equation}
for some constant $b_H>0$; see Lemma~\ref{lemmaB3}, Appendix \ref{sec:appB}. We obtain
\begin{eqnarray*}
\alpha_n(t) & \leq& (u_0+tv_0)^2
K_\mathrm{ w}^n b_{H}^n \Biggl(\int
_{[0,t]^n} \prod_{j=1}^{n}
\beta(\mathbf{t})^{(2-a)/(2H)}\,d\mathbf{t} \Biggr)^{2H}
\\
&=& (u_0+tv_0)^2 K_\mathrm{
w}^n b_{H}^n \biggl(n!\int
_{T_n(t)} \bigl[(t-t_n)\cdots(t_2-t_1)
\bigr]^{(2-a)/(2H)}\,d\mathbf{t} \biggr)^{2H},
\end{eqnarray*}
where $T_n(t)=\{0<t_1<\cdots<t_n<t\}$. By Lemma~3.5 of \cite{BT10}, for
any $h>-1$,
\[
\int_{T_n(t)} \bigl[(t-t_n) (t_n-t_{n-1})
\cdots(t_2-t_1) \bigr]^h d\mathbf{t}=
\frac
{\Gamma(1+h)^{n+1}}{\Gamma((1+h)n+1)}t^{(1+h)n}.
\]
By Stirling's formula,
$\Gamma((1+h)n+1) \sim C_n (n!)^{1+h}$,
where $C_n$ is such that $\lambda^{-n} \leq C_n \leq\lambda^n$ for
some constant $\lambda>1$ depending on $h$; see the proof of Lemma A.1,
Appendix \ref{sec:app}.
Hence
\[
\int_{T_n(t)} \bigl[(t-t_n) (t_n-t_{n-1})
\cdots(t_2-t_1) \bigr]^h \,d\mathbf{t}\leq
\frac
{\Gamma(1+h)^{n} c_0^n }{(n!)^{1+h}}t^{(1+h)n}
\]
for some $c_0>0$.
In our case, $h=(2-a)/(2H)$.
We obtain:
\begin{eqnarray*}
\alpha_n(t) & \leq& (u_0+tv_0)^2
K_\mathrm{ w}^n b_{H}^n \biggl(n!
\frac
{\Gamma(1+h)^{n} c_1^{n}}{(n!)^{1+h}}t^{(1+h)n} \biggr)^{2H}
\\
& = & (u_0+tv_0)^2 K_\mathrm{
w}^n c^n \frac{1}{(n!)^{2-a}}t^{(2H+2-a)n},
\end{eqnarray*}
where $c= b_H \Gamma(1+h)^{2H}c_0^{2H}$ depends on $H$ and $a$.

(b) We use (\ref{moment-2-sol}) and the result from part (a). We obtain
that for any $t > 0$,
\[
E\bigl|u(t,x)\bigr|^2 \leq(u_0+t v_0)^2
\sum_{n \geq0} \frac{c^n K_\mathrm{ w}^n
t^{(2H+2-a)n}}{(n!)^{3-a}}.
\]
Since this series is convergent for any fixed $t > 0$, this proves the
existence result. Now, using Lemma~\ref{series-ineq} (Appendix \ref{sec:app}), we
have that for all $t > 0$
\[
E\bigl|u(t,x)\bigr|^2 \leq c_1 (u_0+t
v_0)^2 \exp \bigl(c_2'
\bigl(cK_\mathrm{ w}t^{2H+2-a} \bigr)^{1/(3-a)} \bigr),
\]
where $c_1>0$ and $c_2'>0$ are some constants depending on $a$. The
conclusion follows taking $c_2=c_2'c^{1/(3-a)}$.
\end{pf}

\section{Hyperbolic case: Upper bound on the moments}
\label{UB-section}

In this section, we give an upper bound for the moments of order $p >
2$ of a solution of equation \eqref{wave} [given by \eqref
{Wiener-chaos-sol}]. This yields the conclusion of Theorem~\ref{wave-main}(b).

Recall that this solution of \eqref{wave} has the Wiener chaos
expansion given by~\eqref{Wiener-chaos-sol}. This means that
$u(t,x)=\sum_{n \geq0}J_n(t,x)$ where $J_n(t,x)$ is in the $n$th
Wiener chaos $\cH_n$ associated to the noise $W$, and
\[
E\bigl|u(t,x)\bigr|^2=\sum_{n \geq0}E\bigl|J_n(t,x)\bigr|^2=
\sum_{n \geq0}\frac
{1}{n!}\alpha_n(t),
\]
where $\alpha_n(t)$ is defined in \eqref{moment-2-sol} and is estimated
by (\ref{UB-alpha}).

The following result is an extension of Theorem~3.2 of \cite
{dalang-mueller09} to the case of the fractional noise in time. 

\begin{proposition}
\label{UB-theorem}
Let $f$ be one of the kernels \textup{(i)--(iv)}, and $\rho_\mathrm{ w},a,K_\mathrm{ w}$
be the constants given by (\ref{def-rho}), (\ref{def-a}), respectively
(\ref{def-K}). Assume that (\ref{cond-mu}) holds. Let $u(t,x)$ be a
solution of \eqref{wave}, given by \eqref{Wiener-chaos-sol}.
Then for any $p \geq2$, for any $x \in\bR^d$ and for any $t >0$,
\[
E\bigl|u(t,x)\bigr|^p \leq c_1^p(u_0+tv_0)^{p}
\exp \bigl(c_2 K_\mathrm{ w}^{1/(3-a)}p^{(4-a)/(3-a)}
t^{\rho_\mathrm{ w}} \bigr),
\]
where $c_1>0$ is a constant depending on $a$, and $c_2>0$ is a constant
depending on $H$ and $a$.
\end{proposition}


\begin{pf}
When $p=2$, the result is given by Propostion
\ref{exist-sol}.

When $p>2$, we use the same idea as in the proof of Theorem~4.1 of
\cite{B12}.
We denote by $\|\cdot\|_{p}$ the $L^{p}(\Omega)$-norm. We use the fact
that for elements in a \emph{fixed} Wiener chaos $\cH_n$, the $\|\cdot
\|
_{p}$-norms are equivalent; see the last line of page 62 of~\cite
{nualart06} with $q=p$ and $p=2$. More precisely,
\[
\bigl\|J_n(t,x)\bigr\|_{p} \leq(p-1)^{n/2}
\bigl\|J_n(t,x)\bigr\|_2 =(p-1)^{n/2} \biggl(
\frac{1}{n!}\alpha_n(t) \biggr)^{1/2}.
\]
Using (\ref{UB-alpha}), we obtain
\[
\bigl\|J_n(t,x)\bigr\|_{p} \leq(u_0+t
v_0) C_{p,K_\mathrm{ w}}^{n} t^{n(2H+2-a)/2}
\frac{1}{(n!)^{(3-a)/2}},
\]
where $C_{p,K_\mathrm{ w}}=(p-1)^{1/2}c^{1/2}K_\mathrm{ w}^{1/2}$ and $c$
depends on $H$ and $a$.

Recall Minkowski's inequality for integrals (see Appendix A.1 of \cite
{stein70}),
\[
\biggl[\int_{Y} \biggl(\int_{X}
\bigl|F(x,y)\bigr| \mu(dx) \biggr)^p \nu(dy) \biggr]^{1/p} \leq\int
_{X} \biggl(\int_{Y}
\bigl|F(x,y)\bigr|^p \nu(dy) \biggr)^{1/p} \mu(dx).
\]
We use this inequality for $(X,\cX)=(\bN, 2^{\bN})$ with $\mu$ the
counting measure, $(Y,\cY,\nu)=(\Omega,\cF,P)$ and $F(n,\omega
)=J_n(\omega,t,x)$.
We have
\begin{eqnarray*}
\bigl\|u(t,x)\bigr\|_{p} & = & \biggl\|\sum_{n \geq0}J_n(t,x)
\biggr\|_{p} \leq\sum_{n
\geq0}
\bigl\|J_n(t,x) \bigr\|_{p}
\\
& \leq& (u_0+t v_0) \sum
_{n \geq0} \frac{C_{p,K_\mathrm{ w}}^{n}
t^{n(2H+2-a)/2}}{(n!)^{(3-a)/2}}.
\end{eqnarray*}
Using Lemma~\ref{series-ineq} (Appendix \ref{sec:app}), we infer that for any $t>0$,
\[
\bigl\|u(t,x)\bigr\|_p \leq c_1 (u_0+t
v_0) \exp \bigl\{c_2'
\bigl(C_{p,K_\mathrm{ w}} t^{(2H+2-a)/2} \bigr)^{2/(3-a)} \bigr\},
\]
where $c_1>0$ and $c_2'>0$ are some constants depending on $a$. The
conclusion follows taking $c_2=c_2' c^{1/(3-a)}$, since $\frac
{2H+2-a}{2} \cdot\frac{2}{3-a}=\rho_\mathrm{ w}$ and
\[
pC_{p,K_\mathrm{ w}}^{2/(3-a)}=p(p-1)^{1/(3-a)} c^{1/(3-a)}K_\mathrm{
w}^{1/(3-a)}.
\]
\upqed\end{pf}

%

\section{Hyperbolic case: FK representation for the second moment}
\label{FK-section}

In this section, we develop a Feynman--Kac (FK) representation for the
second moment of a solution $u(t,x)$ of the wave equation (\ref{wave})
[given by \eqref{Wiener-chaos-sol}], similar to the one obtained in
\cite{DMT08} in the case of white noise in time. Due to the fractional
component of the noise, our representation is based on a Poisson random
measure on $\bR_{+}^2$, rather than a simple Poisson process. This
extension follows the approach of \cite{B09} for the parabolic case.

The following theorem is the main result of this section. This theorem
is valid for any function $f$ for which covariance \eqref{def-cov} of
the noise $W$ is well defined, but may not be valid in case (iv) (since
in this case, $f$ is a distribution). Theorem~\ref{FK-theorem} will be
used in Section~\ref{LB-section} to obtain a lower bound for the second
moment of a solution to \eqref{wave} in cases (i)--(iii). Case (iv)
will be treated differently using an approximation based on case (ii).

\begin{theorem}
\label{FK-theorem}
Suppose that equation (\ref{wave}) with $d \leq3$ has a solution
$u(t,x)$ [given by \eqref{Wiener-chaos-sol}], where $W=\{W(\varphi
);\varphi\in\cH\}$ is a zero-mean Gaussian process with covariance
specified by \eqref{cov-W} and \eqref{def-cov}, and $f$ is a
nonnegative function on $\bR^d$, which is the Fourier transform of a
tempered measure $\mu$ on $\bR^d$. Then for any $t >0,x \in\bR^d$,
\begin{eqnarray*}
&&E\bigl|u(t,x)\bigr|^2
\\
&&\qquad =  e^{t^2} \sum_{n \geq0} \mathop{\sum
_{i_1, \ldots,i_n}}_{\mathrm{distinct}} E_{x} \Biggl[
w_\mathrm{ w} \bigl(t-\tau_n,X_{\tau
_n}^1
\bigr) w_\mathrm{ w} \bigl(t-\tau_n',X_{\tau_n'}^2
\bigr) \prod_{j=1}^{n}(\tau
_j-\tau _{j-1})
\\
& &\hspace*{76pt}\qquad\quad{} \times\prod_{j=1}^{n} \bigl(
\tau_j'-\tau_{j-1}' \bigr)
\prod_{j=1}^{n}f \bigl(X_{T_{i_j}}^1-X_{S_{i_j}}^2
\bigr) \alpha_H^n\\
&&\hspace*{180pt}{}\times \prod_{j=1}^{n}|T_{i_j}-S_{i_j}|^{2H-2}
1_{B_{i_1,\ldots
,i_n}(t)} \Biggr],
\end{eqnarray*}
where, by convention, the term for $n=0$ is taken to be $w_\mathrm{
w}(t,x)^2$, and $w_\mathrm{ w}$ is defined by \eqref{def-w}. Here:
\begin{itemize}
\item$N=\sum_{i \geq1}\delta_{P_i}$ is a Poisson random measure on
$\bR_{+}^2$ of intensity the Lebesgue measure, with $P_i=(T_i,S_i)$;
\item$B_{i_1,\ldots,i_n}(t)$ is the event that $N$ has points
$P_{i_1}, \ldots,P_{i_n}$ in $[0,t]^2$;
\item$\tau_1 \leq\cdots\leq\tau_n$ and $\tau_1'\leq\cdots\leq
\tau
_n'$ are the points $T_{i_1}, \ldots,T_{i_n}$, respectively $S_{i_1},
\ldots,S_{i_n}$ arranged in increasing order;
\item the processes $X^1=(X_s^1)_{s \in[0,t]}$ and $X^2=(X_s^2)_{s \in
[0,t]}$ are defined by (\ref{def-X1}) and~(\ref{def-X2}) below, and we
denote by $P_x$ a probability measure under which $X_0^1=X_0^2=x$.
($E_x$ stands for the expectation with respect to $P_x$.)
\end{itemize}
\end{theorem}

The processes $X^1$ and $X^2$ are constructed as in \cite
{DMT08}, using the coordinates of the points of $N$ on the two axes. We
explain this construction below.
On the event $B_{i_1,\ldots,i_n}(t)$, we arrange the two sets of points
$\{T_{i_1},\ldots,T_{i_n}\}$ and $\{S_{i_1},\ldots,S_{i_n}\}$ in
increasing order as $\tau_1 \leq\cdots\leq\tau_n$, respectively
$\tau
_1' \leq\cdots\leq\tau_n'$. More precisely, if we denote
$U_{j}=T_{i_j}$ and $V_j=S_{i_j}$ for $j=1,\ldots,n$, then there exist
some permutations $\rho$ and $\sigma$ of $\{1,\ldots,n\}$ such that
\[
U_{\rho(n)} \leq U_{\rho(n-1)} \leq\cdots\leq U_{\rho(1)} \qquad\mbox
{and}\qquad V_{\sigma(n)} \leq V_{\sigma(n-1)} \leq\cdots\leq V_{\sigma(1)}.
\]
We let $\tau_j=U_{\rho(n+1-j)}$ and $\tau_{j}'=V_{\sigma(n+1-j)}$ for
any $j=1, \ldots,n$.

We let $(\Theta_i^1)_{i \geq1}$ and $(\Theta_i^2)_{i \geq1}$ be two
independent i.i.d. collections of random variables with the same law as
$\Theta_0$, where $\Theta_0$ is a random variable with values in $\bR
^d$ such that if $d \leq2$, $\Theta_0$ has density function $G_\mathrm{
w}(1,\cdot)$, and if $d=3$, $\Theta_0$ has distribution $G_\mathrm{
w}(1,\cdot)$.
The importance of the variable $\Theta_0$ stems from the fact that for
any $t>0$,
%
\begin{equation}
\label{density-Theta} \frac{G_\mathrm{ w}(t,\cdot)}{t}\qquad \mbox{is the density$/$distribution of $t
\Theta_0$}.
\end{equation}

Using the points $\tau_1 \leq\cdots\leq\tau_n$ and the variables
$(\Theta_i^1)_{i \geq1}$, we construct the process $X^1=(X_s^1)_{s
\in
[0,t]}$ by setting
%
\begin{equation}
\label{def-X1}X_{s}^1 = X_{\tau_i}^1+
(s-\tau_{i})\Theta_{i+1}^1 \qquad\mbox{if }
\tau_i \leq s \leq\tau_{i+1}
\end{equation}
for any $1 \leq i \leq n$, where $\tau_0=0$, $\tau_{n+1}=t$ and
$X_0^1=0$. We use a similar construction for the process
$X^2=(X_s^2)_{s \in[0,t]}$ using the points $\tau_1' \leq\cdots\leq
\tau_n'$ and the variables $(\Theta_i^2)_{i \geq1}$, that is, $\tau
_0'=0$, $\tau_{n+1}'=t$, $X_0^2=0$ and for any $1 \leq i \leq n$,
%
\begin{equation}
\label{def-X2} X_{s}^2 = X_{\tau_i'}^2+
\bigl(s-\tau_{i}' \bigr)\Theta_{i+1}^2\qquad
\mbox{if } \tau_i' \leq s \leq\tau_{i+1}'.
\end{equation}

We now give some remarks about the statement of Theorem~\ref{FK-theorem}.

\begin{remark}
A similar formula can be obtained for $E[u(t,x)u(s,y)]$ using the
points of $N$ in $[0,t] \times[0,s]$ and assuming that $X_0^1=x$ and
$X_0^2=y$.
\end{remark}

\begin{remark}
\label{remark-points-N}
Note that $|T_{i_j}-S_{i_j}|^{2H-2}=\infty$ if the point
$(T_{i_j},S_{i_j})$ falls on the diagonal $D=\{(s,s);0 \leq s \leq t\}$
of the square $[0,t]^2$. This is not a problem since with probability
$1$, $N$ has no points in $D$: $P(N(D)=0)=e^{- \operatorname{Leb}(D)}=1$.
\end{remark}

\begin{remark}
Without loss of generality we may assume that $\tau_1 < \cdots<
\tau_n$ and $\tau_1' < \cdots< \tau_n'$ since the event 
for which $\tau_j=\tau_{j-1}$ (or $\tau_j'=\tau_{j-1}'$) for some $j=1,
\ldots,n$ has probability zero: with probability $1$, no vertical (or
horizontal) line contains two distinct points of $N$; see page 223 of
\cite{resnick07}.
\end{remark}

\begin{remark}
Theorem~\ref{FK-theorem} is valid for any function $f$, not
necessarily as in one of the cases (i)--(iii). In fact, this
representation remains valid if we replace $\alpha_H|t-s|^{2H-2}$ in
(\ref{def-cov}) by a function $\eta(t,s)$, provided that $\langle
\cdot
, \cdot\rangle_{\cH}$ defines an inner product. We only need to assume
that a solution of (\ref{wave}) [given by \eqref{Wiener-chaos-sol}]
exists. In the new representation, $\alpha_H|T_{i_j}-S_{i_j}|^{2H-2}$
is replaced by $\eta(t-T_{i_j},t-S_{i_j})$.
\end{remark}

We now introduce the necessary ingredients for the proof of Theorem~\ref
{FK-theorem}.

Recall first that if $(N_t)_{t \geq0}$ is a Poisson process on $\bR
_{+}$ of rate $1$ with jump times $\tau_1<\tau_2<\cdots,$ then the
conditional distribution of $(\tau_1,\ldots,\tau_n)$ given $N_t=n$
coincides with the distribution of the order statistics of a sample of
size $n$ from the uniform distribution on $[0,t]$. This property lies
at the core of the FK formula obtained in \cite{DMT08} and
can be seen very easily as follows. For any $t>0$ fixed, the process
$(N_s)_{s \in[0,t]}$ can be constructed as $N_s=\sum_{i=1}^{Y}1_{\{X_i
\leq s\}}$, where $(X_i)_{i \geq1}$ are i.i.d. random variables with a
uniform distribution on $[0,t]$, and $Y$ is an independent Poisson
random variable with mean $t$. If $N_t=n$, the jump times of $N$ in
$[0,t]$ coincide with the order statistics $X_{(1)}<\cdots<X_{(n)}$.

A similar property holds for the planar Poisson process. This basic
observation has enabled the first author to obtain in \cite{B09} an FK
formula similar to that of \cite{DMT08} in the case of the heat
equation with fractional noise in time.
In this section, we develop a similar formula for the wave equation
with $d \leq3$.

More precisely, let $N$ be a Poisson random measure as in Theorem~\ref
{FK-theorem}. Since the Lebesgue measure does not have any atoms, $N$
is a.s. simple, that is, \mbox{$N(\{ \mathbf{t}\}) \leq1$} for all $\mathbf{
t}=(t_1,t_2) \in\bR_+^2$ a.s. (see Exercise 2.4 of \cite
{kallenberg83}). This means that with probability $1$, the points
$(P_i)_{i \geq1}$ are distinct. For any $t>0$ fixed, we consider the
event $B_{i_1,\ldots,i_n}(t)$ for distinct indices $i_1, \ldots,i_n
\geq1$.

The following result plays an important role in the present paper; see
also Problem 5.2, page 162 of \cite{resnick07}.

\begin{lemma}
\label{points-uniform}
Let $N=\sum_{i \geq1}\delta_{P_i}$ be a Poisson random measure on
$\bR
_{+}^2$ of intensity the Lebesque measure, with $P_i=(S_i,T_i)$. For
$t>0$ and distinct indices $i_1, \ldots, i_n$, let $B_{i_1, \ldots
,i_n}(t)$ be the event that $N$ has points $P_{i_1}, \ldots,P_{i_n}$ in
$[0,t]^2$.
Given $B_{i_1,\ldots,i_n}(t)$, both vectors
$(P_{i_1}, \ldots,P_{i_n})$ and $(\mathbf{t}-P_{i_1},\ldots, \mathbf{
t}-P_{i_n})$ have a uniform distribution on $[0,t]^{2n}$, where $\mathbf{
t}=(t,t) \in\bR_{+}^2$.
\end{lemma}

\begin{pf} The restriction of $N$ to $[0,t]^2$ can be
constructed\vspace*{1pt} as $N=\sum_{i=1}^{Y}\delta_{X_i}$, where $(X_i)_{i \geq1}$
are i.i.d. random variables with a uniform distribution on $[0,t]^2$,
and $Y$ is an independent Poisson random variable with mean $t^2$. If
$N$ has points $P_{i_1}, \ldots,P_{i_n}$ in $[0,t]^2$, the vector
$(P_{i_1}, \ldots, P_{i_n})$ of the $n$ points coincides with a vector
$(X_{j_1}, \ldots,X_{j_n})$ for some distinct indices $j_1, \ldots
,j_n$, which clearly has a uniform distribution on $[0,t]^{2n}$. The
argument for the vector $(\mathbf{t}-P_{i_1},\ldots,\mathbf{t}-P_{i_n})$ is
similar; see Lemma~2.1 of \cite{B09} for an alternative proof.
\end{pf}

As a consequence of the previous lemma, any $n$-fold integral over
$([0,t]^{2})^{n}$ of a deterministic function $F$ has a stochastic
representation based on the points of~$N$; see page 257 of \cite{B09}
for the proof.

\begin{corollary}
\label{corol-points-uniform}
For any measurable function $F\dvtx [0,t]^{2n} \to\bR$ which is either
bounded or nonnegative, we have
\begin{eqnarray*}
&&\int_{[0,t]^{2n}}F(t_1,s_1,
\ldots,t_n,s_n)\,d\mathbf{t}\,d\mathbf{ s}
\\
& &\qquad= n! e^{t^2} \mathop{\sum_{i_1,\ldots,i_n}}_{\mathrm{distinct}}E
\bigl[F(t-T_{i_1},t-S_{i_1}, \ldots,t-T_{i_1},t-S_{i_n})
1_{B_{i_1,\ldots,i_n}(t)} \bigr],
\end{eqnarray*}
where $\mathbf{t}=(t_1,\ldots,t_n)$ and $\mathbf{s}=(s_1, \ldots,s_n)$ with
$t_i \in[0,t]$ and $s_i \in[0,t]$.
\end{corollary}

The next result gives a stochastic representation for the $n$th term of
series (\ref{moment-2-sol}).

\begin{lemma}
\label{stoch-repr-alpha}
For any $t>0$ and for any integer $n \geq1$, we have
\begin{eqnarray*}
\alpha_n(t)& = & n! \alpha_H^n
e^{t^2} \\
&&{}\times \mathop{\sum_{{i_1, \ldots,i_n}}}_{\mathrm{distinct}}
E \Biggl[ \prod_{j=1}^{n}|T_{i_j}-S_{i_j}|^{2H-2}
\\
& &\hspace*{50pt}{} \times\psi_{n}(t-T_{i_1}, \ldots ,t-T_{i_n},t-S_{i_1},
\ldots,t-S_{i_n}) 1_{B_{i_1, \ldots,i_n}(t)} \Biggr],
\end{eqnarray*}
where $\psi_n(\mathbf{t},\mathbf{s})$ is given by \eqref{def-psi}.
\end{lemma}

\begin{pf} The integral on the right-hand side of (\ref
{def-alpha}) can be represented in the desired form by applying
Corollary~\ref{corol-points-uniform} to the function
\[
F(t_1,s_1,\ldots,t_n,s_n)=
\alpha_H^n \prod_{j=1}^{n}|t_j-s_j|^{2H-2}
\psi _{n}(\mathbf{t},\mathbf{s}).
\]
\upqed\end{pf}

The next result will be used to evaluate the term $\psi_{n}(t-T_{i_1},
\ldots, t-T_{i_n},  t-S_{i_1}, \ldots, t-S_{i_n})$ which
appears in Lemma~\ref{stoch-repr-alpha}. For simplicity, we work first
with some nonrandom points $(t_1,s_1), \ldots, (t_n,s_n)$ in
$[0,t]^{2}$. These points will be replaced later by $(T_{i_1},S_{i_1}),
\ldots, (T_{i_n},S_{i_n})$.

\begin{lemma}
\label{calcul-psi}
Let $(t_1,s_1), \ldots, (t_n,s_n) \in[0,t]^{2}$. Let $\rho,\sigma
\in
S_n$ be such that
\[
0<t_{\rho(n)}< \cdots<t_{\rho(1)}<t \quad\mbox{and}\quad 0<s_{\sigma
(n)}<
\cdots<s_{\sigma(1)}<t.
\]
If $d \leq2$, then
\begin{eqnarray*}
& & \psi_n(t-t_{1},\ldots,t-t_n,t-s_1,
\ldots,t-s_n)
\\
& &\qquad =\int_{\bR^{2nd}} \,d\mathbf{z} \,d\mathbf{w} \prod
_{j=1}^{n} f \Biggl(\sum
_{k=1}^{n+1-\rho^{-1}(j)}z_k-\sum
_{k=1}^{n+1-\sigma^{-1}(j)}w_k \Biggr)
\\
& &\qquad\quad{} \times G_\mathrm{ w}(t_{\rho(n)},z_1)G_\mathrm{
w}(t_{\rho
(n-1)}-t_{\rho(n)},z_2) \cdots G_\mathrm{
w}(t_{\rho(1)}-t_{\rho
(2)},z_n)
\\
& &\qquad\quad{} \times G_\mathrm{ w}(s_{\sigma(n)},w_1)G_\mathrm{
w}(s_{\sigma
(n-1)}-s_{\sigma(n)},w_2) \cdots G_\mathrm{
w}(s_{\sigma(1)}-s_{\sigma
(2)},w_n)
\\
& &\qquad\quad{} \times w \Biggl(t-t_{\rho(1)},x+\sum_{k=1}^{n}z_k
\Biggr) w \Biggl(t-s_{\sigma(1)},x+\sum_{k=1}^{n}w_k
\Biggr),
\end{eqnarray*}
where $\mathbf{z}=(z_1,\ldots,z_n)$ and $\mathbf{w}=(w_1,\ldots,w_n)$ with
$z_i \in\bR^d,w_i \in\bR^d$. A similar relation holds for $d=3$,
replacing $G_\mathrm{ w}(t_{\rho(n)},z_1)\,dz_1$ by $G_\mathrm{ w}(t_{\rho
(n)},dz_1)$, etc.
\end{lemma}

\begin{pf} Assume first that $d \leq2$. We use the
alternative definition \eqref{def-psi1} of $\psi_{n}(\mathbf{t},\mathbf{s})$.
We proceed as in the first part of the proof of Lemma~2.2 of \cite
{B09}. Denote $t_{\rho(n+1)}=s_{\sigma(n+1)}=0$ and $x_{\rho
(n+1)}=y_{\sigma(n+1)}=x$. Note that
\[
0<t-t_{\rho(1)}< \cdots<t-t_{\rho(n)}<t \quad\mbox{and}\quad
0<t-s_{\sigma(1)}< \cdots<t-s_{\sigma(n)}<t
\]
%
and $G_\mathrm{ w}(t,x)=G_\mathrm{ w}(t,-x)$. By definition, $\psi
_n(t-t_{1},\ldots,t-t_n,t-s_1,\ldots,t-s_n)$ is equal to
\begin{eqnarray*}
&&\int_{\bR^{2nd}} \,d\mathbf{x} \,d\mathbf{y} \prod
_{j=1}^{n}G_\mathrm{ w}(t_{\rho(j)}-t_{\rho(j+1)},x_{\rho(j)}-x_{\rho(j+1)})
w(t-t_{\rho
(1)},x_{\rho(1)})
\\
& & \qquad{}\times\prod_{j=1}^{n}G_\mathrm{
w}(s_{\sigma(j)}-s_{\sigma
(j+1)},y_{\sigma(j)}-y_{\sigma(j+1)})
w(t-s_{\sigma(1)},y_{\sigma(1)}) \\
& & \qquad{}\times\prod_{j=1}^{n}f(x_j-y_j).
\end{eqnarray*}

 The result follows by the change of variables
$x_{\rho(j)}-x_{\rho(j+1)}=z_{n+1-j}$ and $y_{\sigma(j)}-y_{\sigma
(j+1)}=w_{n+1-j}$ for $j=1, \ldots,n$.

The same argument works also for $d=3$, using the alternative
definition \eqref{def-psi2} of $\psi_{n}(\mathbf{t},\mathbf{s})$. To see
this, assume for simplicity that $n=2$, $0<t_1<t_2<t$ and
$0<s_2<s_1<t$. (The same argument applies in the general case.)
Then $\psi_2(t-t_1,t-t_2,t-s_1,t-s_2)$ is equal to
\begin{eqnarray*}
&& \int_{\bR^{4d}} h(x_1,x_2,y_1,y_2)
G_\mathrm{ w}(t_1,dx_1-x)G_\mathrm{
w}(t_2-t_1,dx_2-x_1)
\\
&&\qquad{} \times G_\mathrm{ w}(s_2,dy_2-x)G_\mathrm{
w}(s_1-s_2,dy_1-y_2),
\end{eqnarray*}
where $h(x_1,x_2,y_1,y_2)=f(x_1-y_1)f(x_2-y_2)w(t-t_2,x_2)w(t-s_1,y_1)$
and we used the fact that $G_\mathrm{ w}(t,a-dx)=G_\mathrm{ w}(t,dx-a)$. We
claim that for any nonnegative measurable function $\varphi\dvtx \bR^{4d}
\to\bR$,
%
\begin{eqnarray}\label{ch-var}
&& \int_{\bR^{4d}} \varphi(x_1,x_2,y_1,y_2)
G_\mathrm{ w}(t_1,dx_1-x)G_\mathrm{
w}(t_2-t_1,dx_2-x_1)
\nonumber
\\
&&{}\quad \times G_\mathrm{ w}(s_2,dy_2-x)G_\mathrm{
w}(s_1-s_2,dy_1-y_2)
\nonumber
\\[-8pt]
\\[-8pt]
\nonumber
&&\qquad =  \int_{\bR^{4d}} \varphi(x+z_1,x+z_1+z_2,
x+w_1+w_2,x+w_1)
\nonumber
\\
&&\hspace*{18pt}\qquad\quad{} \times G_\mathrm{ w}(t_1,dz_1)G_\mathrm{
w}(t_2-t_1,dz_2)G_\mathrm{
w}(s_2,dw_1)G_\mathrm{ w}(s_1-s_2,dw_2).
\nonumber
\end{eqnarray}
(This means that we can apply informally the change of variables
$x_1-x=z_1, x_2-x_1=z_2$ and $y_2-x=w_1, y_1-y_2=w_2$.) Assuming that
$\varphi(x_1,x_2,y_1,y_2)=\phi_1(x_1)\phi_2(x_2)\psi_1(y_1)\psi
_2(y_2)$, relation (\ref{ch-var}) follows using the fact that for any
nonnegative measurable function $\phi\dvtx \bR^d \to\bR$,
\[
\int_{\bR^d}\phi(x)G_\mathrm{ w}(t,dx-a)=\int
_{\bR^d}\phi(a+y)G_\mathrm{ w}(t,dy).
\]
The case of an arbitrary function $\varphi$ follows by approximation.
The conclusion follows applying (\ref{ch-var}) to the function
$\varphi
=h$.
\end{pf}

\begin{remark}
\label{heat-calcul-psi}
In the case of the heat equation, $G_\mathrm{ h}(t-s,\cdot)$ is the
density of $B_t^1-B_s^1$, where $(B_t^1)_{t \geq0}$ is a
$d$-dimensional Brownian motion, and
the product
\[
G_\mathrm{ h}(t_{\rho(n)},z_1)G_\mathrm{
h}(t_{\rho(n-1)}-t_{\rho(n)},z_2) \cdots G_\mathrm{
h}(t_{\rho(1)}-t_{\rho(2)},z_n),
\]
which appears in Lemma~\ref{calcul-psi} is the density of the random
vector
\[
\bigl(B_{t_{\rho(n)}}^1, B_{t_{\rho(n-1)}}^1-B_{t_{\rho(n)}}^1,
\ldots,B_{t_{\rho
(1)}}^1-B_{t_{\rho(2)}}^1 \bigr).
\]
Applying a similar argument for the other $n$-term product (depending
on $s$) and using an independent Brownian motion $(B_t^2)_{t \geq0}$,
we infer that
\begin{eqnarray*}
&&\psi_n(t\mathbf{e}-\mathbf{t},t\mathbf{e}-\mathbf{s})
\\
& &\qquad=  E \Biggl[w \bigl(t-t_{\rho(1)}, x+B_{t_{\rho(1)}}^1
\bigr)w \bigl(t-s_{\sigma(1)}, x+B_{s_{\sigma(1)}}^2 \bigr)\prod
_{j=1}^{n}f \bigl(B_{t_j}^1-B_{s_j}^2
\bigr) \Biggr],
\end{eqnarray*}
where $\mathbf{e}=(1,\ldots,1) \in\bR^n$, $\mathbf{t}=(t_1,\ldots,t_n)$ and
$\mathbf{s}=(s_1, \ldots,s_n)$. Something similar will happen in the case
of the wave equation, conditionally on $N$.
\end{remark}


\begin{remark}
\label{cond-density-Y}
Due to (\ref{density-Theta}), when $d \leq2$, the product
\[
\frac{G_\mathrm{ w}(\tau_1,z_1)}{\tau_1}\cdot\frac{G_\mathrm{ w}(\tau
_2-\tau
_1,z_2)}{\tau_2-\tau_1} \cdots\frac{G_\mathrm{ w}(\tau
_n-\tau
_{n-1},z_n)}{\tau_n-\tau_{n-1}}
\]
is the conditional density of $\mathbf{Y}^1=(X_{\tau_1}^1,X_{\tau
_2}^1-X_{\tau_1}^1, \ldots,X_{\tau_n}^1-X_{\tau_{n-1}}^1)$ given $N$.
Let $\mathbf{Y}^2=(X_{\tau_1'}^2,X_{\tau_2'}^2-X_{\tau_1'}^2, \ldots
,X_{\tau_{n}'}^2-X_{\tau_{n-1}'}^2)$.
Since $X^1$ and $X^2$ are conditionally independent given $N$,
\[
\prod_{j=1}^{n}\frac{G_\mathrm{ w}(\tau_j-\tau_{j-1},z_j)}{\tau_j-\tau
_{j-1}} \prod
_{j=1}^{n}\frac{G_\mathrm{ w}(\tau_j'-\tau
_{j-1}',w_j)}{\tau
_j'-\tau_{j-1}'}
\]
is the conditional density of $(\mathbf{Y}^1,\mathbf{Y}^2)$ given $N$. A
similar thing happens when $d=3$.
Therefore, for the wave equation, the processes $X^1,X^2$ play the same
role (conditionally on $N$), as the Brownian motions $B^1,B^2$ for the
heat equation; see Remark~\ref{heat-calcul-psi}.
\end{remark}

\begin{pf*}{Proof of Theorem~\ref{FK-theorem}}
By applying Lemma~\ref{calcul-psi} to the points $(t_j,s_j)=(T_{i_j},S_{i_j})$ we obtain
that on the event $B_{i_1,\ldots,i_n}(t)$,
\begin{eqnarray*}
&&\psi_n(t-T_{i_1}, \ldots,t-T_{i_n},t-S_{i_1},
\ldots ,t-S_{i_n})
\\
& & \qquad=\int_{\bR^{2nd}} \,d\mathbf{z} \,d\mathbf{w} \prod
_{j=1}^{nd}f \Biggl(\sum
_{k=1}^{n+1-\rho^{-1}(j)}z_k-\sum
_{k=1}^{n+1-\sigma
^{-1}(j)}w_k \Biggr)
\\
& &\qquad\quad{} \times G_\mathrm{ w}(\tau_1,z_1)G_\mathrm{
w}(\tau_2-\tau _1,z_2) \cdots
G_\mathrm{ w}(\tau_n-\tau_{n-1},z_n)
\\
& &\qquad\quad{} \times G_\mathrm{ w} \bigl(\tau_1',w_1
\bigr)G_\mathrm{ w} \bigl(\tau_2'-\tau
_1',w_2 \bigr) \cdots G_\mathrm{
w} \bigl(\tau_n'-\tau_{n-1}',w_n
\bigr)
\\
& &\qquad\quad{} \times w_\mathrm{ w} \Biggl(t-\tau_n,x+\sum
_{k=1}^{n}z_k \Biggr)
w_\mathrm{ w} \Biggl(t-\tau_n',x+\sum
_{k=1}^{n}w_k \Biggr),
\end{eqnarray*}
assuming that $d \leq2$. A similar identity holds for $d=3$ replacing
$G_\mathrm{ w}(\tau_1,z_1)\,dz_1$ by
$G_\mathrm{ w}(\tau_1,dz_1)$, and so on. Inside this integral, we multiply
and divide by $\prod_{j=1}^{n}(\tau_j-\tau_{j-1}) \prod_{j=1}^{n}(\tau
_j'-\tau_{j-1}')$.

We assume that $X_0^1=X_0^2=0$. Using Remark~\ref{cond-density-Y}, we
infer that on the event $B_{i_1,\ldots,i_n}(t)$, $\psi_n(t-T_{i_1},
\ldots,t-T_{i_n},t-S_{i_1}, \ldots,t-S_{i_n})$ is equal to the
conditional expectation of
\begin{eqnarray*}
& & \prod_{j=1}^{n}f \Biggl(\sum
_{k=1}^{n+1-\rho^{-1}(j)} \bigl(X_{\tau
_k}^1-X_{\tau_{k-1}}^1
\bigr)- \sum_{k=1}^{n+1-\rho^{-1}(j)}
\bigl(X_{\tau_k'}^2-X_{\tau
_{k-1}'}^2 \bigr)
\Biggr)
\\
& &\qquad{} \times w_\mathrm{ w} \Biggl(t-\tau_n,x+\sum
_{k=1}^{n} \bigl(X_{\tau
_k}^1-X_{\tau_{k-1}}^1
\bigr) \Biggr) w_\mathrm{ w} \Biggl(t-\tau_n',x+
\sum_{k=1}^{n} \bigl(X_{\tau_k'}^2-X_{\tau_{k-1}'}^2
\bigr) \Biggr)
\\
& &\qquad{} \times\prod_{j=1}^{n}(
\tau_j-\tau_{j-1}) \prod_{j=1}^{n}
\bigl(\tau _j'-\tau_{j-1}'
\bigr)
\end{eqnarray*}
given $N$. Note that
\[
\sum_{k=1}^{n+1-\rho^{-1}(j)} \bigl(X_{\tau_k}^1-X_{\tau_{k-1}}^1
\bigr)= X_{\tau_{n+1-\rho^{-1}(j)}}^1 \quad\mbox{and}\quad \sum
_{k=1}^{n} \bigl(X_{\tau_k}^1-X_{\tau_{k-1}}^1
\bigr)= X_{\tau_{n}}^1
\]
(these are telescopic sums whose first term is $X_{\tau_0}^1=0$).
Recall that $\tau_{k}=U_{\rho(n+1-k)}$ for any $k=1, \ldots,n$ (where
$U_j=T_{i_j}$). Hence
\[
\tau_{n+1-\rho^{-1}(j)}=U_{\rho(n+1-n-1+\rho^{-1}(j))}=U_{\rho
(\rho^{-1}(j))}= U_j=T_{i_j}.
\]
A similar argument applies to the terms depending on $X^2$. We obtain
that on the event $B_{i_1,\ldots,i_n}(t)$,
\begin{eqnarray*}
&&\psi_n(t-T_{i_1}, \ldots,t-T_{i_n},t-S_{i_1},
\ldots ,t-S_{i_n})
\\
& &\qquad= E \Biggl[\prod_{j=1}^{n}f
\bigl(X_{T_{i_j}}^1-X_{S_{i_j}}^2 \bigr)
w_\mathrm{ w} \bigl(t-\tau_n,x+X_{\tau_n}^1
\bigr) w_\mathrm{ w} \bigl(t-\tau_n',x+X_{\tau_n'}^2
\bigr)
\\
& &\hspace*{122pt}\qquad\quad{} \times\prod_{j=1}^{n}(
\tau_j-\tau_{j-1}) \prod_{j=1}^{n}
\bigl(\tau_j'-\tau_{j-1}'
\bigr) \bigg\vert N \Biggr].
\end{eqnarray*}

Looking now back at the representation of $\alpha_n(t)$ (Lemma~\ref
{stoch-repr-alpha}), we obtain
\begin{eqnarray*}
&&\frac{1}{n!}\alpha_n(t)
\\
& &\qquad=  e^{t^2} \mathop{\sum_{{i_1,\ldots,i_n}}}_{\mathrm{ distinct}}
E \Biggl[1_{B_{i_1, \ldots,i_n}(t)}\prod_{j=1}^{n}|T_{i_j}-S_{i_j}|^{2H-2}
\\
& &\hspace*{56pt}\quad\qquad{} \times E \Biggl[\prod_{j=1}^{n}f
\bigl(X_{T_{i_j}}^1-X_{S_{i_j}}^2 \bigr)
w_\mathrm{ w} \bigl(t-\tau_n,x+X_{\tau
_n}^1
\bigr)\\
&&\hspace*{129pt}{}\times w_\mathrm{ w} \bigl(t-\tau_n',x+X_{\tau_n'}^2
\bigr)
\\
& &\hspace*{98pt}\qquad\quad{} \times\prod_{j=1}^{n}(
\tau_j-\tau _{j-1}) \prod_{j=1}^{n}
\bigl(\tau_j'-\tau_{j-1}'
\bigr) \bigg\vert N \Biggr] \Biggr].
\end{eqnarray*}

 Note that $ 1_{B_{i_1, \ldots,i_n}(t)} \prod_{j=1}^{n}|T_{i_j}-S_{i_j}|^{2H-2}$ is measurable with respect to $N$,
and so, this term goes inside the conditional expectation with respect
to $N$. The result follows using the fact that $E[E[ \cdot|N]]=E[
\cdot]$ and taking the sum over $n \geq1$. In the final step, the
values $x+X_{\tau_n}^1$ and $x+X_{\tau_n'}^2$ are replaced by
$X_{\tau
_n}^1$, respectively $X_{\tau_n'}^2$, under the probability measure $P_x$.
\end{pf*}

\section{Hyperbolic case: Lower bound on the moment of order 2}
\label{LB-section}

In this section, we give a lower bound for the second moment of a
solution $u$ to \eqref{wave} [given by \eqref{Wiener-chaos-sol}], when
$f$ is a kernel of cases (i)--(iv). This yields the conclusion of
Theorem~\ref{wave-main}(c).


For cases (i)--(iii), we follow the approach of Dalang and Mueller
\cite
{dalang-mueller09}.
This means that for any $x,y \in\bR^d$ with $x \neq y$, we consider
the solid (infinite) cone $C(x,y)$ in $\bR^d$, with vertex $y$, axis
oriented in the direction of the vector $x-y$ and an angle of $\pi/4$
between the axis and any lateral side.
This cone has the following properties:
\begin{longlist}[(iii)]
\item[(i)] if $|z-y| \leq\delta$, $|y-x|\leq\delta$ and $z \in
C(x,y)$, then $|z-x| \leq\delta$;
\item[(ii)] $y+z \in C(x,y)$ if and only if $y+rz \in C(x,y)$ for
any $r>0$;
\item[(iii)] $C(x,y)+z=C(x+z,y+z)$.
\end{longlist}
%

\subsection{Case \textup{(i)}: Spatially smooth noise}
In this case, since $f$ is continuous at~$0$, $\lim_{x \to
0}f(x)=f(0)=\mu(\bR^d)=K_\mathrm{ w}$. Letting $\alpha_0=K_\mathrm{ w}/2$, we
infer that there exists $\delta>0$ such that
%
\begin{equation}
\label{AssC-DM} f(x) \geq\alpha_0\qquad \mbox{for all } x \in
\bR^d,|x| \leq2\delta;
\end{equation}
that is, $f$ satisfies Assumption C of \cite{dalang-mueller09}. We
assume that $\delta$ is a rational number.

The next result corresponds to Theorem~\ref{wave-main}(c), in case
(i). Its proof relies on the Feynman--Kac formula developed in Section~\ref{FK-section}.

\begin{theorem}
\label{LB-th-f-bounded}
Let $f$ be a kernel of case \textup{(i)}. Then, for any $x \in\bR^d$ and for
any $t>0$,
\[
E\bigl|u(t,x)\bigr|^2 \geq c_3 u_0^2
\exp \bigl(c_5 K_\mathrm{ w}^{1/3} t^{\rho_\mathrm{ w}}
\bigr),
\]
where $c_3>0$ and $c_5>0$ are some constants depending on $H$, and the
constants $\rho_\mathrm{ w}$ and $K_\mathrm{ w}$ are given by (\ref{def-rho}),
respectively (\ref{def-K}).
\end{theorem}

\begin{pf} We proceed as in the proof of Theorem~4.1 of
\cite{dalang-mueller09}. To facilitate the comparison with the proof of
these authors, we use the same notation; that is, we denote $n$ by $k$
in the statement of Theorem~\ref{FK-theorem} above. We let $N_t=N([0,t]^2)$.

\textit{Step} 1. \emph{First step for the lower bound of $E|u(t,x)|^2$.}

Let $k \in\bZ_{+}$ be a large enough value (depending on $t$) such that
%
\begin{equation}
\label{def-m} m:=k\delta\in\bZ_{+},
\end{equation}
where $\delta$ is given by \eqref{AssC-DM}. (The precise value of $k$
will be given in step 8 below.) Notice that for all $t > 0$ and $x \in
\bR^d$, $w_\mathrm{ w}(t,x) = u_0 + v_0 t \geq u_0$, since $v_0 \geq0$.
Hence, by Theorem~\ref{FK-theorem},
%
\begin{eqnarray}\label{LB-step1}
E\bigl|u(t,x)\bigr|^2 &\geq& u_0^2
e^{t^2} \alpha_H^k\nonumber\\
&&{}\times  \mathop{\sum
_{
{i_1, \ldots,i_k}}}_{\mathrm{ distinct}} E_x \Biggl[\prod
_{j=1}^{k}(\tau_j-
\tau_{j-1}) \prod_{j=1}^{k}
\bigl( \tau_j'-\tau_{j-1}'
\bigr)
\\
& &\hspace*{56pt}{} \times\prod_{j=1}^{k}f
\bigl(X_{T_{i_j}}^1-X_{S_{i_j}}^2 \bigr)
\prod_{j=1}^{k}|T_{i_j}-S_{i_j}|^{2H-2}
1_{B_{i_1,\ldots,i_k}(t)} \Biggr].\nonumber
\end{eqnarray}

\textit{Step} 2. \emph{The event $D(t)$.}

We consider the event $D(t)=D^1(t) \cap D^2(t)$, where
\begin{eqnarray*}
D^1(t)&= & \bigl\{X_{\tau_j}^1+
\Theta_{j+1}^1 \in C \bigl(x,X_{\tau_j}^1
\bigr) \mbox{ for all } j=1,\ldots, k-1 \bigr\} \cap B_{i_1, \ldots,i_k}(t),
\\
D^2(t)&= & \bigl\{X_{\tau_j'}^2+
\Theta_{j+1}^2 \in C \bigl(x,X_{\tau_j'}^2
\bigr) \mbox{ for all } j=1,\ldots, k-1 \bigr\} \cap B_{i_1,\ldots,i_k}(t).
\end{eqnarray*}

On the event $D^{1}(t)$, if we assume that $\tau_j-\tau_{j-1} \leq
\delta$ for all $j=1, \ldots, k$, then
%
\begin{equation}
\label{X-close-x} \bigl|X_{\tau_j}^1-x\bigr| \leq\delta\qquad\mbox{for all
$j=1, \ldots,k$}.
\end{equation}

We first prove \eqref{X-close-x} by induction on $j$, using the
properties of the cone. The argument is the same as in \cite
{dalang-mueller09}. We include it for the sake of completeness. If
$j=1$, then $X_{\tau_1}^1=x+\tau_1 \Theta_1^1$ and $|X_{\tau
_1}^{1}-x|=\tau_1|\Theta_1^1| \leq\delta$ since $|\Theta_1^1| \leq1$.
Assume now that $|X_{\tau_j}^1-x| \leq\delta$. We use property (i) of
the cone for points $x'=0,y'=X_{\tau_j}^1-x$ and $z'=X_{\tau
_{j+1}}^1-x$. We note that
$|z'-y'|=|X_{\tau_{j+1}}^1-X_{\tau_j}^1|=(\tau_{j+1}-\tau_j)|\Theta
_{j+1}^1| \leq\delta$, and $|y'-x'|=|X_{\tau_j}^1-x| \leq\delta$ by
the induction hypothesis.
We also have $z' \in C(x',y')$, that is, $X_{\tau_{j+1}}^1-x \in
C(0,X_{\tau_j}^1-x)$. [This is equivalent to $X_{\tau_{j+1}}^1 \in
C(0,X_{\tau_j}^1-x)+x=C(x,X_{\tau_j}^1)$, using property (iii) of the
cone for the last equality, which is in turn equivalent to $X_{\tau
_j}^1+(\tau_{j+1}-\tau_j) \Theta_{j+1}^1 \in C(x,X_{\tau_{j}}^{1})$,
using the definition of $X_{\tau_{j+1}}^1$. By property (ii) of the
cone, this last property is equivalent to $X_{\tau_j}^1 +\Theta_{j+1}^1
\in C(x,X_{\tau_j}^{1})$, which holds true on the event $D^1(t)$.] By
property (i) of the cone, it follows that $|z'-x'| \leq\delta$, that
is, $|X_{\tau_{j+1}}^1-x| \leq\delta$.
This completes the proof of \eqref{X-close-x}.

Recall that $\tau_j=U_{\rho(k+1-j)}$ for some permutation $\rho$ of
$\{
1,\ldots,k\}$, where $U_{j}=T_{i_j}$. As $j$ runs through the set $\{1,
\ldots,k\}$, so does the value $\rho(k+1-j)$. Therefore, on the event
$D^{1}(t)$, if we assume that $\tau_j-\tau_{j-1} \leq\delta$ for all
$j=1, \ldots, k$, then
$|X_{T_{i_j}}^1-x| \leq\delta$ for all $j=1, \ldots, k$, by \eqref
{X-close-x}.
A similar property holds for $X^2$ on the event $D^2(t)$. Hence, on the
event $ D(t)$, if we assume that $\tau_j-\tau_{j-1} \leq\delta$ and
$\tau_j'-\tau_{j-1}' \leq\delta$ for all $j=1, \ldots, k$, then
\[
\bigl|X_{T_{i_j}}^1-X_{S_{i_j}}^2\bigr| \leq2 \delta\qquad
\mbox{for all } j=1, \ldots,k
\]
and so
%
\begin{equation}
\label{LB-f}f \bigl(X_{T_{i_j}}^1-X_{S_{i_j}}^2
\bigr) \geq\alpha_0\qquad \mbox {for all } j=1, \ldots,k.
\end{equation}

\begin{figure}[t]

\includegraphics{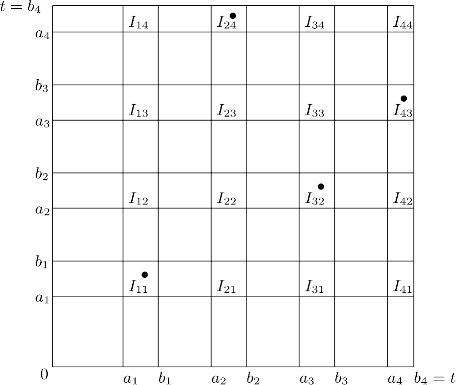}

\caption{The islands $I_{j,l}$ (for $k=4$) with points situated on the
islands $I_{11},I_{24},I_{32},I_{43}$ corresponding to the permutation
$(l_1,l_2,l_3,l_4)=(1,4,2,3)$.}
\label{island-figure}
\end{figure}

\textit{Step} 3. \emph{The islands $(I_{j,l})_{1 \leq j,l \leq k}$.}

The idea of the proof is to build some small islands around the $k$
points of the process $N$ in the region $[0,t]^2$. Figure~\ref{island-figure} shows these islands for $k=4$. To define these islands,
we let $\varepsilon=\frac{\delta t}{m+1}$ and $t_j=j \varepsilon$ for
any $j=1, \ldots,k$.
Due to (\ref{def-m}), we have
\[
t_k=k\varepsilon=\frac{m}{m+1}t \approx t \qquad\mbox{if } m \mbox
{ is large}.
\]
We consider the intervals $I_{j}=[a_j,b_j]$ with $j=1, \ldots,k$, where
$a_j=t_{j}-\varepsilon/4$ for $j=1, \ldots,k$, $b_j=t_j+\varepsilon/4$
if $j \leq k-1$, and $b_k=t$. For any $j,l=1, \ldots,k$, we define
\[
I_{j,l}=I_j \times I_l.
\]
The area of each square island $I_{j,l}$ is greater than $(\varepsilon
/4)^2$. In both the horizontal and vertical directions, the islands are
separated by intervals of length $\varepsilon/2$.

\textit{Step} 4. \emph{The event $C_{i_1,\ldots,i_k}(t)$.}

Let $C_{i_1,\ldots,i_k}(t)$ be the event that $N$ has points $P_{i_1},
\ldots,P_{i_k}$ in $[0,t]^2$ located on the islands $I_{1,l_1}, \ldots,
I_{k,l_k}$, for some permutation $(l_1, \ldots,l_k)$ of $\{1, \ldots
,k\}$.

Clearly, $C_{i_1,\ldots,i_k}(t)$ is included in $B_{i_1,\ldots
,i_k}(t)$. Notice that on the event $C_{i_1, \ldots,i_k}(t)$, it is not
possible to have two points $(T_{i_p},S_{i_p})$ and $(T_{i_q},S_{i_q})$
of $N$ in $[0,t]^2$ such that $T_{i_p},T_{i_q}$ are in the same
interval $I_j$ or $S_{i_p},S_{i_q}$ are in the same interval $I_l$.
Therefore, on the event $C_{i_1, \ldots,i_k}(t)$,
for any $j=1, \ldots,k$, we have $\tau_j \in I_j, \tau_j' \in I_j$,
and hence
%
\begin{equation}
\label{bound-tau} \frac{\varepsilon}{2} \leq\tau_{j}-\tau_{j-1}
\leq2 \varepsilon\quad \mbox{and}\quad \frac{\varepsilon}{2} \leq\tau_{j}'-
\tau_{j-1}' \leq 2 \varepsilon.
\end{equation}
In particular, if
%
\begin{equation}
\label{m-large} m \geq m_0(t):=[2t-1],
\end{equation}
then $\tau_j-\tau_{j-1} \leq\delta$ and $\tau_j'-\tau_{j-1}' \leq
\delta$ for all $j=1, \ldots,k$. It follows that
%
\begin{equation}
\label{LB-f-holds} \mbox{inequality (\ref{LB-f}) holds on the event $D(t) \cap
C_{i_1,
\ldots,i_k}(t)$},
\end{equation}
provided that $m \geq m_0=m_0(t)$.

\textit{Step} 5. \emph{Second step for the lower bound of $E|u(t,x)|^2$.}\vspace*{1pt}

On the event $B_{i_1,\ldots,i_k}(t)$, we define $\widetilde
{Z}_t=\prod_{j=1}^{k}(\tau_j-\tau_{j-1}) \prod_{j=1}^{k}(\tau_j'-\tau_{j-1}')$.
Using (\ref{LB-step1}) and (\ref{LB-f-holds}), we obtain
\begin{eqnarray*}
&&E\bigl|u(t,x)\bigr|^2
\\
&&\qquad \geq u_0^2 e^{t^2} \alpha_H^k
\mathop{\sum_{i_1,\ldots, i_k}}_{\mathrm{ distinct}}
E_{x} \Biggl[\widetilde{Z}_{t} \prod
_{j=1}^{k}f \bigl(X_{T_{i_j}}^1-X_{S_{i_j}}^2
\bigr)
\\
&&\hspace*{83pt}\qquad\quad{} \times\prod_{j=1}^{k}|T_{i_j}-S_{i_j}|^{2H-2}
1_{D(t)} 1_{C_{i_1,\ldots
,i_k}(t)} \Biggr]
\\
&&\qquad \geq u_0^2 e^{t^2} \alpha_H^k
\alpha_0^k \mathop{\sum_{i_1, \ldots, i_k}}_\mathrm{distinct} E_x \Biggl[\widetilde{Z}_t
\prod_{j=1}^{k}|T_{i_j}-S_{i_j}|^{2H-2}
1_{D(t)} 1_{C_{i_1, \ldots
,i_k}(t)} \Biggr]
\\
&&\qquad= u_0^2 e^{t^2} \alpha_H^k
\alpha_0^k \mathop{\sum_{i_1, \ldots, i_k}}_{\mathrm{distinct}} E_x \Biggl[\widetilde{Z}_t
\prod_{j=1}^{k}|T_{i_j}-S_{i_j}|^{2H-2}
1_{C_{i_1, \ldots
,i_k}(t)}P_x \bigl[D(t)|N \bigr] \Biggr].
\end{eqnarray*}

 Since the events $D^1(t)$ and $D^2(t)$ are conditionally
independent given $N$,
\[
P_x \bigl[D(t)|N \bigr]=P_x \bigl[D^1(t)|N
\bigr] P_x \bigl[D^2(t)|N \bigr].
\]
Using the properties of the cone and the independence of $(\Theta
_i^1)_{i \geq1}$, it can be shown that
$P_x[D^1(t)|N]=\gamma^{N_t-1}$,
where $\gamma=P(y+\Theta_0 \in C(0,y)) \in(0,1)$ does not depend on $y
\in\bR^d$. Note that $\gamma$ depends on $d$. A similar property holds
for $D^2(t)$. Hence,
\[
P_{x} \bigl[D(t)|N \bigr]=\gamma^{2(N_t-1)} >
\gamma^{2N_t}.
\]
Combining this with the previous lower bound for $E|u(t,x)|^2$, we obtain
\begin{eqnarray*}
E\bigl|u(t,x)\bigr|^2 &\geq& u_0^2 e^{t^2}
\alpha_H^k \alpha_0^k
\gamma^{2k} \\
&&{}\times \mathop{\sum_{i_1, \ldots, i_k}}_{\mathrm{distinct}} E_x \Biggl[\widetilde{Z}_t
\prod_{j=1}^{k}|T_{i_j}-S_{i_j}|^{2H-2}
1_{C_{i_1, \ldots
,i_k}(t)} \Biggr].
\end{eqnarray*}

We define the conditional expectation of a random variable $X$ with
respect to an event $B$ by $E[X|B]=E[X 1_{B}]/P(B)$. (This is not the
same as $E[X|\cG]$, where $\cG=\sigma(\{B\})=\{\varnothing
,B,B^c,\Omega\}
$ since $E[X|\cG]=E[X|B]1_{B}+E[X|B^c]1_{B^c}$.) In our case, $X$ is
the random variable appearing in the expectation above and
$B=B_{i_1,\ldots,i_k}(t)$. We obtain
\begin{eqnarray*}
E\bigl|u(t,x)\bigr|^2& \geq& u_0^2
e^{t^2} \alpha_H^k \alpha_0^k
\gamma ^{2k}
\\
& & {}\times\mathop{\sum_{i_1, \ldots, i_k}}_{\mathrm{ distinct}}
E_x \Biggl[\widetilde{Z}_t \prod
_{j=1}^{k}|T_{i_j}-S_{i_j}|^{2H-2}
1_{C_{i_1, \ldots,i_k}(t)}\bigg \vert B_{i_1,\ldots,i_k}(t) \Biggr] \\
&&\hspace*{38pt}{}\times P_x
\bigl(B_{i_1,\ldots,i_k}(t) \bigr).
\end{eqnarray*}

Note that by (\ref{bound-tau}), on the event $C_{i_1, \ldots,i_k}(t)$,
we have
$\widetilde{Z}_t \geq(\varepsilon/2)^{2k}$. Using the fact that
$\delta
=m/k$ [by the definition (\ref{def-m}) of $m$], we see that
%
\begin{equation}
\label{LB-epsilon} \frac{\varepsilon}{2}= \frac{\delta t}{2(m+1)} =\frac{m}{m+1} \cdot
\frac{t}{2k} \geq\frac{ct}{k}
\end{equation}
with $c=1/8$. Hence $\widetilde{Z}_t \geq(ct/k)^{2k}$ and
\begin{eqnarray*}
E\bigl|u(t,x)\bigr|^2 &\geq &e^{t^2} \alpha_H^k
\alpha_0^k \gamma^{2k} \biggl(
\frac{ct}{k} \biggr)^{2k}
\\
& &{} \times\mathop{\sum_{i_1, \ldots, i_k}}_{\mathrm{ distinct}}
E_x \Biggl[ \prod_{j=1}^{k}|T_{i_j}-S_{i_j}|^{2H-2}
1_{C_{i_1, \ldots,i_k}(t)} \bigg\vert B_{i_1,\ldots,i_k}(t) \Biggr]\\
&&\hspace*{38pt}{}\times P_x
\bigl(B_{i_1,\ldots,i_k}(t) \bigr).
\end{eqnarray*}

Since both $T_{i_j}$ and $S_{i_j}$ are in $[0,t]$, we obviously have
$|T_{i_j} - S_{i_j}| < t$. Thus since $2H-2 < 0$,
\[
\prod_{j=1}^{k}|T_{i_j}-S_{i_j}|^{2H-2}
> t^{(2H-2)k}.
\]
This turns out to be enough for our purposes. With this bound, we have
%
\begin{eqnarray}
\label{LB-step5} E\bigl|u(t,x)\bigr|^2 & \geq& u_0^2
e^{t^2} \alpha_H^k \alpha_0^k
\gamma^{2k} \biggl(\frac{ct}{k} \biggr)^{2k}
t^{(2H-2)k}
\nonumber
\\[-8pt]
\\[-8pt]
\nonumber
& &{} \times\mathop{\sum_{i_1, \ldots, i_k}}_{\mathrm{ distinct}}
P_x \bigl(C_{i_1, \ldots,i_k}(t)|B_{i_1,\ldots,i_k}(t)
\bigr)P_x \bigl(B_{i_1,\ldots
,i_k}(t) \bigr).
\end{eqnarray}

\textit{Step} 6. \emph{The conditional probability $P_x(C_{i_1,\ldots
,i_k}(t)|B_{i_1, \ldots,i_k}(t))$.}

Let $S_k$ be the set of all permutations $(l_1,\ldots,l_k)$ of $\{
1,\ldots,k\}$. By the definition of the event $C_{i_1,\ldots,i_k}(t)$,
\[
P_x \bigl(C_{i_1,\ldots,i_k}(t)|B_{i_1, \ldots,i_k}(t) \bigr)=\sum
_{(l_1,\ldots
,l_k) \in S_k}P_x \bigl(A_{i_1,\ldots,i_k}
\bigl(t,(l_1,\ldots ,l_k) \bigr)|B_{i_1,\ldots,i_k}(t)
\bigr),
\]
where $A_{i_1,\ldots,i_k}(t,(l_1,\ldots,l_k))$ is the event that $N$
has points $P_{i_1}, \ldots,P_{i_k}$ in $[0,t]^2$ located on the
islands $I_{1,l_1}, \ldots, I_{k,l_k}$. Note that
\[
A_{i_1,\ldots,i_k} \bigl(t,(l_1,\ldots,l_k) \bigr)=
\bigcup_{(j_1,\ldots,j_k)
\in
S_k}\{P_{i_1} \in
I_{j_1,l_1}, \ldots, P_{i_k} \in I_{j_k,l_k}\}.
\]

Given $B_{i_1,\ldots,i_k}(t)$, $(P_{i_1},\ldots,P_{i_k})$ has a uniform
distribution on $[0,t]^{2k}$. Hence
\begin{eqnarray*}
&&P_x \bigl(P_{i_1} \in I_{j_1,l_1}, \ldots,
P_{i_k} \in I_{j_k,l_k}|B_{i_1,\ldots,i_k}(t) \bigr) \\
&&\qquad =
\frac{\operatorname{ Leb}(I_{j_1,l_1}
\times\cdots\times I_{j_k,l_k})}{\operatorname{ Leb}([0,t]^{2k})}
\\
& &\qquad\geq \frac{1}{t^{2k}} \biggl( \frac{\varepsilon}{4} \biggr)^{2k}.
\end{eqnarray*}
Since the last quantity does not depend on the permutations
$(j_1,\ldots
,j_k)$ and $(l_1,\ldots,l_k)$, we obtain that
%
\begin{equation}
\label{cond-prob-C-B} P_x \bigl(C_{i_1,\ldots,i_k}(t)|B_{i_1, \ldots,i_k}(t)
\bigr)=(k!)^2 \frac
{1}{t^{2k}} \biggl( \frac{\varepsilon}{4}
\biggr)^{2k} \geq(k!)^2 \biggl(\frac{c}{k}
\biggr)^{2k},
\end{equation}
using (\ref{LB-epsilon}) for the inequality. Relation (\ref
{cond-prob-C-B}) is the analogue of (4.7) of \cite{dalang-mueller09}
(with $n=2$) for the fractional noise.

\textit{Step} 7. \emph{Third step for the lower bound of $E|u(t,x)|^2$.}

Combining (\ref{LB-step5}) and (\ref{cond-prob-C-B}), we get
\begin{eqnarray*}
E\bigl|u(t,x)\bigr|^2 & \geq& u_0^2 e^{t^2}
\alpha_H^k \alpha_0^k
\gamma^{2k} \biggl(\frac{ct}{k} \biggr)^{2k}
t^{(2H-2)k} (k!)^2 \biggl(\frac{c}{k}
\biggr)^{2k}
\\
& & {}\times\mathop{\sum_{i_1, \ldots,i_k}}_{\mathrm{ distinct}}
P_x \bigl(B_{i_1,\ldots,i_k}(t) \bigr).
\end{eqnarray*}
We now use the fact that $\{N_t=k\}$ is the disjoint union of all
events $B_{i_1,\ldots,i_k}(t)$ for all sets $\{i_1, \ldots,i_k\}$ of
cardinality $k$. Moreover, $N_t$ has a Poisson distribution with mean
$t^2$. Hence $P(N_t=k)= e^{-t^2} t^{2k}/k!$ and
\begin{eqnarray*}
E\bigl|u(t,x)\bigr|^2 & \geq& u_0^2 e^{t^2}
\alpha_H^k \alpha_0^k
\gamma^{2k} \biggl(\frac{ct}{k} \biggr)^{2k}
t^{(2H-2)k} (k!)^2 \biggl(\frac{c}{k}
\biggr)^{2k} e^{-t^2}\frac{t^{2k}}{k!}
\\
&=& u_0^2 \bigl(\alpha_0
\alpha_H \gamma^{2}c^4 \bigr)^k
t^{(2H+2)k} \frac
{1}{k^{4k}} k!.
\end{eqnarray*}
By Stirling's formula, there exists some $k_0 \geq1$ such that
$k! \geq e^{-k} k^k$ for all $k \geq k_0$. It follows that if $k \geq
k_0$, then
%
\begin{equation}
\label{final-LB} E\bigl|u(t,x)\bigr|^2 \geq u_0^2
\biggl(\alpha_0 c_H \frac
{t^{2H+2}}{k^{3}}
\biggr)^k,
\end{equation}
where $c_H=\alpha_H \gamma^{2}c^4 e^{-1}$ depends on $H$. ($c_H$
depends also on $d$, through $\gamma$.)

\textit{Step} 8. \emph{The choice of $k$.}

Let
\[
k= \bigl[e^{-1/3} \alpha_0^{1/3}
c_H^{1/3}t^{(2H+2)/3} \bigr],
\]
where $[x]=k \in\bZ$ if $k \leq x <k+1$.
Since $k \leq e^{-1/3} \alpha_0^{1/3} c_H^{1/3}t^{(2H+2)/3}$, it
follows that $e \leq\alpha_0 c_H t^{2H+2}/k^3$.
On the other hand, letting
\[
k_1=\tfrac{1}{2} \bigl(e^{-1}\alpha_0
c_H \bigr)^{1/3}=\alpha_0^{1/3}
\cdot \tfrac
{1}{2} \bigl(\alpha_H\gamma^2
c^4 e^{-2} \bigr)^{1/3}=:\alpha_0^{1/3}c_1^*,
\]
we have
$k > 2k_1 t^{(2H+2)/3} -1 \geq k_1 t^{(2H+2)/3}$ if $k_1 t^{(2H+2)/3}
\geq1$.
Using (\ref{final-LB}), we infer that
\[
E\bigl|u(t,x)\bigr|^2 \geq u_0^2 e^k
\geq u_0^2 \exp \bigl(k_1 t^{(2H+2)/3}
\bigr)\qquad \mbox{if } \alpha_0 t^{2H+2} \geq t_1':=
\bigl(c_1^{*} \bigr)^{-3}.
\]
Note that $k \geq k_0$ if $\alpha_0 t^{2H+2} \geq t_1'':=k_0^3 e c_H$.
We take $t_1=t_1' \vee t_1''$.

Let $c_4=c_1^* \alpha_0^{1/3}=c_1^* 2^{-1/3}K_\mathrm{ w}^{1/3}$. This
proves that for any $t>0$ such that $\alpha_0 t^{2H+2} \geq t_1$ (i.e.,
for all $t \geq t_0$ for some $t_0>0$),
\[
E\bigl|u(t,x)\bigr|^2 \geq u_0^2 \exp
\bigl(c_4 t^{\rho_\mathrm{ w}} \bigr).
\]

\textit{Step} 9. \emph{Extension to all $t>0$.}

Using \eqref{moment-2-sol} and the fact that $u_0>0$ and $v_0>0$, we
infer that for any $0<t<t_0$,
\[
E\bigl|u(t,x)\bigr|^2 \geq w(t,x)^2=(u_0+tv_0)^2
\geq c_3^{*} u_0^2 \exp
\bigl(c_4 t_0^{\rho_\mathrm{ w}} \bigr) \geq
c_3^* u_0^2 \exp \bigl(c_4
t^{\rho_\mathrm{ w}} \bigr),
\]
where $c_3^{*} =\exp(-c_4 t_0^{\rho_\mathrm{ w}})$. Finally, we let
$c_3=\min(1,c_3^{*})$ and $c_5=c_1^* 2^{-1/3}$.
\end{pf}


\subsection{Cases \textup{(ii)} and \textup{(iii)}: Fractional noise in space}

These cases are treated similar to case (i), using Theorem~\ref{FK-theorem}.
The difference is that instead of \eqref{AssC-DM}, we use the fact that
for any $\delta>0$,
%
\begin{equation}
\label{AssC-DM-caseii} f(x) \geq\alpha_0(\delta):=(2\delta)^{-a}
\qquad\mbox{for all } x \in \bR^d,|x|\leq2\delta,
\end{equation}
where $a$ is given by (\ref{def-a}).

The next result corresponds to Theorem~\ref{wave-main}(c), in cases
(ii)--(iii).

\begin{theorem}
\label{LB-th-f-Riesz}
Let $f$ be a kernel of either case \textup{(ii)} or \textup{(iii)}. If \eqref{cond-mu-1}
holds, then for any $x \in\bR^d$ and for any $t > 0$,
\[
E\bigl|u(t,x)\bigr|^2 \geq c_3 u_0^2
\exp \bigl(c_4 t^{\rho_\mathrm{ w}} \bigr),
\]
where $c_3>0$ and $c_4>0$ are some constants depending on $H$ and $a$,
and the constants
$\rho_\mathrm{ w}$ and $a$ are given by (\ref{def-rho}), respectively
(\ref{def-a}).
\end{theorem}

\begin{pf} We use the same argument as in the proof of
Theorem~\ref{LB-th-f-bounded}, but with a different method of
specifying the parameters.

More precisely, we let $k \in\bZ_{+}$ be a large enough value
(depending on $t$) which will be chosen later.
We choose $\delta=m/k$ where $m=[2t]$. This ensures that (\ref{def-m})
and (\ref{m-large}) are satisfied. Note that $\delta$ depends on $t/k$.

Let $c_H=\alpha_H \gamma^2 c^4 e^{-1}$. Relation (\ref{final-LB}) says
that if $k \geq k_0$, then
\begin{eqnarray*}
E\bigl|u(t,x)\bigr|^2 & \geq& u_0^2 \biggl(
c_H (2\delta)^{-a}\frac
{t^{2H+2}}{k^3}
\biggr)^k = u_0^2 \biggl( c_H
2^{-a} \biggl(\frac{m}{k} \biggr)^{-a}
\frac{t^{2H+2}}{k^3} \biggr)^k
\\
& \geq& u_0^2 \biggl( c_H
2^{-a} \biggl(\frac{2t}{k} \biggr)^{-a}
\frac
{t^{2H+2}}{k^3} \biggr)^k = u_0^2
\biggl( c_H^*\frac
{t^{2H+2-a}}{k^{3-a}} \biggr)^k,
\end{eqnarray*}
where $c_H^*=c_H 4^{-a}$. We let
\[
k= \bigl[ \bigl(e^{-1} c_H^* t^{2H+2-a}
\bigr)^{1/(3-a)} \bigr].
\]
(This choice will ensure that $\delta$ is small since $\delta\approx
2t/k \approx C t^{1-\rho_\mathrm{ w}}$ and $\rho_\mathrm{ w} > 1$.)
Then $e \leq c_H^* t^{2H+2-a}/k^{3-a}$. On the other hand, letting
%
\begin{equation}
\label{def-c2} c_4=\tfrac{1}{2} \bigl(e^{-1}c_H^*
\bigr)^{1/(3-a)},
\end{equation}
we have
$k > 2c_4 t^{\rho_\mathrm{ w}} -1 \geq c_4 t^{\rho_\mathrm{ w}}$ if
$c_4t^{\rho
_\mathrm{ w}} \geq1$.
Hence
\[
E\bigl|u(t,x)\bigr|^2 \geq u_0^2 e^k
\geq u_0^2 \exp \bigl(c_4 t^{\rho_\mathrm{
w}}
\bigr) \qquad\mbox{for all } t \geq t_0,
\]
where
%
\begin{equation}
\label{def-t2} t_0= \bigl(e c_H^{-1}
2^{3+a} \bigr)^{1/(2H+2-a)}.
\end{equation}

 For $0<t<t_0$, we argue as in step 9 of the proof of
Theorem~\ref{LB-th-f-bounded}.
\end{pf}

\subsection{Case \textup{(iv)}: White noise in space}

In this case, we cannot apply directly Theorem~\ref{FK-theorem} since
$f$ is not a function. Instead of this, we will use an approximation
technique based on case (ii).

The fact that we use this approximation may be surprising, since in
many instances, it is easier to deal with the white noise than a
correlated noise. This is due to the fact that our method for proving
the lower bound relies on the representation given by Theorem~\ref
{FK-theorem}. Obtaining a similar representation in the case $f=\delta
_0$ is more delicate. (The Dirac distribution would have to be
approximated in some sense, so that the representation make sense.)
Instead of this, we decided to use an approximation directly for
obtaining the lower bound.

Our procedure can be viewed as another method to smoothen the noise,
paralleling the method used in \cite{hu_nualart} and \cite{HHNT14}.
The fact that we approximate $\delta_0$ by a Riesz kernel allows us to
use the result that we proved for case (ii). A~more standard procedure
in the literature is to approximate $\delta_0$ by the heat kernel
$p_{\varepsilon}(x)=(2\pi\varepsilon)^{-1/2} \exp
(-|x|^2/(2\varepsilon
))$ as $\varepsilon\downarrow0$. This is a kernel of case (i), with
$\lim_{x \to0}p_{\varepsilon}(x)=(2\pi\varepsilon)^{-1/2}$. Denoting
by $u_{\varepsilon}(t,x)$ the solution of \eqref{wave} driven by a
noise $W_{\varepsilon}$ with spatial covariance $p_{\varepsilon}$, one
infers by Theorem~\ref{LB-th-f-bounded} that $E|u_{\varepsilon}(t,x)|^2
\geq u_0^2 \exp(c_1 \alpha_0(\varepsilon) t^{(2H+2)/3})$,
with $\alpha_0(\varepsilon)=(2\pi\varepsilon)^{-1/2}/2$. However, this
approximation is not suitable for our purposes, since $\lim_{\varepsilon\downarrow0}\alpha_0(\varepsilon)=\infty$.

We begin to explain this approximation technique.
For any $a \in(0,1)$, let $W_{a}=\{W_{a}(\varphi);\varphi\in\cH
_{a}\}
$ be an isonormal Gaussian noise with covariance $E[W_{a}(\varphi
)\times  W_{a}(\psi)]=\langle\varphi,\psi\rangle_{\cH_{a}}$ where
$\langle
\cdot, \cdot\rangle_{\cH_{a}}$ is given by
(\ref{def-cov}) with $f(x)$ replaced by $f_{a}(x)=|x|^{-a}$. Note that
$f=\cF\mu_{a}$ where $\mu_a(\xi)=(2\pi)^{-1}|\xi|^{a-1}\,d\xi$.
Let $u_{a}(t,x)$ be the solution of the equation
\[
\frac{\partial^2 u}{\partial t^2}=\Delta u+u \dot{W}_{a}\qquad (t>0,x \in\bR)
\]
with initial conditions $u(0,x)=u_0$ and $\frac{\partial u}{\partial
t}(0,x)=v_0$. This solution has the Wiener chaos expansion
$u_{a}(t,x)=\sum_{n \geq0}I_{n,a}(f_n(\cdot,t,x))$ where $I_{n,a}$
denotes the multiple Wiener integral with respect to $W_{a}$. Hence
\[
E\bigl|u_{a}(t,x)\bigr|^2=\sum_{n \geq0}
\frac{1}{n!} \alpha_{n,a}(t),
\]
where
%
\begin{equation}
\label{def-alpha-ve} \alpha_{n,a}(t) = \alpha_H^n
\int_{[0,t]^{2n}} \prod_{j=1}^{n}|t_j-s_j|^{2H-2}
\psi_{n,a}(\mathbf{t},\mathbf{s})\,d\mathbf{t} \,d\mathbf{s}
\end{equation}
and
%
\begin{eqnarray}\quad
&&\psi_{n,a}(\mathbf{t},\mathbf{s})
\nonumber
\\[-8pt]
\\[-8pt]
\nonumber
& &\qquad =\int_{\bR^{2n}} g_\mathbf{t}^{(n)}(x_1,
\ldots,x_n,t,x) g_\mathbf{s}^{(n)}(y_1,
\ldots,y_n,t,x) \prod_{j=1}^{n}f_{a}(x_j-y_j)\,d
\mathbf{x} \,d\mathbf{y}
\end{eqnarray}
and $g_\mathbf{t}^{(n)}(x_1, \ldots,x_n,t,x)=\prod_{j=1}^{n}G_\mathrm{
w}(t_{\rho(j+1)}-t_{\rho(j)},x_{\rho(j+1)}-x_{\rho(j)}) w(t_{\rho
(1)}, x_{\rho(1)})$ if $t_{\rho(1)}<\cdots<t_{\rho(n)}$.

\begin{lemma}
\label{conv-psi}
For any integer $n \geq1$ and for any $\mathbf{t},\mathbf{s} \in[0,t]^n$,
\[
\lim_{a \uparrow1} \psi_{n,a}(\mathbf{t},\mathbf{s})=\psi(
\mathbf{t},\mathbf{s}),
\]
where $\psi_n(\mathbf{t},\mathbf{s})$ is given by \eqref{def-psi1} with $d=1$
and $f=\delta_0$, that is,
\begin{eqnarray*}
\psi_n(\mathbf{t},\mathbf{s}) &= & \int_{\bR^n}
\prod_{j=1}^{n}G_\mathrm{
w}(t_{\rho(j+1)}-t_{\rho(j)},x_{\rho(j+1)}-x_{\rho(j)})w_\mathrm{
w}(t_{\rho
(1)},x_{\rho(1)})
\\
& &\quad{}\times \prod_{j=1}^{n}G_\mathrm{
w}(s_{\sigma(j+1)}-s_{\sigma(j)},x_{\sigma
(j+1)}-x_{\sigma(j)})w_\mathrm{
w}(s_{\sigma(1)},x_{\sigma(1)}) \,d\mathbf{x},
\end{eqnarray*}
where the permutations $\rho,\sigma\in S_n$ are chosen such that
$t_{\rho(1)}< \cdots<t_{\rho(n)}$ and $s_{\sigma(1)}< \cdots
<s_{\sigma
(n)}$, $t_{\rho(n+1)}=s_{\sigma(n+1)}=t$ and $x_{\rho
(n+1)}=x_{\sigma(n+1)}=x$.
\end{lemma}

\begin{pf} Note that for any $g,h \in L^1(\bR) \cap
L^2(\bR)$,
\begin{eqnarray*}
&&\lim_{a \uparrow1}\int_{\bR} \int
_{\bR} g(x)h(y) f_{a}(x-y)\,dx\,dy \\
&&\qquad= \lim
_{a \uparrow1}\frac{1}{2\pi} \int_{\bR} \cF g(
\xi) \overline{\cF h(\xi )}|\xi|^{a-1}\,d\xi\\
&&\qquad=
\frac{1}{2\pi} \int_{\bR}\cF g(\xi) \overline{\cF h(
\xi)}\,d\xi =\int_{\bR
} g(x)h(x)\,dx,
\end{eqnarray*}
by the dominated convergence theorem.
To justify the application of this theorem, we note that for $a>1/2$,
the integrand $|\cF g(\xi)| |\cF h(\xi)| |\xi|^{a-1}$ is bounded by the
integrable function
\[
\|g\|_1 \|h\|_1 |\xi|^{-1/2}
1_{\{|\xi| \leq1\}} +\bigl|\cF g(\xi)\bigr| \bigl|\cF h(\xi)\bigr| 1_{\{|\xi| \geq1\}}.
\]

From here we infer that
for any $g,h \in L^1(\bR^n) \cap L^2(\bR^n)$,
\[
\lim_{a \uparrow1}\int_{\bR^n} \int
_{\bR^n} g(\mathbf{x})h(\mathbf{y}) \prod
_{j=1}^{n}f_{a}(x_j-y_j)\,d
\mathbf{x} \,d\mathbf{y} = \int_{\bR^n} g(\mathbf{ x})h(
\mathbf{x})\,d\mathbf{x},
\]
with $\mathbf{x}=(x_1, \ldots,x_n)$ and $\mathbf{y}=(y_1, \ldots,y_n)$. We
apply this to $g=g_\mathbf{t}^{(n)}(\cdot,t,x)$ and $h=g_\mathbf{
s}^{(n)}(\cdot,t,x)$, using the fact that
\[
\psi_n(\mathbf{t},\mathbf{s})=\int_{\bR^n}g_\mathbf{t}^{(n)}(
\mathbf{x},t,x) g_\mathbf{s}^{(n)}(\mathbf{x},t,x)\,d\mathbf{x}.
\]
\upqed\end{pf}

\begin{lemma}
\label{conv-alpha}
For any $t>0$ and for any integer $n \geq1$,
\[
\lim_{a \uparrow1}\alpha_{n,a}(t)=\alpha_n(t).
\]
\end{lemma}

\begin{pf} This follows by Lemma~\ref{conv-psi} and the
dominated convergence theorem. It remains to justify the application of
this theorem. For this, we note that $\psi_{n,a}(\mathbf{t},\mathbf{s}) \leq
\psi_{n,a}(\mathbf{t},\mathbf{t})^{1/2} \psi_{n,a}(\mathbf{s},\mathbf{s})^{1/2}$.
Let $u_j=t_{\rho(j+1)}-t_{\rho(j)}$. As in the proof of Lemma~\ref
{estimate-psi}, it follows that for any $t \geq1$,
\begin{eqnarray*}
\psi_{n,a}(\mathbf{t},\mathbf{t}) & \leq& (u_0+t
v_0)^2 \frac{1}{(2\pi)^n} (4K_{a})^n(u_1
,\ldots, u_n)^{2-a}
\\
& \leq& (u_0+t v_0)^2 \frac{1}{(2\pi)^n}
(4K_{a})^n t^{n(1-a)} u_1 ,\ldots,
u_n
\\
& \leq& (u_0+t v_0)^2 \frac{t^n}{(2\pi)^n}
(4K_{a})^n u_1 ,\ldots, u_n,
\end{eqnarray*}
where $K_a:=K(\mu_{a})$ is given by (\ref{def-K-mu}). We now prove that
%
\begin{equation}
\label{K-equal-L} K_a=L_a:=\int_{\bR}
\frac{1}{1+|\xi|^2}\mu_a(d\xi).
\end{equation}

 To see this, note first that $L_a \leq K_a$. On the other
hand, for any $\eta\in\bR$,
\[
\int_{\bR} \frac{1}{1+|\xi-\eta|^2}\mu_{a}(d\xi)=\int
_{\bR
}e^{i\eta
x}p(x)|x|^{-a}\,dx,
\]
where $p(x)=(4\pi)^{-1/2}\int_0^{\infty}e^{-u} u^{-1/2}
e^{-|x|^2/(4u)}\,du$; see (3.4) of \cite{dalang-mueller03}. Taking the
modulus on both sides and using (3.5) of \cite{dalang-mueller03}, we
obtain that for any $\eta\in\bR$,
\[
\int_{\bR} \frac{1}{1+|\xi-\eta|^2}\mu_{a}(d\xi) \leq
\int_{\bR} p(x)|x|^{-a}\,dx=L_a.
\]
Taking the supremum over $\eta\in\bR$, we obtain that $K_a \leq L_a$.
This proves (\ref{K-equal-L}).

By considering separately the regions $\{|\xi| \leq1\}$ and $\{|\xi|
\geq1\}$, we see that $L_a \leq2(a^{-1}+(2-a)^{-1})$. Hence $L_a \leq
6$ if $a>1/2$.

Denote $\beta(\mathbf{t})=\prod_{j=1}^{n}(t_{\rho(j+1)}-t_{\rho(j)})$. It
follows that for any $a \in(1/2,1)$,
%
\begin{equation}
\label{estim-psi-ve} \psi_{n,a}(\mathbf{t},\mathbf{s}) \leq(u_0+t
v_0)^2 \frac{t^n}{(2\pi)^n} 24^n \bigl[\beta(
\mathbf{t}) \beta(\mathbf{s}) \bigr]^{1/2}.
\end{equation}
The claim is justified since $\int_{[0,t]^{2n}} \prod_{j=1}^{n}|t_j-s_j|^{2H-2} [\beta(\mathbf{t}) \beta(\mathbf{s})]^{1/2}\,d\mathbf{
t}\,d\mathbf{s}<\infty$; see the proof of Theorem~\ref{exist-sol}.
\end{pf}

\begin{lemma}
\label{conv-u}
For any $t>0$ and for any $x \in\bR^d$,
\[
\lim_{a \uparrow1}E\bigl|u_{a}(t,x)\bigr|^2=E\bigl|u(t,x)\bigr|^2.
\]
\end{lemma}

\begin{pf} The result follows by Lemma~\ref{conv-alpha} and
the dominated convergence theorem. We justify the application of this
theorem. By (\ref{def-alpha-ve}) and (\ref{estim-psi-ve}),
\begin{eqnarray*}
\alpha_{n,a}(t) & \leq& (u_0+t v_0)^2
\frac{t^n}{(2\pi)^n} 24^n \int_{[0,t]^{2n}} \prod
_{j=1}^{n}|t_j-s_j|^{2H-2}
\bigl[\beta(\mathbf{t}) \beta (\mathbf{s}) \bigr]^{1/2}\,d\mathbf{t}\,d
\mathbf{s}\,d\mathbf{t}\,d\mathbf{s}
\\
& \leq& (u_0+t v_0)^2 c^n
\frac{1}{n!} t^{(2H+1)n},
\end{eqnarray*}
for any $a \in(1/2,1)$, where the last inequality follows as in the
proof of Theorem~\ref{exist-sol}. Since $\sum_{n} c^n
t^{(2H+1)n}/(n!)^2<\infty$, the proof is complete.
\end{pf}

The next result corresponds to Theorem~\ref{wave-main}(c), in case (iv).

\begin{theorem}
\label{LB-th-f-white}
Let $f$ be the kernel of case \textup{(iv)}. Then, for any $x \in\bR$ and for
any $t >0$, we have
\[
E\bigl|u(t,x)\bigr|^2 \geq c_3 u_0^2\exp
\bigl(c_4 t^{\rho_\mathrm{ w}} \bigr),
\]
where $c_3>0$ and $c_4>0$ are some constants depending on $H$, and
$\rho
_\mathrm{ w}$ is given by~(\ref{def-rho}).
\end{theorem}

\begin{pf} By Theorem~\ref{LB-th-f-Riesz}, for any $x \in
\bR^d$ and for any $t \geq t_a$,
%
\begin{equation}
\label{ineq-u-a} E\bigl|u_a(t,x)\bigr|^2 \geq u_0^2
\exp \bigl(c_{a} t^{(2H+2-a)/(3-a)} \bigr),
\end{equation}
where the constants $c_a>0$ and $t_a>0$ are given by (\ref{def-c2}) and
(\ref{def-t2}), that is,
\[
c_{a}=\tfrac{1}{2} \bigl(e^{-1} c_H
4^{-a} \bigr)^{1/(3-a)} \quad\mbox{and}\quad t_a= \bigl(e
c_H ^{-1}2^{3+a} \bigr)^{1/(2H+2-a)}.
\]
Let
$c_4=\lim_{a \uparrow1}c_a$ and $t_0'=\lim_{a \uparrow1}t_a$. Then
$t_a \leq2t_0'=:t_0$ for all $a \in(a_0,1)$. Fix $t \geq t_0$.
We let $a \uparrow1$ in (\ref{ineq-u-a}). Using Lemma~\ref{conv-u}, we
infer that
\[
E\bigl|u(t,x)\bigr|^2 \geq u_0^2\exp
\bigl(c_4 t^{\rho_\mathrm{ w}} \bigr)\qquad \mbox{for all } t \geq
t_0.
\]
For $0<t<t_0$, we argue as in step 9 of the proof of Theorem~\ref
{LB-th-f-bounded}.
\end{pf}


Summarizing the results of this section, we can say that Theorems \ref
{LB-th-f-bounded} and \ref{LB-th-f-Riesz} generalize Theorem~4.1 of \cite{dalang-mueller09} (with $p=2$) to the case of the
fractional noise in time. However, reference \cite{dalang-mueller09}
does not contain a result analogous to Theorem~\ref{LB-th-f-white} for
the case $H=1/2$, that is, when $W$ is a space--time white noise.

\section{Parabolic case: Proof of Theorem \texorpdfstring{\protect\ref{heat-main}}{2.2}} \label{sec:PAM}

In this section, we examine equation \eqref{heat}. We state and sketch
the proof of two results, which together give the conclusion of Theorem~\ref{heat-main}. The proofs are similar to those presented above in the
hyperbolic case, taking $G = G_\mathrm{ h}$ and $w = w_\mathrm{ h}$. For the
lower bound, we use a FK representation similar to the one given in
\cite{B09}, except that here we work with processes $X^1,X^2$ defined
by \eqref{def-X1-heat} and \eqref{def-X2-heat} below, instead of
Brownian motions $B^1,B^2$.

We define a different constant $K_\mathrm{ h}$ than in the hyperbolic case, namely
%
\begin{equation}
\label{def-K-heat} K_\mathrm{ h} =\cases{ %
\mu
\bigl( \bR^d \bigr) ,& \quad $\mbox{in case (i),}$
\vspace*{2pt}\cr
K(\mu), & \quad$\mbox{in cases (ii) and (iii),}$
\vspace*{2pt}\cr
\sqrt{\pi}, &\quad $\mbox{in case (iv).}$}
\end{equation}
%
%
We recall that in the parabolic case, the spatial dimension $d \geq1$
is arbitrary. The first result gives the existence of the solution and
an upper bound for its moments of order $p \geq2$.

\begin{proposition}
\label{exist-sol-heat}
Let $f$ be a kernel of cases \textup{(i)--(iv),} and $\rho_\mathrm{ h},a,K_\mathrm{ h}$
be the constants given by (\ref{def-rho}), (\ref{def-a}), respectively
(\ref{def-K-heat}). Assume that (\ref{cond-mu}) holds. Then:
\begin{longlist}[(a)]
\item[(a)] for any $t>0$ and for any integer $n \geq1$,
\[
\alpha_n(t) \leq u_0^2 K_\mathrm{
h}^n c^n (n!)^{a/2} t^{(4H-a)n/2},
\]
where $c$ is a constant depending on $H$ and $a$;

\item[(b)] equation (\ref{heat}) has a unique solution $u(t,x)$ which
has the following property: for any $p \geq2$, for any $x \in\bR^d$
and for any $t>0$, 
\[
E\bigl|u(t,x)\bigr|^p \leq c_1^p u_0^p
\exp \bigl(c_2 K_\mathrm{ h}^{2/(2-a)} p^{(4-a)/(2-a)}
t^{\rho_\mathrm{ h}} \bigr),
\]
where $c_1>0$ is a constant depending on $a$, and $c_2>0$ is a constant
depending on $H$ and $a$.
\end{longlist}
\end{proposition}

\begin{pf}(a) Similar to Lemma~\ref{estimate-psi}, it
can be shown that
\begin{eqnarray*}
&&\psi_n(\mathbf{t},\mathbf{t})
\\
&&\qquad= u_0^2 \int_{\bR^{nd}} \exp
\bigl(-u_1|\xi_1|^2 \bigr) \cdots\exp
\bigl(-u_n|\xi _1+\cdots+\xi_n|^2
\bigr) \mu(d\xi_1) \cdots\mu(d\xi_n)
\\
&&\qquad \leq u_0^2 K_\mathrm{ h}^n
(u_1 ,\ldots, u_n)^{-a/2}.
\end{eqnarray*}
To prove this in cases (ii) and (iii), one uses the following inequality:
\[
\int_{\bR^d} \exp \bigl(-t|\xi-\eta|^2 \bigr)
\mu(d \xi) \leq K(\mu) t^{-a/2}.
\]
The conclusion follows as in the proof of Proposition~\ref
{exist-sol}(a). Note that
\[
\psi_{n}(\mathbf{t},\mathbf{s}) \leq u_0^2
K_\mathrm{ h}^n \bigl[\beta(\mathbf{t}) \beta (\mathbf{s})
\bigr]^{-a/4}.
\]

(b) The conclusion follows as in the proof of Proposition~\ref
{exist-sol}(b) (case $p=2$), respectively Proposition~\ref{UB-theorem}
(case $p>2$).
\end{pf}

For the lower bound, we use the following representation for the second
moment of the solution to \eqref{heat}, which can be obtained as in
Section~\ref{FK-section}, assuming that $f$ is a function: 
\begin{eqnarray*}
\hspace*{-6pt}&&E\bigl|u(t,x)\bigr|^2
\\
\hspace*{-6pt}&&\qquad= e^{t^2} u_0^2 \sum
_{n \geq0}\mathop{ \sum_{i_1,\ldots, i_n}}_{\mathrm{distinct}} E_x \Biggl[\prod
_{j=1}^{n}f \bigl(X_{T_{i_j}}^{1}-X_{S_{i_j}}^{2}
\bigr) \alpha_H^n \prod_{j=1}^{n}|T_{i_j}-S_{i_j}|^{2H-2}1_{B_{i_1 ,\ldots, i_n}(t)}
\Biggr].
\end{eqnarray*}
Here, the event $B_{i_1, \ldots,i_n}(t)$ and the points
$(T_{i_j},S_{i_j})$ are defined as in Section~\ref{FK-section}, but the
processes $X^1$ and $X^2$ are given by
%
\begin{eqnarray}
\label{def-X1-heat} X_{s}^1 &= &X_{\tau_i}^1+
\sqrt{s-\tau_{i}}\Theta_{i+1}^1\qquad \mbox {if }
\tau_i \leq s \leq\tau_{i+1},
\\
\label{def-X2-heat} X_{s}^2 &=& X_{\tau_i'}^2+
\sqrt{s-\tau_{i}'}\Theta_{i+1}^2\qquad
\mbox {if } \tau_i' \leq s \leq\tau_{i+1}',
\end{eqnarray}
where $(\Theta_i^1)_{i \geq1}$ and $(\Theta_{i}^2)_{i \geq1}$ are two
independent collections of i.i.d. random variables with values in $\bR
^d$ with the same law as $\Theta_0$, and $\Theta_0$ has a
$d$-dimensional standard normal distribution.
Note that in this case,
%
\begin{equation}
\label{distr-heat} G_\mathrm{ h}(t,\cdot) \qquad\mbox{is the density of } \sqrt{t}
\Theta_0.
\end{equation}
(Alternatively, $X^1,X^2$ can be two independent $d$-dimensional
standard Brownian motions; see Remark~\ref{cond-density-Y} and \cite{B09}.)

\begin{proposition}
Let $f$ be a kernel of cases \textup{(i)--(iv)}, and $\rho_\mathrm{ h}$ be the
constant given by (\ref{def-rho}). Then for any $x \in\bR^d$ and for
any $t>0$,
\[
E\bigl|u(t,x)\bigr|^2 \geq c_3 u_0^2
\exp \bigl(c_4 t^{\rho} \bigr),
\]
where $c_3>0$ and $c_4>0$ are some constants depending on $H$ and $a$.
\end{proposition}

\begin{pf} In case (i), the argument is similar to the one
used in Theorem~\ref{LB-th-f-bounded}. One difference is that we
replace $\delta$ by $\delta^2$. This is essentially due to the use of a
parabolic rather than hyperbolic scaling; compare \eqref{distr-heat}
with \eqref{density-Theta}. In addition, in the events $D^1(t),D^2(t)$,
we add the condition $|\Theta_{j+1}^1| \leq1$, respectively $|\Theta
_{j+1}^2| \leq1$, for all $j=1, \ldots,k-1$. Note that the variable
$\widetilde{Z}_t$ (in step 5) is replaced by $1$. Instead of (\ref
{final-LB}), we obtain that for all $k \geq k_0$,
\[
E\bigl|u(t,x)\bigr|^2 \geq u_0^2 \biggl(
\alpha_0 c_H \frac{t^{2H}}{k} \biggr)^{k}.
\]
The argument for cases (ii)--(iii) is similar to the proof of Theorem~\ref{LB-th-f-Riesz}, leading to the following lower bound: there exists
$k_0>0$ such that for all $k \geq k_0$,
\[
E\bigl|u(t,x)\bigr|^2 \geq u_0^2 \biggl(
c_H^* \frac
{t^{2H-a/2}}{k^{1-a/2}} \biggr)^k.
\]
Choosing $k$ appropriately completes the proof. The argument for case
(iv) is similar to the proof of Theorem~\ref{LB-th-f-white}. In all
cases, the argument is extended to all $t>0$, as in step 9 of the proof
of Theorem~\ref{LB-th-f-bounded}.
\end{pf}

\begin{appendix}
\section{An elementary result}
\label{sec:app}

\begin{lemma}
\label{series-ineq}
For any $a>0$, we have
\[
\sum_{n \geq0}\frac{x^n}{(n!)^a} \leq c_1
\exp \bigl(c_2 x^{1/a} \bigr) \qquad\mbox {for all } x >0,
\]
where $c_1>0$ and $c_2>0$ are some constants depending on $a$.
\end{lemma}

\begin{pf} Note that for any $a>0$, we have
%
\begin{equation}
\label{stirling} \lim_{n \to\infty}\frac{\Gamma(an+1)}{(n!)^aa^{an+1/2}(2\pi
n)^{(1-a)/2}}=1;
\end{equation}
see also (3.19) of \cite{hu_nualart}. To see this, we use Stirling's
formula in the following format:
\[
\lim_{x \to\infty}\frac{\Gamma(x+1)}{x^{x} e^{-x} \sqrt{2\pi x}}=1
\]
(see, e.g., Corollary~3 of \cite{li06}), from which we infer that
\[
\Gamma(an+1) \sim(an)^{an} e^{-an} (2\pi an)^{1/2}
\quad\mbox{and}\quad n!=\Gamma(n+1) \sim n^{n} e^{-n} (2\pi
n)^{1/2}.
\]
Here we use the notation $a_n \sim b_n$ to indicate that $a_n/b_n \to
1$ as $n \to\infty$.

From \eqref{stirling}, it follows that there exists a constant
$C_1>0$ depending on $a$, such that
\[
\frac{\Gamma(an+1)}{C_n (n!)^a} \leq C_1\qquad \mbox{for all } n \geq0,
\]
where $C_n=a^{an}n^{(1-a)/2}$. Clearly, we can choose a constant
$\lambda>1$ (depending on~$a$) such that $\lambda^{-n} \leq C_n \leq
\lambda^n$ for all $n$. Therefore,
%
\begin{equation}
\label{lemmaA-ineq1} \sum_{n \geq0}\frac{x^n}{(n!)^{a}} \leq
C_1 E_{a}(\lambda x),
\end{equation}
where $E_a(x)=\sum_{n \geq0}x^n/\Gamma(an+1)$ denotes the
Mittag--Leffler function.

We now use the asymptotic behavior of $E_a(x)$ for $x>0$:
\[
\lim_{x \to\infty}\frac{E_a(x)}{\exp(x^{1/a})}=\frac{1}{a} \qquad\mbox {for
all } a>0
\]
(see Theorem~1 of \cite{gerhold12}). Hence there exist some constants
$C_2>0$ and $x_0>0$ depending on $a$ such that
\[
E_a(x) \leq C_2 \exp \bigl(x^{1/a} \bigr)
\qquad\mbox{for all } x \geq x_0.
\]
If $0<x<x_0$, then $x^n \leq x_0^n$ and $E_{a}(x) \leq C_3 \leq C_3
\exp
(x^{1/a})$, where
$C_3=E_a(x_0)$ depends only on $a$. Taking $C_4=\max(C_2,C_3)$, it
follows that
%
\begin{equation}
\label{lemmaA-ineq2} E_{a}(x) \leq C_4 \exp
\bigl(x^{1/a} \bigr) \qquad\mbox{for all } x>0.
\end{equation}
The conclusion follows from \eqref{lemmaA-ineq1} and \eqref
{lemmaA-ineq2}.
\end{pf}

\section{A fundamental inequality}
\label{sec:appB}
In this section, we prove inequality \eqref{MMV-ineq} which is used in
the proof of Proposition~\ref{exist-sol}. Note that this inequality is
a simplified form of (2.5) of \cite{hu_nualart}.

We first recall the Hardy--Littlewood--Sobolev theorem.

\begin{theorem}[(Theorem~1, page 119 of \cite{stein70})]
Let $0<\alpha<n$ and $1<p<\infty$. Let $q>p$ be such that
$1/p-1/q=\alpha/n$. For any $\varphi\in L^p(\bR^n)$, the integral
\[
(I_{\alpha} \varphi) (x):=\int_{\bR^n}
\varphi(y)|x-y|^{-n+\alpha}\,dy
\]
converges absolutely for almost all $x \in\bR^n$, and
%
\begin{equation}
\label{HLS} \|I_{\alpha} \varphi\|_{L^q(\bR^n)} \leq
C_{n,\alpha,p}\|\varphi\| _{L^p(\bR^n)},
\end{equation}
where $C_{n,\alpha,p}>0$ is a constant depending on $n,\alpha$ and $p$.
\end{theorem}

The following inequality is due to \cite{MMV01}. We include its proof
for the sake of completeness.

\begin{lemma}
\label{lemmaB2}
Let $H \in(1/2,1)$ and $\alpha_H=H(2H-1)$. For any $f,g \in
L^{1/H}(\bR)$,
%
\begin{eqnarray}
\label{MMV-ineq-1} &&\int_{\bR} \int_{\bR}\bigl|f(t)\bigr|
\bigl|g(s)\bigr||t-s|^{2H-2}\,dt\,ds
\nonumber
\\[-8pt]
\\[-8pt]
\nonumber
&&\qquad
\leq C_H \biggl(\int_{\bR}\bigl|f(t)\bigr|^{1/H}\,dt
\biggr)^{H} \biggl(\int_{\bR
}\bigl|g(t)\bigr|^{1/H}\,dt
\biggr)^{H},
\end{eqnarray}
where $C_H>0$ is the constant from \eqref{HLS} with $n=1, \alpha=2H-1$
and $p=1/H$.
\end{lemma}

\begin{pf} Using H\"older's inequality with $p=1/H$ and
$q=1/(1-H)$, we infer that the left-hand side of \eqref{MMV-ineq-1} is
smaller than
\begin{eqnarray*}
& & \biggl(\int_{\bR}\bigl|f(t)\bigr|^{1/H}\,dt
\biggr)^{H} \biggl[\int_{\bR} \biggl( \int
_{\bR}|g(s)||t-s|^{2H-2}\,ds \biggr)^{1/(1-H)} \,dt
\biggr]^{1-H}
\\
& &\qquad = \|f\|_{L^{1/H}(\bR)} \cdot \biggl\{\int_{\bR}
\bigl[\bigl(I_{2H-1} |g|\bigr) (t) \bigr]^{1/(1-H)} \,dt \biggr
\}^{1-H}
\\
& & \qquad = \|f\|_{L^{1/H}(\bR)}\bigl \|I_{2H-1}|g|\bigr\|_{L^{1/(1-H)}(\bR)}.
\end{eqnarray*}
The conclusion now follows by \eqref{HLS} with $n=1, \alpha=2H-1$,
$p=1/H$ and $q=1/(1-H)$.
\end{pf}

\begin{lemma}
\label{lemmaB3}
For any $\varphi\in L^{1/H}(\bR^n)$,
%
\begin{equation}
\label{ineq} \int_{\bR^n} \int_{\bR^n}
\varphi(\mathbf{t}) \varphi(\mathbf{s})\prod_{i=1}^{n}|t_i-s_i|^{2H-2}\,d
\mathbf{t} \,d\mathbf{s} \leq C_H^n \biggl(\int
_{\bR
^n}\bigl|\varphi(\mathbf{t})\bigr|^{1/H}\,d\mathbf{t}
\biggr)^{2H},
\end{equation}
where $C_H>0$ is the constant from Lemma~\ref{lemmaB2}, and we denote
$\mathbf{t}=(t_1, \ldots,t_n)$ and $\mathbf{s}=(s_1, \ldots,s_n)$.
\end{lemma}

\begin{pf} We proceed by induction on $n$. For $n=1$, the
result holds by Lemma~\ref{lemmaB2}. Suppose that \eqref{ineq} holds
for $n-1$. By applying Lemma~\ref{lemmaB2} to the functions $f(\cdot
)=\varphi(t_1,\ldots,t_{n-1},\cdot)$ and $g(\cdot)=\varphi
(s_1,\ldots
,s_{n-1},\cdot)$ for fixed $(t_1, \ldots,t_{n-1}) \in\bR^{n-1}$ and
$(s_1, \ldots,s_{n-1}) \in\bR^{n-1}$, we obtain that
\begin{eqnarray*}
&&\int_{\bR} \int_{\bR}\bigl |
\varphi(t_1,\ldots,t_{n-1},t_n)\bigr|\bigl|\varphi
(s_1,\ldots,s_{n-1},s_n)\bigr||t_n-s_n|^{2H-2}\,dt_n
\,ds_n \\
&&\qquad\leq
C_H \bigl\|\varphi(t_1,\ldots,t_{n-1},\cdot)
\bigr\|_{1/H} \bigl\|\varphi (s_1,\ldots ,s_{n-1},\cdot)
\bigr\|_{1/H},
\end{eqnarray*}
where $\|\cdot\|_{1/H}$ denotes the $L^{1/H}(\bR)$-norm. [By Fubini's
theorem, the functions $f$ and $g$ are in $L^{1/H}(\bR)$ for almost all
$(t_1,\ldots,t_{n-1}) \in\bR^{n-1}$ and $(s_1,\ldots,\break s_{n-1}) \in
\bR^{n-1}$.]
Hence, the left-hand side of \eqref{ineq} is less that
%
\begin{eqnarray}
\label{ineq1} &&C_H \int_{\bR}\int
_{\bR}\bigl\|\varphi(t_1,\ldots,t_{n-1},
\cdot)\bigr\| _{1/H} \bigl\| \varphi(s_1,\ldots,s_{n-1},
\cdot)\bigr\|_{1/H}
\nonumber
\\[-8pt]
\\[-8pt]
\nonumber
&&\hspace*{18pt}\qquad{}\times \prod_{i=1}^{n-1}|t_i-s_i|^{2H-2}\,d
\mathbf{t}_{n-1} \,d\mathbf{s}_{n-1},
\end{eqnarray}
where $\mathbf{t}_{n-1}=(t_1, \ldots,t_{n-1})$ and $\mathbf{s}_{n-1}=(s_1,
\ldots,s_{n-1})$.

By the induction hypothesis, \eqref{ineq1} is less than
\begin{eqnarray*}
&&C_H^{n} \biggl(\int_{\bR^{n-1}}\bigl\|
\varphi(t_1,\ldots,t_{n-1},\cdot )\bigr\| _{1/H}^{1/H}
\,dt_1,\ldots,dt_{n-1} \biggr)^{2H}
\nonumber
\\[-8pt]
\\[-8pt]
\nonumber
&&\qquad=C_H^n
\biggl(\int_{\bR
}\int_{\bR}\bigl |\varphi(
\mathbf{t})\bigr|^{1/H}\,d\mathbf{t} \biggr)^{2H},
\end{eqnarray*}
where $\mathbf{t}=(t_1,\ldots,t_{n})$. This proves \eqref{ineq}.
\end{pf}
\end{appendix}

\section*{Acknowledgments}
The authors would like to thank Jian Song (and co-authors) for sending
them the preprint of \cite{CHSX}. The authors are grateful to two
anonymous referees who read the paper very carefully and made numerous
suggestions for improving the presentation.






\printaddresses

\begin{thebibliography}{50}

\bibitem{alos-nualart03}
\begin{barticle}[mr]
\bauthor{\bsnm{Al{\`o}s},~\bfnm{Elisa}\binits{E.}} \AND
\bauthor{\bsnm{Nualart},~\bfnm{David}\binits{D.}}
(\byear{2003}).
\btitle{Stochastic integration with respect to the fractional {B}rownian motion}.
\bjournal{Stoch. Stoch. Rep.}
\bvolume{75}
\bpages{129--152}.
\bid{doi={10.1080/1045112031000078917}, issn={1045-1129}, mr={1978896}}
\end{barticle}
%
\bptok{imsref}%
\endbibitem

\bibitem{B09}
\begin{barticle}[mr]
\bauthor{\bsnm{Balan},~\bfnm{Raluca M.}\binits{R.~M.}}
(\byear{2009}).
\btitle{A note on a {F}enyman--{K}ac-type formula}.
\bjournal{Electron. Commun. Probab.}
\bvolume{14}
\bpages{252--260}.
\bid{doi={10.1214/ECP.v14-1468}, issn={1083-589X}, mr={2516260}}
\end{barticle}
%
\bptok{imsref}%
\endbibitem

\bibitem{B12}
\begin{barticle}[mr]
\bauthor{\bsnm{Balan},~\bfnm{Raluca~M.}\binits{R.~M.}}
(\byear{2012}).
\btitle{The stochastic wave equation with multiplicative fractional noise: A~{M}alliavin calculus approach}.
\bjournal{Potential Anal.}
\bvolume{36}
\bpages{1--34}.
\bid{doi={10.1007/s11118-011-9219-z}, issn={0926-2601}, mr={2886452}}
\end{barticle}
%
\bptok{imsref}%
\endbibitem

\bibitem{B12-Fourier}
\begin{barticle}[mr]
\bauthor{\bsnm{Balan},~\bfnm{Raluca~M.}\binits{R.~M.}}
(\byear{2012}).
\btitle{Linear {SPDE}s driven by stationary random distributions}.
\bjournal{J. Fourier Anal. Appl.}
\bvolume{18}
\bpages{1113--1145}.
\bid{doi={10.1007/s00041-012-9240-7}, issn={1069-5869}, mr={3000977}}
\end{barticle}
%
\bptok{imsref}%
\endbibitem

\bibitem{BC13}
\begin{barticle}[mr]
\bauthor{\bsnm{Balan},~\bfnm{Raluca~M.}\binits{R.~M.}} \AND
\bauthor{\bsnm{Conus},~\bfnm{Daniel}\binits{D.}}
(\byear{2014}).
\btitle{A note on intermittency for the fractional heat equation}.
\bjournal{Statist. Probab. Lett.}
\bvolume{95}
\bpages{6--14}.
\bid{doi={10.1016/j.spl.2014.08.001}, issn={0167-7152}, mr={3262944}}
\end{barticle}
%
\bptok{imsref}%
\endbibitem

\bibitem{BT10}
\begin{barticle}[mr]
\bauthor{\bsnm{Balan},~\bfnm{Raluca~M.}\binits{R.~M.}} \AND
\bauthor{\bsnm{Tudor},~\bfnm{Ciprian~A.}\binits{C.~A.}}
(\byear{2010}).
\btitle{Stochastic heat equation with multiplicative fractional-colored noise}.
\bjournal{J. Theoret. Probab.}
\bvolume{23}
\bpages{834--870}.
\bid{doi={10.1007/s10959-009-0237-3}, issn={0894-9840}, mr={2679959}}
\end{barticle}
%
\bptok{imsref}%
\endbibitem

\bibitem{BT10-SPA}
\begin{barticle}[mr]
\bauthor{\bsnm{Balan},~\bfnm{Raluca~M.}\binits{R.~M.}} \AND
\bauthor{\bsnm{Tudor},~\bfnm{Ciprian~A.}\binits{C.~A.}}
(\byear{2010}).
\btitle{The stochastic wave equation with fractional noise: A random field approach}.
\bjournal{Stochastic Process. Appl.}
\bvolume{120}
\bpages{2468--2494}.
\bid{doi={10.1016/j.spa.2010.08.006}, issn={0304-4149}, mr={2728174}}
\end{barticle}
%
\bptok{imsref}%
\endbibitem

\bibitem{BQS}
\begin{barticle}[mr]
\bauthor{\bsnm{Bal{\'a}zs},~\bfnm{M.}\binits{M.}},
\bauthor{\bsnm{Quastel},~\bfnm{J.}\binits{J.}} \AND
\bauthor{\bsnm{Sepp{\"a}l{\"a}inen},~\bfnm{T.}\binits{T.}}
(\byear{2011}).
\btitle{Fluctuation exponent of the KPZ/stochastic {B}urgers equation}.
\bjournal{J. Amer. Math. Soc.}
\bvolume{24}
\bpages{683--708}.
\bid{doi={10.1090/S0894-0347-2011-00692-9}, issn={0894-0347}, mr={2784327}}
\end{barticle}
%
\bptok{imsref}%
\endbibitem

\bibitem{bertini_cancrini}
\begin{barticle}[mr]
\bauthor{\bsnm{Bertini},~\bfnm{Lorenzo}\binits{L.}} \AND
\bauthor{\bsnm{Cancrini},~\bfnm{Nicoletta}\binits{N.}}
(\byear{1995}).
\btitle{The stochastic heat equation: {F}eynman--{K}ac formula and intermittence}.
\bjournal{J. Stat. Phys.}
\bvolume{78}
\bpages{1377--1401}.
\bid{doi={10.1007/BF02180136}, issn={0022-4715}, mr={1316109}}
\end{barticle}
%
\bptok{imsref}%
\endbibitem

\bibitem{caithamer05}
\begin{barticle}[mr]
\bauthor{\bsnm{Caithamer},~\bfnm{Peter}\binits{P.}}
(\byear{2005}).
\btitle{The stochastic wave equation driven by fractional {B}rownian noise and temporally correlated smooth noise}.
\bjournal{Stoch. Dyn.}
\bvolume{5}
\bpages{45--64}.
\bid{doi={10.1142/S0219493705001286}, issn={0219-4937}, mr={2118754}}
\end{barticle}
%
\bptok{imsref}%
\endbibitem

\bibitem{carmona_molchanov}
\begin{barticle}[mr]
\bauthor{\bsnm{Carmona},~\bfnm{Ren{\'e}~A.}\binits{R.~A.}} \AND
\bauthor{\bsnm{Molchanov},~\bfnm{S.~A.}\binits{S.~A.}}
(\byear{1994}).
\btitle{Parabolic {A}nderson problem and intermittency}.
\bjournal{Mem. Amer. Math. Soc.}
\bvolume{108}
\bpages{viii+125}.
\bid{doi={10.1090/memo/0518}, issn={0065-9266}, mr={1185878}}
\bptnote{check volume, check pages}%
\end{barticle}
%
\bptok{imsref}%
\endbibitem

\bibitem{Chen_Dalang}
\begin{bmisc}[mr]
\bauthor{\bsnm{Chen},~\bfnm{Le}\binits{L.}} \AND
\bauthor{\bsnm{Dalang},~\bfnm{Robert~C.}\binits{R.~C.}}
(\byear{2015}).
\bhowpublished{Moments and
growth indices for the nonlinear stochastic heat equation with rough
initial conditions. \textit{Ann. Probab.}  \textbf{43} 3006--3051.}
\bid{mr={3433576}}
\end{bmisc}
%
\bptok{imsref}%
\endbibitem

\bibitem{chen14}
\begin{barticle}[auto:parserefs-M02]
\bauthor{\bsnm{Chen},~\bfnm{X.}\binits{X.}}
(\byear{2016}).
\btitle{Spatial asymptotics for the parabolic Anderson models with generalized times--space Gaussian noise}.
\bjournal{Ann. Probab.}
\bnote{To appear}.
\end{barticle}
%
\bptok{imsref}%
\endbibitem

\bibitem{CHS14}
\begin{bmisc}[auto:parserefs-M02]
\bauthor{\bsnm{Chen},~\bfnm{X.}\binits{X.}},
\bauthor{\bsnm{Hu},~\bfnm{Y.}\binits{Y.}} \AND
\bauthor{\bsnm{Song},~\bfnm{J.}\binits{J.}}
(\byear{2014}).
\bhowpublished{Feynman--Kac formula for fractional heat equation driven by fractional white noise.
Preprint. Available at \arxivurl{arXiv:1203.0477}.}
\end{bmisc}
%
\bptok{imsref}%
\endbibitem

\bibitem{CHSX}
\begin{barticle}[auto:parserefs-M02]
\bauthor{\bsnm{Chen},~\bfnm{X.}\binits{X.}},
\bauthor{\bsnm{Hu},~\bfnm{Y.}\binits{Y.}},
\bauthor{\bsnm{Song},~\bfnm{J.}\binits{J.}} \AND
\bauthor{\bsnm{Xing},~\bfnm{F.}\binits{F.}}
(\byear{2016}).
\btitle{Exponential asymptotics for time--space Hamiltonians}.
\bjournal{Ann. Inst. Henri Poincar\'{e} Probab. Stat.}
\bvolume{51}
\bpages{1529--1561}.
\bid{mr={3414457}}
\end{barticle}
%
\bptok{imsref}%
\endbibitem

\bibitem{CD08}
\begin{barticle}[mr]
\bauthor{\bsnm{Conus},~\bfnm{Daniel}\binits{D.}} \AND
\bauthor{\bsnm{Dalang},~\bfnm{Robert~C.}\binits{R.~C.}}
(\byear{2008}).
\btitle{The non-linear stochastic wave equation in high dimensions}.
\bjournal{Electron. J. Probab.}
\bvolume{13}
\bpages{629--670}.
\bid{doi={10.1214/EJP.v13-500}, issn={1083-6489}, mr={2399293}}
\end{barticle}
%
\bptok{imsref}%
\endbibitem

\bibitem{CJK}
\begin{barticle}[mr]
\bauthor{\bsnm{Conus},~\bfnm{Daniel}\binits{D.}},
\bauthor{\bsnm{Joseph},~\bfnm{Mathew}\binits{M.}} \AND
\bauthor{\bsnm{Khoshnevisan},~\bfnm{Davar}\binits{D.}}
(\byear{2013}).
\btitle{On the chaotic character of the stochastic heat equation, before the onset of intermitttency}.
\bjournal{Ann. Probab.}
\bvolume{41}
\bpages{2225--2260}.
\bid{doi={10.1214/11-AOP717}, issn={0091-1798}, mr={3098071}}
\end{barticle}
%
\bptok{imsref}%
\endbibitem

\bibitem{CJKS}
\begin{bincollection}[mr]
\bauthor{\bsnm{Conus},~\bfnm{Daniel}\binits{D.}},
\bauthor{\bsnm{Joseph},~\bfnm{Mathew}\binits{M.}},
\bauthor{\bsnm{Khoshnevisan},~\bfnm{Davar}\binits{D.}} \AND
\bauthor{\bsnm{Shiu},~\bfnm{Shang-Yuan}\binits{S.-Y.}}
(\byear{2013}).
\btitle{Intermittency and chaos for a nonlinear stochastic wave equation in dimension 1}.
In \bbooktitle{Malliavin Calculus and Stochastic Analysis}.
\bseries{Springer Proc. Math. Stat.}
\bvolume{34}
\bpages{251--279}.
\bpublisher{Springer},
\blocation{New York}.
\bid{doi={10.1007/978-1-4614-5906-4_11}, mr={3070447}}
\end{bincollection}
%
\bptok{imsref}%
\endbibitem

\bibitem{CJKS2}
\begin{barticle}[mr]
\bauthor{\bsnm{Conus},~\bfnm{Daniel}\binits{D.}},
\bauthor{\bsnm{Joseph},~\bfnm{Mathew}\binits{M.}},
\bauthor{\bsnm{Khoshnevisan},~\bfnm{Davar}\binits{D.}} \AND
\bauthor{\bsnm{Shiu},~\bfnm{Shang-Yuan}\binits{S.-Y.}}
(\byear{2013}).
\btitle{On the chaotic character of the stochastic heat equation, {II}}.
\bjournal{Probab. Theory Related Fields}
\bvolume{156}
\bpages{483--533}.
\bid{doi={10.1007/s00440-012-0434-3}, issn={0178-8051}, mr={3078278}}
\end{barticle}
%
\bptok{imsref}%
\endbibitem

\bibitem{CK12}
\begin{barticle}[mr]
\bauthor{\bsnm{Conus},~\bfnm{Daniel}\binits{D.}} \AND
\bauthor{\bsnm{Khoshnevisan},~\bfnm{Davar}\binits{D.}}
(\byear{2012}).
\btitle{On the existence and position of the farthest peaks of a family of stochastic heat and wave equations}.
\bjournal{Probab. Theory Related Fields}
\bvolume{152}
\bpages{681--701}.
\bid{doi={10.1007/s00440-010-0333-4}, issn={0178-8051}, mr={2892959}}
\end{barticle}
%
\bptok{imsref}%
\endbibitem

\bibitem{dalang}
\begin{barticle}[mr]
\bauthor{\bsnm{Dalang},~\bfnm{Robert~C.}\binits{R.~C.}}
(\byear{1999}).
\btitle{Extending the martingale measure stochastic integral with applications to spatially homogeneous S.P.D.E.'s}.
\bjournal{Electron. J. Probab.}
\bvolume{4}
\bpages{29 pp. (electronic)}.
\bid{doi={10.1214/EJP.v4-43}, issn={1083-6489}, mr={1684157}}
\bptnote{check pages}%
\end{barticle}
%
\bptok{imsref}%
\endbibitem

\bibitem{dalang-frangos98}
\begin{barticle}[mr]
\bauthor{\bsnm{Dalang},~\bfnm{Robert~C.}\binits{R.~C.}} \AND
\bauthor{\bsnm{Frangos},~\bfnm{N.~E.}\binits{N.~E.}}
(\byear{1998}).
\btitle{The stochastic wave equation in two spatial dimensions}.
\bjournal{Ann. Probab.}
\bvolume{26}
\bpages{187--212}.
\bid{doi={10.1214/aop/1022855416}, issn={0091-1798}, mr={1617046}}
\end{barticle}
%
\bptok{imsref}%
\endbibitem

\bibitem{spde_book}
\begin{bbook}[auto:parserefs-M02]
\bauthor{\bsnm{Dalang},~\bfnm{R.~C.}\binits{R.~C.}},
\bauthor{\bsnm{Khohsnevisan},~\bfnm{D.}\binits{D.}},
\bauthor{\bsnm{Mueller},~\bfnm{C.}\binits{C.}},
\bauthor{\bsnm{Nualart},~\bfnm{D.}\binits{D.}} \AND
\bauthor{\bsnm{Xiao},~\bfnm{Y.}\binits{Y.}}
(\byear{2006}).
\btitle{A Minicourse in Stochastic Partial Differential Equations}.
\bseries{Lecture Notes in Math.}
\bvolume{1962}.
\bpublisher{Springer},
\blocation{Berlin}.
\end{bbook}
%
\bptok{imsref}%
\endbibitem

\bibitem{dalang-mueller03}
\begin{barticle}[mr]
\bauthor{\bsnm{Dalang},~\bfnm{Robert~C.}\binits{R.~C.}} \AND
\bauthor{\bsnm{Mueller},~\bfnm{Carl}\binits{C.}}
(\byear{2003}).
\btitle{Some non-linear S.{P}.{D}.{E}.'s that are second order in time}.
\bjournal{Electron. J. Probab.}
\bvolume{8}
\bpages{21 pp. (electronic)}.
\bid{doi={10.1214/EJP.v8-123}, issn={1083-6489}, mr={1961163}}
\bptnote{check pages}%
\end{barticle}
%
\bptok{imsref}%
\endbibitem

\bibitem{dalang-mueller09}
\begin{barticle}[mr]
\bauthor{\bsnm{Dalang},~\bfnm{Robert~C.}\binits{R.~C.}} \AND
\bauthor{\bsnm{Mueller},~\bfnm{Carl}\binits{C.}}
(\byear{2009}).
\btitle{Intermittency properties in a hyperbolic {A}nderson problem}.
\bjournal{Ann. Inst. Henri Poincar\'e Probab. Stat.}
\bvolume{45}
\bpages{1150--1164}.
\bid{doi={10.1214/08-AIHP199}, issn={0246-0203}, mr={2572169}}
\end{barticle}
%
\bptok{imsref}%
\endbibitem

\bibitem{DMT08}
\begin{barticle}[mr]
\bauthor{\bsnm{Dalang},~\bfnm{Robert~C.}\binits{R.~C.}},
\bauthor{\bsnm{Mueller},~\bfnm{Carl}\binits{C.}} \AND
\bauthor{\bsnm{Tribe},~\bfnm{Roger}\binits{R.}}
(\byear{2008}).
\btitle{A {F}eynman--{K}ac-type formula for the deterministic and stochastic wave equations and other P.{D}.{E}.'s}.
\bjournal{Trans. Amer. Math. Soc.}
\bvolume{360}
\bpages{4681--4703}.
\bid{doi={10.1090/S0002-9947-08-04351-1}, issn={0002-9947}, mr={2403701}}
\end{barticle}
%
\bptok{imsref}%
\endbibitem

\bibitem{dalang-sanzsole08}
\begin{barticle}[mr]
\bauthor{\bsnm{Dalang},~\bfnm{Robert~C.}\binits{R.~C.}} \AND
\bauthor{\bsnm{Sanz-Sol{\'e}},~\bfnm{Marta}\binits{M.}}
(\byear{2009}).
\btitle{H\"older--{S}obolev regularity of the solution to the stochastic wave equation in dimension three}.
\bjournal{Mem. Amer. Math. Soc.}
\bvolume{199}
\bpages{vi+70}.
\bid{doi={10.1090/memo/0931}, issn={0065-9266}, mr={2512755}}
\bptnote{check pages, check year}%
\end{barticle}
%
\bptok{imsref}%
\endbibitem

\bibitem{FK09}
\begin{barticle}[mr]
\bauthor{\bsnm{Foondun},~\bfnm{Mohammud}\binits{M.}} \AND
\bauthor{\bsnm{Khoshnevisan},~\bfnm{Davar}\binits{D.}}
(\byear{2009}).
\btitle{Intermittence and nonlinear parabolic stochastic partial differential equations}.
\bjournal{Electron. J. Probab.}
\bvolume{14}
\bpages{548--568}.
\bid{doi={10.1214/EJP.v14-614}, issn={1083-6489}, mr={2480553}}
\bptnote{check pages}%
\end{barticle}
%
\bptok{imsref}%
\endbibitem

\bibitem{FK13}
\begin{barticle}[mr]
\bauthor{\bsnm{Foondun},~\bfnm{Mohammud}\binits{M.}} \AND
\bauthor{\bsnm{Khoshnevisan},~\bfnm{Davar}\binits{D.}}
(\byear{2013}).
\btitle{On the stochastic heat equation with spatially-colored random forcing}.
\bjournal{Trans. Amer. Math. Soc.}
\bvolume{365}
\bpages{409--458}.
\bid{doi={10.1090/S0002-9947-2012-05616-9}, issn={0002-9947}, mr={2984063}}
\end{barticle}
%
\bptok{imsref}%
\endbibitem

\bibitem{gerhold12}
\begin{barticle}[mr]
\bauthor{\bsnm{Gerhold},~\bfnm{Stefan}\binits{S.}}
(\byear{2012}).
\btitle{Asymptotics for a variant of the {M}ittag--{L}effler function}.
\bjournal{Integral Transforms Spec. Funct.}
\bvolume{23}
\bpages{397--403}.
\bid{doi={10.1080/10652469.2011.596151}, issn={1065-2469}, mr={2929183}}
\end{barticle}
%
\bptok{imsref}%
\endbibitem

\bibitem{hairer}
\begin{barticle}[mr]
\bauthor{\bsnm{Hairer},~\bfnm{Martin}\binits{M.}}
(\byear{2013}).
\btitle{Solving the KPZ equation}.
\bjournal{Ann. of Math. (2)}
\bvolume{178}
\bpages{559--664}.
\bid{doi={10.4007/annals.2013.178.2.4}, issn={0003-486X}, mr={3071506}}
\end{barticle}
%
\bptok{imsref}%
\endbibitem

\bibitem{hu01}
\begin{barticle}[mr]
\bauthor{\bsnm{Hu},~\bfnm{Y.}\binits{Y.}}
(\byear{2001}).
\btitle{Heat equations with fractional white noise potentials}.
\bjournal{Appl. Math. Optim.}
\bvolume{43}
\bpages{221--243}.
\bid{doi={10.1007/s00245-001-0001-2}, issn={0095-4616}, mr={1885698}}
\end{barticle}
%
\bptok{imsref}%
\endbibitem

\bibitem{HHNT14}
\begin{bmisc}[auto:parserefs-M02]
\bauthor{\bsnm{Hu},~\bfnm{Y.}\binits{Y.}},
\bauthor{\bsnm{Huang},~\bfnm{J.}\binits{J.}},
\bauthor{\bsnm{Nualart},~\bfnm{D.}\binits{D.}} \AND
\bauthor{\bsnm{Tindel},~\bfnm{S.}\binits{S.}}
(\byear{2014}).
\bhowpublished{Stochastic heat equations with general multiplicative Gaussian noises: H\"older continuity and intermittency.
Preprint. Available at \arxivurl{arXiv:1402.2618}}.
\end{bmisc}
%
\bptok{imsref}%
\endbibitem

\bibitem{HLN12}
\begin{barticle}[mr]
\bauthor{\bsnm{Hu},~\bfnm{Yaozhong}\binits{Y.}},
\bauthor{\bsnm{Lu},~\bfnm{Fei}\binits{F.}} \AND
\bauthor{\bsnm{Nualart},~\bfnm{David}\binits{D.}}
(\byear{2012}).
\btitle{Feynman--{K}ac formula for the heat equation driven by fractional noise with {H}urst parameter {$H<1/2$}}.
\bjournal{Ann. Probab.}
\bvolume{40}
\bpages{1041--1068}.
\bid{doi={10.1214/11-AOP649}, issn={0091-1798}, mr={2962086}}
\end{barticle}
%
\bptok{imsref}%
\endbibitem

\bibitem{hu_nualart}
\begin{barticle}[mr]
\bauthor{\bsnm{Hu},~\bfnm{Yaozhong}\binits{Y.}} \AND
\bauthor{\bsnm{Nualart},~\bfnm{David}\binits{D.}}
(\byear{2009}).
\btitle{Stochastic heat equation driven by fractional noise and local time}.
\bjournal{Probab. Theory Related Fields}
\bvolume{143}
\bpages{285--328}.
\bid{doi={10.1007/s00440-007-0127-5}, issn={0178-8051}, mr={2449130}}
\end{barticle}
%
\bptok{imsref}%
\endbibitem

\bibitem{HNS11}
\begin{barticle}[mr]
\bauthor{\bsnm{Hu},~\bfnm{Yaozhong}\binits{Y.}},
\bauthor{\bsnm{Nualart},~\bfnm{David}\binits{D.}} \AND
\bauthor{\bsnm{Song},~\bfnm{Jian}\binits{J.}}
(\byear{2011}).
\btitle{Feynman--{K}ac formula for heat equation driven by fractional white noise}.
\bjournal{Ann. Probab.}
\bvolume{39}
\bpages{291--326}.
\bid{doi={10.1214/10-AOP547}, issn={0091-1798}, mr={2778803}}
\end{barticle}
%
\bptok{imsref}%
\endbibitem

\bibitem{kallenberg83}
\begin{bbook}[mr]
\bauthor{\bsnm{Kallenberg},~\bfnm{Olav}\binits{O.}}
(\byear{1983}).
\btitle{Random Measures},
\bedition{3rd} ed.
\bpublisher{Academic Press},
\blocation{London}.
\bid{mr={0818219}}
\end{bbook}
%
\bptok{imsref}%
\endbibitem

\bibitem{KPZ}
\begin{barticle}[auto:parserefs-M02]
\bauthor{\bsnm{Kardar},~\bfnm{M.}\binits{M.}},
\bauthor{\bsnm{Parisi},~\bfnm{G.}\binits{G.}} \AND
\bauthor{\bsnm{Zhang},~\bfnm{Y.-C.}\binits{Y.-C.}}
(\byear{1986}).
\btitle{Dynamic scaling of growing interfaces}.
\bjournal{Phys. Rev. Lett.}
\bvolume{56}
\bpages{889--892}.
\end{barticle}
%
\bptok{imsref}%
\endbibitem

\bibitem{li06}
\begin{barticle}[mr]
\bauthor{\bsnm{Li},~\bfnm{Yuan-Chuan}\binits{Y.-C.}}
(\byear{2006}).
\btitle{A note on an identity of the gamma function and {S}tirling's formula}.
\bjournal{Real Anal. Exchange}
\bvolume{32}
\bpages{267--271}.
\bid{issn={0147-1937}, mr={2329236}}
\bptnote{check pages, check year}%
\end{barticle}
%
\bptok{imsref}%
\endbibitem

\bibitem{MMV01}
\begin{barticle}[mr]
\bauthor{\bsnm{Memin},~\bfnm{J.}\binits{J.}},
\bauthor{\bsnm{Mishura},~\bfnm{Yuliya}\binits{Y.}} \AND
\bauthor{\bsnm{Valkeila},~\bfnm{Esko}\binits{E.}}
(\byear{2001}).
\btitle{Inequalities for the moments of Wiener integrals with respect to
fractional Brownian motions}.
\bjournal{Statist. Probab. Lett.}
\bvolume{55}
\bpages{421--430}.
\end{barticle}
%
\bptok{imsref}%
\endbibitem

\bibitem{millet-sanzsole99}
\begin{barticle}[mr]
\bauthor{\bsnm{Millet},~\bfnm{Annie}\binits{A.}} \AND
\bauthor{\bsnm{Sanz-Sol{\'e}},~\bfnm{Marta}\binits{M.}}
(\byear{1999}).
\btitle{A stochastic wave equation in two space dimension: Smoothness of the law}.
\bjournal{Ann. Probab.}
\bvolume{27}
\bpages{803--844}.
\bid{doi={10.1214/aop/1022677387}, issn={0091-1798}, mr={1698971}}
\end{barticle}
%
\bptok{imsref}%
\endbibitem

\bibitem{nualart98}
\begin{bincollection}[mr]
\bauthor{\bsnm{Nualart},~\bfnm{David}\binits{D.}}
(\byear{1998}).
\btitle{Analysis on {W}iener space and anticipating stochastic calculus}.
In \bbooktitle{Lectures on Probability Theory and Statistics ({S}aint-{F}lour, 1995)}.
\bseries{Lecture Notes in Math.}
\bvolume{1690}
\bpages{123--227}.
\bpublisher{Springer},
\blocation{Berlin}.
\bid{doi={10.1007/BFb0092538}, mr={1668111}}
\bptnote{check pages}%
\end{bincollection}
%
\bptok{imsref}%
\endbibitem

\bibitem{nualart06}
\begin{bbook}[mr]
\bauthor{\bsnm{Nualart},~\bfnm{David}\binits{D.}}
(\byear{2006}).
\btitle{The {M}alliavin Calculus and Related Topics},
\bedition{2nd} ed.
\bpublisher{Springer},
\blocation{Berlin}.
\bid{mr={2200233}}
\end{bbook}
%
\bptok{imsref}%
\endbibitem

\bibitem{Nualart-QS07}
\begin{barticle}[mr]
\bauthor{\bsnm{Nualart},~\bfnm{David}\binits{D.}} \AND
\bauthor{\bsnm{Quer-Sardanyons},~\bfnm{Llu{\'{\i}}s}\binits{L.}}
(\byear{2007}).
\btitle{Existence and smoothness of the density for spatially homogeneous {SPDE}s}.
\bjournal{Potential Anal.}
\bvolume{27}
\bpages{281--299}.
\bid{doi={10.1007/s11118-007-9055-3}, issn={0926-2601}, mr={2336301}}
\end{barticle}
%
\bptok{imsref}%
\endbibitem

\bibitem{QS-SanzSole04a}
\begin{barticle}[mr]
\bauthor{\bsnm{Quer-Sardanyons},~\bfnm{L.}\binits{L.}} \AND
\bauthor{\bsnm{Sanz-Sol{\'e}},~\bfnm{M.}\binits{M.}}
(\byear{2004}).
\btitle{Absolute continuity of the law of the solution to the 3-dimensional stochastic wave equation}.
\bjournal{J. Funct. Anal.}
\bvolume{206}
\bpages{1--32}.
\bid{doi={10.1016/S0022-1236(03)00065-X}, issn={0022-1236}, mr={2024344}}
\end{barticle}
%
\bptok{imsref}%
\endbibitem

\bibitem{QT07}
\begin{barticle}[mr]
\bauthor{\bsnm{Quer-Sardanyons},~\bfnm{Llu{\'{\i}}s}\binits{L.}} \AND
\bauthor{\bsnm{Tindel},~\bfnm{Samy}\binits{S.}}
(\byear{2007}).
\btitle{The 1-d stochastic wave equation driven by a fractional {B}rownian sheet}.
\bjournal{Stochastic Process. Appl.}
\bvolume{117}
\bpages{1448--1472}.
\bid{doi={10.1016/j.spa.2007.01.009}, issn={0304-4149}, mr={2353035}}
\end{barticle}
%
\bptok{imsref}%
\endbibitem

\bibitem{resnick07}
\begin{bbook}[mr]
\bauthor{\bsnm{Resnick},~\bfnm{Sidney~I.}\binits{S.~I.}}
(\byear{2007}).
\btitle{Heavy-Tail Phenomena: Probabilistic and Statistical Modeling}.
\bpublisher{Springer},
\blocation{New York}.
\bid{mr={2271424}}
\end{bbook}
%
\bptok{imsref}%
\endbibitem

\bibitem{sanzsole-sarra02}
\begin{bincollection}[mr]
\bauthor{\bsnm{Sanz-Sol{\'e}},~\bfnm{M.}\binits{M.}} \AND
\bauthor{\bsnm{Sarr{\`a}},~\bfnm{M.}\binits{M.}}
(\byear{2002}).
\btitle{H\"older continuity for the stochastic heat equation with spatially correlated noise}.
In \bbooktitle{Seminar on {S}tochastic {A}nalysis, {R}andom {F}ields and {A}pplications, III ({A}scona, 1999)}.
\bseries{Progress in Probability}
\bvolume{52}
\bpages{259--268}.
\bpublisher{Birkh\"auser},
\blocation{Basel}.
\bid{mr={1958822}}
\end{bincollection}
%
\bptok{imsref}%
\endbibitem

\bibitem{song12}
\begin{barticle}[mr]
\bauthor{\bsnm{Song},~\bfnm{Jian}\binits{J.}}
(\byear{2012}).
\btitle{Asymptotic behavior of the solution of heat equation driven by fractional white noise}.
\bjournal{Statist. Probab. Lett.}
\bvolume{82}
\bpages{614--620}.
\bid{doi={10.1016/j.spl.2011.11.017}, issn={0167-7152}, mr={2887479}}
\end{barticle}
%
\bptok{imsref}%
\endbibitem

\bibitem{stein70}
\begin{bbook}[mr]
\bauthor{\bsnm{Stein},~\bfnm{Elias~M.}\binits{E.~M.}}
(\byear{1970}).
\btitle{Singular Integrals and Differentiability Properties of Functions}.
\bseries{Princeton Mathematical Series}
\bvolume{30}.
\bpublisher{Princeton Univ. Press},
\blocation{Princeton, NJ}.
\bid{mr={0290095}}
\end{bbook}
%
\bptok{imsref}%
\endbibitem

\end{thebibliography}
\end{document}